%This is a corrected version of the paper "Discrete series characters..."
%(by Goresky, Kottwitz and MacPherson). We followed all the referee's 
%suggestions. Let me know if you need a copy of this file on a diskette.
%I'm also mailing you a hardcopy of this version.

\input amstex 
\documentstyle{amsppt}
\NoRunningHeads
\magnification=1200 
%\hcorrection{1truein} \vcorrection{1truein} 
\parskip
6pt 
\NoBlackBoxes \pagewidth{5.3in} \hfuzz=5pt
\refstyle{A}
\widestnumber\key{LS2}

\topmatter

\title Discrete series characters and the Lefschetz formula for Hecke operators \endtitle

\author M. Goresky, R. Kottwitz and R. MacPherson  \endauthor

\thanks Partially supported by NSF Grants DMS-9303550, DMS-9203380,
DMS-9106522 respectively. \endthanks 

\endtopmatter

\redefine\C{{ \Bbb C }} \redefine\r{{ \Bbb R }} \define\Z{{ \Bbb Z }}
\define\a{{ \alpha }} \redefine\b{{ \beta }} \redefine\d{ \delta }
\define\e{{ \epsilon }} \redefine\g{\gamma} 
\redefine\D{{ \Delta }} \redefine\l{\lambda} 
   
\redefine\det{\operatorname{det}}
\define\({ \left( } \define\){ \right) } \define\[{ \left[ } \define\]{
\right] } \define\<{ \langle } \define\>{ \rangle } \redefine\^{ \wedge }
  \define\tr{\operatorname{tr}}

\define\Int{\operatorname{Int}}

\define\Hom{\operatorname{Hom}}
\define\Lie{\operatorname{Lie}}

\define\Q{\Bbb Q}

\define\Gal{ \operatorname{Gal} }

\define\A{\Bbb A}

\define\G{\Gamma}

\define\der{\operatorname{der}}
\define\ad{\operatorname{ad}}
\define\vl{\operatorname{vol}}

\define\Fix{\operatorname{Fix}}

\define\Cent{\operatorname{Cent}}

\define\sc{\operatorname{sc}}

\define\Ad{\operatorname{Ad}}

\define\im{\operatorname{im}}

\define\reg{\operatorname{reg}}

\define\Sbar{\overline{S}}
\define\fA{{\frak A}}
\define\spn{\operatorname{span}}
\define\inC{\overset \circ \to C}
\define\intF{\overset \circ \to F}
\define\Cbar{\overline{C}}
\define\Ebar{\overline{\bold E}}
\define\sgn{\operatorname{sgn}}
\define\inI{\overset \circ \to I}

\document

This paper consists of three independent but related parts. In the first
part (\S\S1--6) we give a combinatorial formula for the constants appearing
in the ``numerators'' of characters of stable discrete series
representations of real groups (see \S3) as well as an analogous formula
for individual discrete series representations (see \S6). Moreover we
give an explicit formula (Theorem 5.1) for certain stable
virtual characters on real groups; by Theorem 5.2 these include the stable
discrete series characters, and thus we recover the results of \S3 in a
more natural way. 

In the second part (\S7) we use the character formula given in Theorem 5.1
to rewrite the Lefschetz formula of \cite{GM} (for the local contribution
at a single fixed point component to the trace of a Hecke operator on
weighted cohomology) in the same spirit as that of Arthur's Lefschetz
formula \cite{A}: in terms of stable virtual characters on real groups (see
Theorem 7.14.B). We then sum the contributions of the various fixed point
components and show that, in the case of middle weighted cohomology, the
resulting global Lefschetz fixed point formula agrees with Arthur's
Lefschetz formula. This gives a topological proof of Arthur's formula. 

The third part of the paper (Appendices A and B) is purely
combinatorial. In Appendix A we develop the combinatorics of convex
polyhedral cones on which our results on characters of real groups are
based. The methods of Appendix A are also used in Appendix B to prove a
generalization of a combinatorial lemma of Langlands.

The formula for stable discrete series constants given in Theorem 3.1 is
redundant, since it follows easily from Theorems 5.1 and 5.2. Nevertheless
the proof of Theorem 3.1 is instructive and should probably not be skipped
by readers interested in the case of individual discrete series
constants. Theorem 3.2 is not redundant and in fact provides the link
between our results on stable discrete series constants and individual
discrete series constants (we return to this point later in the
introduction). Because of the redundancy built into the paper, the reader
who is mainly interested in the Lefschetz formula only needs to read
\S\S5,7 and a little bit of Appendix A.

Let $G$ be a connected reductive group over $\Q$ and let $A_G$ denote the
maximal $\Q$-split torus in the center of $G$. Let $K_G$ be a maximal
compact subgroup of $G(\r)$ and let $X_G$ denote the homogeneous space
$$
   G(\r)/(K_G \cdot A_G(\r)^0). $$
Let $K$ be a suitably small compact open subgroup of $G(\A_f)$. We denote
by $S_K$ the space
$$
    G(\Q)\backslash [(G(\A_f)/K) \times X_G]. $$
Let $E$ be an irreducible representation of the algebraic group $G$ on a
finite dimensional complex vector space. Then $E$ gives rise to a local
system $\bold E_K$ on $S_K$.

Let $P_0=M_0N_0$ be a minimal parabolic subgroup of $G$, with Levi
component $M_0$ and unipotent radical $N_0$. As usual by a standard
parabolic subgroup of $G$ we mean one that contains $P_0$. For any standard
parabolic subgroup $P$ we write $P=MN$ where $M$ is the unique Levi
component of $P$ containing $M_0$ and $N$ is the unipotent radical of $P$.
The reductive Borel-Serre compactification $\Sbar_K$ of $S_K$ is a
stratified space whose strata are indexed by the standard parabolic
subgroups of $G$. The stratum indexed by $G$ is the space $S_K$. The
stratum $S^P_K$ indexed by standard $P=MN$ is a finite union of spaces of
the same type as $S_K$, but for the group $M$ rather than $G$.

Let
$$
    j:S_K \hookrightarrow \Sbar_K $$
denote the inclusion. Consider the object $\bold R j_*\bold E_K$ in the
derived category of $\Sbar_K$. The restriction to the stratum $S^P_K$ of
the $i$-th cohomology sheaf of $\bold R j_* \bold E_K$ is the local system
on $S^P_K$ associated to the representation of $M$ on
$$
    H^i(\Lie(N),E). $$

For any standard $P=MN$ we write $\fA_M$ for the real vector space
$$
    X_*(A_M) \otimes_\Z \r $$ 
and $\fA_M^*$ for its dual. Let $\nu \in \fA_{M_0}^*$ and suppose that the
restriction of $\nu$ to $\fA_G$ coincides with the element of $X^*(A_G)$ by
which $A_G$ acts on $E$. For any standard parabolic subgroup $P=MN$ we
write $\nu_P \in \fA^*_M$ for the restriction of $\nu$ to the subspace
$\fA_M$ of $\fA_{M_0}$. Then $\nu$ determines a weight profile and hence a
weighted cohomology complex (see \cite{GHM}) $\Ebar_K$ on $\Sbar_K$ (an
object in the derived category of $\Sbar_K$). The restriction to $S^P_K$ of
the $i$-th cohomology sheaf of $\Ebar_K$ is the local system on $S^P_K$
associated to the representation of $M$ on
$$
   H^i(\Lie(N),E)_{\ge \nu_P}, $$
the subspace of  $H^i(\Lie(N),E)$ on which $A_M$ acts by weights $\ge
\nu_P$ (a weight $\mu \in X^*(A_M) \subset \fA_M^* $ is $\ge\nu_P$ if
$\mu-\nu_P$ takes non-negative values on the chamber in $\fA_M$ determined
by $P$).

The main result of \cite{GM} is an explicit version of the Lefschetz
formula for the alternating sum of the traces of the self-maps on
$$
  H^i(\Sbar_K,\Ebar_K) $$
induced by a Hecke correspondence. One of our goals in this paper (see
Theorem 7.14.B) is to rewrite the Lefschetz formula in \cite{GM} in terms
involving a stable virtual character $\Theta_\nu$ on the group $G(\r)$. If
$\nu$ is so positive that $\Ebar_K$ coincides with the extension by zero of
$\bold E_K$, then $\Theta_\nu$ is just 
the character of the contragredient $E^*$
of $E$ (see 7.17). There is a similar (but more complicated) statement in
case $\nu$ is sufficiently negative (see 7.18). 
In general $\Theta_\nu$ is given by
    $$
\Theta_\nu=\sum_P (-1)^{\dim(A_M/A_G)} i^G_M(\d_P^{-1/2}\otimes
(E^\nu_P)^*) $$ 
where $E^\nu_P$ is the following virtual finite dimensional representation
of $M$ 
$$
   \sum_i (-1)^i H^i(\Lie(N),E)_{\ge\nu_P}, $$
and $\d_P$ is the usual modulus character
$$
   \d_P(x)=|\det(x;\Lie(N))| $$
on $M(\r)$; the sum is taken over all standard parabolic subgroups $P=MN$
and $i^G_M(\cdot)$ denotes normalized 
parabolic induction from $M(\r)$ to $G(\r)$.

Since there are simple formulas for characters that are
parabolically induced from finite dimensional representations of Levi
subgroups, it is possible to determine the character $\Theta_\nu$
explicitly (Theorem 5.1). In fact Theorem 5.1 is just what is needed to
rewrite the Lefschetz formula of \cite{GM} in terms of $\Theta_\nu$ (the
numbers $L^\nu_Q(\g)$ that go into the definition of $L^\nu_M(\g)$ occur as
factors in the local Lefschetz numbers). If $G(\r)$ has a discrete series
and if $\nu$ is the ``upper middle'' weight profile $\nu_m$ of \S5, then (see
Theorem 5.2) $\Theta_\nu$ agrees on all relevant maximal tori with
$$
    (-1)^{q(G)} \sum_{\pi \in \Pi} \Theta_\pi $$
where $\Pi$ is the L-packet of discrete series representations of $G(\r)$
having the same infinitesimal and central characters as $E^*$ (and in fact
$\Theta_\nu$ is equal to this virtual character if $P_0$ remains minimal
over $\r$), and our formula essentially coincides with 
Arthur's formula for $L^2$-Lefschetz numbers of
Hecke operators \cite{A} (see the remarks at the end of 7.19 for a detailed
comparison with Arthur's formula). This provides evidence for the agreement
of middle weighted cohomology and $L^2$-cohomology in this degree of
generality (in the Hermitian symmetric case this agreement is known: middle
weighted cohomology agrees with intersection cohomology of the Baily-Borel
compactification \cite{GHM} and this in turn
agrees with $L^2$-cohomology \cite{L},\cite{SS}).

This concludes our discussion of the global results in this paper. However
it remains to summarize the results on real groups. The two theorems just
mentioned (Theorems 5.1 and 5.2) together give a simple formula for the
stable discrete series character
$$
    \sum_{\pi \in \Pi} \Theta_\pi $$
on any maximal torus in $G$ over $\r$. In particular we obtain a simple
formula (Theorem 3.1) for the constants $d(w)$ $(w \in W)$ appearing in
stable discrete series character formulas (actually we prove Theorem 3.1 by
a different method). Here $W$ is the Weyl group of a root system $R$ such
that $-1 \in W$, for which we have fixed a system of positive roots. The
formula for $d(w)$ is expressed as a sum over $W$, each term in the sum
being $1,-1$ or $0$, and bears no obvious relation to the formula of Herb
\cite{He} for $d(w)$ in terms of two-systems. 
The terms in the sum depend on the finite dimensional representation $E$,
although their sum does not, so that we in fact get finitely many different
formulas, one for each cone in a certain decomposition of the positive Weyl
chamber. 
In Theorem 3.2 we prove an
unexpected symmetry for the function $d(w)$:
$$
    d(w^{-1})=(-1)^{q(R)}\epsilon(w) d(w), $$
where $\epsilon$ is the usual sign function on the Weyl group. This
symmetry, together with the formula for $d(w^{-1})$ as a sum over $W$,
gives a second formula for $d(w)$ as a sum over $W$. Our fixed system of
positive roots determines a certain subgroup $W_c$ of $W$ (the ``compact''
Weyl group). If in the second formula for $d(w)$ we replace the sum over
$W$ by a sum over a coset of $W_c$ in $W$, the resulting expression turns
out to be a formula (see \S 6) for the constants appearing in the
characters of \it individual \rm discrete series representations.
Of course these constants were already known, implicitly by work of
Harish-Chandra and explicitly by work of Hirai (or by combining
formulas for the stable constants---Herb's or ours---with Shelstad's theory
of endoscopy); what is perhaps interesting
is the simplicity of our formula (again we
in fact get finitely many different formulas, all giving the same result). 

More should be said about Theorem 5.2, which expresses stable discrete
series characters as linear combinations of characters induced from finite
dimensional representations of Levi subgroups. J.~Adams informs us that a
result of this kind was known to G.~Zuckerman when he wrote his 1974
Princeton thesis (certain examples are treated in the thesis, but a
precise general statement is not given there). We do not know if our result
is the one Zuckerman had in mind, though it seems unlikely that there could
be two essentially different formulas. What we do know is that an inversion
procedure due to Langlands (a simple special case of his combinatorial
lemma, often applied in Arthur's work on the trace formula) allows one to
obtain from Theorems 5.2 and 5.3 an expression for the character of a
finite dimensional representation as a linear combination of standard
characters, and it is not hard to see that this inverted formula coincides
with the one due to Zuckerman \cite{Z} (see also \cite{V}). Langlands's
inversion procedure works in both directions, so that one could also invert
Zuckerman's theorem to obtain our Theorem 5.2. However, this would result
in a more complicated proof (ours uses only elementary
combinatorics and Harish-Chandra's characterization of discrete series
characters). 

We wish to bring to the reader's attention some related results of
J.~Franke \cite{F}, G.~Harder \cite{H2}, A.~Nair \cite{N} and M.~Stern
\cite{St}. We would like to thank G.~Harder for useful discussions
concerning topological trace formulas and J.~Adams for his comments on the
history of Theorem 5.2. We are also indebted to the referee for finding a
mistake in our original formulation of Lemma 1.1(b). 
The first author would like to thank the
Institute for Advanced Study in Princeton for its hospitality and support
which was partially provided by NSF Grant DMS-9304580.

In this paper we use the following notation. For a finite set $S$ we write
$|S|$ for the cardinality of $S$. For a subset $A$ of a set $S$ we write
$\xi_A$ for the characteristic function of $A$. For a subgroup $H$ of a
group $G$ we write $N_G(H)$ (respectively, $\Cent_G(H)$) for the normalizer
(respectively, centralizer) of $H$ in $G$ (sometimes we allow $H$ to be a
subgroup of a bigger group of which $G$ is also a subgroup). For a free
abelian group $X$ of finite rank we write $X_\r$ for the real vector space
$X\otimes_\Z \r$. Given an endomorphism $A$ of a finite dimensional vector
space $V$, we write $\det(A;V)$ (respectively, $\tr(A;V)$) for the
determinant (respectively, trace) of $A$ whenever the name of $A$ leaves
doubt about which vector space $V$ we have in mind.
\heading 1. The function $\psi_R(C_0,x,\l)$ \endheading
In this section we consider a root system $(X,X^*,R,R^\vee)$. Here $X$ is a
real vector space, $X^*$ its dual vector space, $R \subset X^*$ a root
system in $X^*$ that spans $X^*$, and $R^\vee \subset X$ the coroot system
in $X$. We write $W=W(R)$ for the Weyl group of the root system $R$. For
any Weyl chamber $C$ in $X$ we write $C^\vee$ for the corresponding Weyl
chamber in $X^*$. In this section we will define integers
$\psi_R(C_0,x,\l)$; in \S3 we will see that in case $-1 \in W$ these
integers are (essentially) the ones appearing in the formulas for stable
discrete series characters on real groups. 

Let $C_0$ be a Weyl chamber in $X$. We write $\Cbar_0$ for the closure of
$C_0$. Let $\omega \in \Cbar_0$ be a non-zero element in 
 a 1-dimensional face of
$\Cbar_0$ (thus $\omega$ is, up to a positive scalar, a fundamental
coweight for $C_0$). Put 
$$R_\omega:=\{\a \in R \, | \, \a(\omega)=0 \}.$$
For any chamber $C$ in $X$ relative to $R$ let $\tilde C$ denote the unique
chamber in $X$ relative to $R_\omega$ that contains $C$. For two chambers
$C_1,C_2$ in $X$ relative to $R$ we write $l_R(C_1,C_2)$ or just
$l(C_1,C_2)$ for the number of root hyperplanes in $X$ separating $C_1$ and
$C_2$; thus $l(C_1,C_2)$ is the length with respect to $C_1$ of the unique
element $w \in W$ such that $wC_1=C_2$. We write $R^+$ for the set of roots
in $R$ that are positive on $C_0$. Finally we write $R^+_\omega$ for
$R_\omega \cap R^+$, the set of roots in $R_\omega$ that are positive on
$\tilde C_0$.
\proclaim{Lemma 1.1}\rom{(a)} The map $C \mapsto \tilde C$ yields a bijection
from the set of chambers in $X$ relative to $R$ whose closures contain
$\omega$ to the set of chambers in $X$ relative to $R_\omega$. If $C_1,C_2$
are chambers in $X$ relative to $R$ that contain $\omega$, then 
$$
l(C_1,C_2)=l_{R_\omega}(\tilde C_1,\tilde C_2). $$
\rom{(b)} Consider the difference $|R^+|-|R^+_\omega|$. If $\r\omega$ 
contains a
coroot, this difference is odd. If $\r\omega$ does not contain a coroot and
if $-1_{X/\r\omega} \in W(R_\omega)$,
then this difference is even. \newline
\rom{(c)} There exists a unique chamber $C'_0$ in $X$ such that $-\omega \in
\Cbar'_0$ and $\tilde C_0=\tilde C'_0$. Moreover 
$$
\{ \a \in R^+ \, | \, \text{$\ker(\a)$ separates
$C_0,C'_0$}\}=R^+\setminus R^+_\omega. $$
In particular
$$
l(C_0,C'_0)=|R^+|-|R^+_\omega|. $$
\rom{(d)} Suppose that there exists a positive scalar $c$ such that $c\omega$ is
a coroot $\a^\vee \in R^\vee$. Note that the corresponding root $\a$
belongs to $R^+$. Suppose that $C''$ is a chamber in $X$ satisfying the
following three conditions: \roster \item $\a$ takes non-negative values on
$C''$, \item $\ker(\a)$ is a wall of $C''$, \item $\tilde C''=\tilde C_0$.
\endroster Then
$$
l(C_0,C'')=(|R^+|-|R^+_\omega|-1)/2 $$
\rom(note that by \rom{(b)}
the quantity on the right-hand side is an integer\rom).
\endproclaim
The assertion (a) is standard. Now we prove (b). First suppose that
$\r\omega$ contains a coroot $\a^\vee$ and let $w \in W$ denote the
reflection in $\a^\vee$. We consider the action of
$w$ on the set $R/\pm$ obtained from $R$ by taking the quotient
by the action of the group $\{\pm1\}$. Since $w^2=1$, we see that
$|R^+|$ has the same parity as the number of fixed points of $w$
on $R/\pm$. Let $\b \in R$. Then
$w\b=\pm \b$ if and only if $\b \in R_\omega$ or $\b^\vee \in \r\omega$
(of course these alternatives are mutually exclusive). Therefore the number
of fixed points of $w$ on $R/\pm$ is $|R^+_\omega|+1$, which shows that
$|R^+|-|R^+_\omega|$ is odd, as desired.

Now suppose that $-1_{X/\r\omega} \in W(R_\omega)$. In other words we are
supposing that there exists an element $w \in W$ of order 2 whose $+1$
eigenspace is $\r\omega$ and whose $-1$ eigenspace is the span of
$R^\vee_\omega$. Again we consider the action of $w$ on $R/\pm$. Let $\b
\in R$. Then $w\b=\pm \b$ if and only if $\b \in R_\omega$ or $\b^\vee \in
\r\omega$. Suppose further that $\r\omega$ contains no coroot. Then the
number of fixed points of $w$ on $R/\pm$ is $|R^+_\omega|$, which shows
that $|R^+|-|R^+_\omega|$ is even, as desired. 
	
Now we consider (c). The existence and uniqueness of $C'_0$ follow from the
first statement in (a), applied to both $\omega$ and $-\omega$. Next we
prove the second statement in (c). Let $\a \in R^+$. Suppose first that $\a$
(strictly speaking, $\ker(\a))$ separates $C_0,C'_0$. Since $\tilde
C_0=\tilde C'_0$ it follows that $\a \notin R^+_\omega$. Conversely,
suppose that $\a \notin R^+_\omega$. Then $\a(\omega)\ne0$, so that $\a$
strictly separates $\omega,-\omega$. Since $\omega \in \Cbar_0$ and
$-\omega \in \Cbar'_0$ it follows that $\a$ separates $C_0,C'_0$.

Finally we consider (d). We are interested in roots $\b \in R$ that
separate $C_0,C''$. Using \therosteritem1 and \therosteritem3, we see that
any such $\b$ belongs to $R_0:=R\setminus(R_\omega \cup \{\pm\a\})$. To
prove (d) we must show that exactly half of the elements of $R_0$ separate
$C_0,C''$. Let $s \in W$ be the reflection in the root $\a$. Then $s$
preserves both $R_\omega$ and $\{\pm\a\}$ and hence preserves $R_0$ as
well. Let $\b \in R_0$. We will be done if we can show that $\b$ separates
$C_0,C''$ if and only if $s\b$ does not separate $C_0,C''$. Let $C'_0$ be
as in (c) and note that $C'_0=sC_0$. Of course $C_0,C''$ are separated by
$\b$ if and only if $C'_0=sC_0,sC''$ are separated by $s\b$. By
\therosteritem2 $sC''$ and $C''$ are separated only by $\pm\a$; therefore
$sC'',C''$ are not separated by $s\b$. Moreover $s\b$ does separate
$C_0,C'_0$ (use (c)). Therefore $C_0,C''$ are separated by $\b$ if and only
if $C_0,C''$ are not separated by $s\b$, as we wished to show. The proof of
the lemma is now complete.

Now we prepare to define the main object of study in this section. For any
two chambers $C_1$ and $C_2$ in $X$ we write $\e(C_1,C_2)$ for
$(-1)^{l(C_1,C_2)}$. For any chamber $C$ in $X$ we define a function
$\psi_C$ on $X \times X^*$ as follows. Let $\a_1,\dots,\a_n$ ($n=\dim(X)$)
be the simple roots in $R$ relative to $C$, and let
$\omega_1,\dots,\omega_n \in X$ be the basis for $X$ dual to the basis
$\a_1,\dots,\a_n$ for $X^*$. Let $x \in X$, $\l \in X^*$ and write
$$ \align
&x=a_1\omega_1+\dots+a_n\omega_n \\
&\l=b_1\a_1+\dots+b_n\a_n. \endalign $$
Let $I=\{1,\dots,n\}$ and define two subsets $I_x,I_\l$ of $I$ by
$$ \align
&I_x=\{i \in I \, | \, a_i \ge 0 \} \\
&I_\l = \{ i \in I \, | \, b_i \ge 0 \}. \endalign $$
Then define $\psi_C(x,\l)$ by
$$
\psi_C(x,\l)=   \cases (-1)^{|I_\lambda|}
&\text{if $I_\lambda=I \setminus I_x$}, \\
0 &\text{otherwise}. \endcases $$
Note that this function $\psi_C$ coincides with the function denoted by
$\psi_{\Cbar}$ in Appendix A (see Lemma A.1).

Now let $C_0$ be a chamber in $X$. Define a function
$\psi(C_0,\cdot,\cdot)$ on $X \times X^*$ by
$$
    \psi(C_0,x,\l)=\sum_C \e(C_0,C)\psi_C(x,\l), $$
where $C$ runs through the set of chambers in $X$. 
When it is necessary to stress the root system $R$ we will write
$\psi_R(C_0,x,\l)$ instead. We have the following obvious property:
$$ \align
\psi(wC_0,x,\l)&=\e(w)\psi(C_0,x,\l) \tag{1.1}\\
&=\psi(C_0,wx,w\l) \endalign $$
for any $w \in W$, where $\e(w)$ denotes the sign of $w$.

As usual we say that an element $x \in X$ is \it regular \rm if it lies on
no root hyperplane. We say that an element $\l \in X^*$ is $R$-\it regular \rm
if it lies on no hyperplane of the form $\{\l \in X^* \, | \,
\l(\omega)=0\}$, where $\omega$ is a non-zero element of a 1-dimensional
face of some closed Weyl chamber in $X$. Of course this notion of
regularity in $X^*$ is in general different from the usual one, and we
refer to the connected components in the set of $R$-regular elements in
$X^*$ as $R$-\it chambers \rm in $X^*$ to avoid confusion with the usual
Weyl chambers in $X^*$.

Suppose that $\omega$ is a non-zero element in a 1-dimensional face of
$\Cbar_0$. We adopt the notation of Lemma 1.1 and the discussion preceding
it (\it e.g.\rm, $R_\omega$, $R^+$,$R^+_\omega$,$\tilde C$,$C'_0$). Let $Z$
denote the hyperplane $\{\l \in X^* \,|\, \l(\omega)=0 \}$ in $X^*$. We
have the root system $(X/\r\omega,Z,R_\omega,R^\vee_\omega)$. Note that
$\tilde C_0=C_0+\r\omega$ has the same image as $C_0$ in $X/\r\omega$; we
denote this chamber in $X/\r\omega$ by $C^\omega_0$. 
\proclaim{Lemma 1.2} Let $x$ be a regular element in  $X$. 
The function $\psi(C_0,x,\cdot)$ on $X^*$
is constant on $R$-chambers. Suppose that $\l,\l'$ are $R$-regular elements
of $X^*$ lying in adjacent $R$-chambers separated only by the hyperplane
$Z$, and suppose further that $\l(\omega)>0$, $\l'(\omega)<0$. Then
$$
\psi(C_0,x,\l)-\psi(C_0,x,\l')=\cases -2\psi_{R_\omega}(C^\omega_0,\tilde
x,\tilde \lambda) &\text{if $\Bbb R\omega$ contains a coroot,}\\0
&\text{otherwise.} \endcases $$
Here $\tilde x$ denotes the image of $x$ in $X/\r\omega$ and $\tilde \l$
denotes the unique point of $Z$ lying on the line segment joining $\l$ and
$\l'$. Moreover $\tilde x$ is regular relative to $R_\omega$ and 
$\tilde \l$ is $R_\omega$-regular.  \endproclaim
It is clear that $\psi(C_0,x,\cdot)$ is constant on $R$-chambers. The
statement regarding the regularity of $\tilde x$ and $\tilde \l$ is easy
and will be left to the reader. By Corollary A.3
$$
  \psi(C_0,x,\l)-\psi(C_0,x,\l') $$
is equal to the sum over all chambers $C$ such that $\Cbar$ contains
$\omega$ or $-\omega$ of terms
$$
\pm\e(C_0,C)\psi_{\tilde C}(\tilde x,\tilde \l), $$
where the sign is $-$ if $\Cbar$ contains $\omega$ and $+$ if $\Cbar$
contains $-\omega$. We have abused notation slightly by writing $\tilde C$
when we mean its image in $X/\r\omega$; since $\Cbar$ contains $\omega$ or
$-\omega$ this image coincides with the image of $C$ in $X/\r\omega$. By
Lemma 1.1(c), for each chamber $C$ such that $\Cbar$ contains $\omega$,
there exists a unique chamber $C'$ such that $\Cbar'$ contains $-\omega$
and $\tilde C=\tilde C'$. Combining the terms for $C,C'$, we get 
$$
-\e(C_0,C)\(1-\e(C,C')\)\psi_{\tilde C}(\tilde x,\tilde \l). $$
From  Lemma 1.1(c) we see that $\e(C,C')$ is $-1$ if
$|R^+|-|R^+_\omega|$ is odd and is 1 otherwise. In the latter case each of
the combined terms is 0 and so is their sum.  In the former
case the sum of the combined terms is
$$
-2\sum_C \e(C_0,C)\psi_{\tilde C}(\tilde x,\tilde \l), $$
where $C$ ranges through the set of chambers in $X$ containing $\omega$. It
follows from Lemma 1.1(a) that this expression coincides with
$-2\psi_{R_\omega}(C^\omega_0,\tilde x,\tilde \l)$. 

Thus we have shown that 
$$
\psi(C_0,x,\l)-\psi(C_0,x,\l')=\cases -2\psi_{R_\omega}(C^\omega_0,\tilde
x,\tilde \lambda) &\text{if $|R^+|-|R^+_\omega|$ is odd,}\\0
&\text{otherwise.} \endcases \tag{1.2}$$
From the equality (1.2) and Lemma 1.1(b) we see that Lemma 1.2 holds
whenever $-1_{X/\r\omega} \in W(R_\omega)$. Using only equality (1.2), we
will prove Corollary 1.3 below. But then the general case of Lemma 1.2 will
follow, since $\psi_{R_\omega}(C_0^\omega,\tilde x,\tilde \l)=0$ if
$-1_{X/\r\omega} \notin W(R_\omega)$ (by Corollary 1.3).

\proclaim{Corollary 1.3} Suppose that $-1_X \notin W$. Then
$\psi(C_0,x,\l)=0$ for all regular $x \in X$ and all $R$-regular $\l \in
X^*$. \endproclaim
We prove this by induction on $\dim(X)$. 
If $\dim(X)=0$, the statement is trivially true. Now assume that $\dim(X)
>0$. Fix a regular element $x \in X$.
There exists $R$-regular $\l_0 \in X^*$ such that $\l_0(x)>0$. By
Proposition A.5 $\psi(C_0,x,\l_0)=0$. Therefore, to prove the corollary it
would be enough to show that
$$
\psi (C_0,x,\l)-\psi(C_0,x,\l') $$
vanishes whenever $\l,\l'$ lie in adjacent $R$-chambers separated by the
hyperplane 
$$
   Z=\{ \l \in X^* \, | \, \l(\omega)=0 \} $$
determined by a non-zero element $\omega$ of some 
1-dimensional face of the closure
of some chamber in $X$; by (1.1) it is harmless to assume that $\omega \in
\Cbar_0$. By equality (1.2), Lemma 1.1(b) and our induction hypothesis,
 this difference does vanish unless $-1_{X/\r\omega} \in W(R_\omega)$ and 
$\r\omega$ contains some coroot $\a^\vee$. 
But the product of $-1_{X/\r\omega}$ and reflection in the coroot $\a^\vee$
is equal to $-1_X$, and we are assuming that $-1_X \notin W$. We conclude
that the difference always vanishes, as desired.

We need to introduce more notation. Let $P$ (respectively, $Q$) denote
the lattice of coweights in $X$ (respectively, the lattice in $X$ generated
by the coroots). For any chamber $C$ in $X$ we denote by $\d_C \in P$ the
half-sum of the coroots that are positive for $C$. Put
$$ \align
&\hat A_{\sc}=\Hom(P,\C^\times) \\
&\hat A_{\ad}=\Hom(Q,\C^\times). \endalign  $$
The inclusion $Q \subset P$ induces a surjection 
$$\hat A_{\sc} \to \hat A_{\ad} $$ 
of complex tori, whose kernel we denote by $Z^\vee$, so that we get an
exact sequence
$$
  1 \to Z^\vee \to \hat A_{\sc} \to \hat A_{\ad} \to 1. $$
There are natural $\C^\times$-valued pairings $\langle\cdot,\cdot\rangle$
between $P$ and $\hat A_{\sc}$ and between $Q$ and $\hat A_{\ad}$. Let $s \in
\hat A_{\sc}$ and suppose that $s^2 \in Z^\vee$. Define a root system $R_s$
by $$
 \align &R^\vee_s=\{\a^\vee \in R^\vee \, | \, \langle\a^\vee,s\rangle=1\}
\\ &R_s=\{\a \in R \, | \, \a^\vee \in R^\vee_s \}. \endalign $$
When we defined $\psi_R(C_0,\cdot,\cdot)$ we insisted that $R$ generate
$X^*$. Of course this was just a matter of convenience. In the general case
the intersection of all the root hyperplanes in $X$ is a linear subspace
$X_0$ in $X$. Defining $\psi(C_0,\cdot,\cdot)$ as before, we find from
(A.2) that $\psi(C_0,x,\l)$ is 0 unless $\l$ vanishes on $X_0$, in which
case
$$
\psi(C_0,x,\l)=(-1)^{\dim(X_0)}\psi(\tilde C_0,\tilde x,\l) $$
where $\tilde C_0$ (respectively, $\tilde x$) denotes the image of $C_0$
(respectively, $x$) in $X/X_0$; note that on the right-hand side $\l$ is
regarded as an element of $(X/X_0)^*$. These remarks allow us to consider
the function $\psi_{R_s}(\tilde C_0,x,\l)$ obtained from
$(X,X^*,R_s,R^\vee_s)$, where $\tilde C_0$ denotes the unique chamber for
$R_s$ in $X$ containing $C_0$. 
\proclaim{Lemma 1.4} For all regular $x \in X$ and all $\l \in X^*$ there
is an equality
$$
\sum_C \e(C_0,C)\langle \d_C-\d_{C_0},s\rangle
\psi_C(x,\l)=\psi_{R_s}(\tilde C_0,x,\l), $$ 
in which the sum runs over all chambers $C$ in $X$. 
In particular, if $s^2 \ne 1$,
then the left-hand side of this equality vanishes for regular $x$ and
$R$-regular $\l$. \endproclaim
Since $s^2 \in Z^\vee$, the image of $s^2$ in $\hat A_{\ad}$ is 1, and
therefore $\langle \a^\vee,s \rangle=\pm 1$ for every coroot $\a^\vee$. By the
definition of $R_s$ we have $\langle \a^\vee,s \rangle=1$ if $\a \in R_s$
and $\langle \a^\vee,s \rangle = -1$ if $ \a \notin R_s$. Since
$\d_C-\d_{C_0}$ is the sum of the coroots in $R^\vee$ that are positive on
$C$ and negative on $C_0$, we see that
$$ \e(C_0,C)\langle \d_C-\d_{C_0},s \rangle = \e_{R_s}(\tilde C_0,\tilde C)
$$
where $\tilde C$ denotes the unique chamber of $(X,R_s)$ containing $C$.
Therefore the left-hand side of the equality we are trying to prove is
equal to the sum over chambers $D$ for $(X,R_s)$ of $\e(\tilde C_0,D)$ times
$$
\sum _C \psi_C(x,\l) $$
where $C$ runs through the chambers for $(X,R)$ contained in $D$. By
Proposition A.4 the difference between 
$$
   \sum_C \psi_C(x,\l) $$
and
$$
 \psi_D(x,\l) $$
is a sum of terms of the form $\pm\psi_F(x,\l)$ where $F$ is a proper face
of the closure of some chamber $C$ of $(X,R)$ contained in $D$. Here
$\psi_F$ is as in Appendix A. But $\psi_F(x,\l)$ vanishes unless $x \in
\spn(F)$ (see (A.2)). Therefore for regular $x$ the left-hand side of the
equality we are trying to prove is 
$$
   \sum_D \e(\tilde C_0,D)\psi_D(x,\l) $$
which, by definition, is $\psi_{R_s}(\tilde C_0,x,\l)$.

It remains to prove the second statement of the lemma. If the rank of the
root system $R_s$ is smaller than that of $R$, then $\psi_{R_s}(\tilde
C_0,x,\l)$ is 0 unless $\l$ belongs to the proper linear subspace
$\spn(R_s)$ of $X^*$. But the left-hand side of the equality of the lemma is
constant on $R$-chambers in $X^*$; therefore it vanishes for regular $x$
and $R$-regular $\l$. If $s^2\ne 1$ and $R_s$ has the same rank as $R$,
then $-1_X \notin W(R_s)$ (since $-1_X$ sends $s$ to $s^{-1}$ while all
elements of $W(R_s)$ fix $s$), and therefore by Corollary 1.3
$\psi_{R_s}(\tilde C_0,x,\l)$ vanishes for all regular $x \in X$ (regular
for $R_s$) and all $R_s$-regular $\l \in X^*$. Any $x \in X$ that is
regular for $R$ is regular for $R_s$, and again using that the left-hand
side of the equality of the lemma is constant on $R$-chambers in $X^*$, we
see that it vanishes for regular $x$ and $R$-regular $\l$. This completes
the proof of the lemma.

There is another result of this kind. With $P,Q$ as before now put
$$
\align &A_{\sc}=Q\otimes \C^\times \\ &A_{\ad}=P \otimes \C^\times.
\endalign $$
The inclusion $ Q \subset P $ induces a surjection 
$$
A_{\sc} \to A_{\ad} $$
of complex tori, whose kernel we denote by $Z$, so that we get an exact
sequence 
$$
    1 \to Z \to A_{\sc} \to A_{\ad} \to 1. $$
There are natural $\C^\times$-valued pairings $\langle \cdot,\cdot \rangle$
between $Q^*$ and $A_{\sc}$ and between $P^*$ and $A_{\ad}$ ($P^*,Q^*$ are
the free abelian groups dual to $P,Q$ respectively). Note that $Q^*$ is the
lattice of weights in $X^*$ and that $P^*$ is the lattice in $X^*$
generated by the roots. For any chamber $C$ in $X$ we write $\rho_C \in
Q^*$ for the half-sum of the roots that are positive for $C$. 

Let $a \in A_{\sc}$ and suppose that $a^2 \in Z$. Define a root system
$R_a$ by 
$$
   R_a=\{ \a \in R \, | \, \langle \a,a \rangle=1 \}. $$
Let $\tilde C_0$ denote the unique chamber for $R_a$ in $X$ containing
$C_0$. 
\proclaim{Lemma 1.5}For all regular $x \in X$ and all $\l \in X^*$ there is
an equality
$$
\sum_C \e(C_0,C)\langle
\rho_C-\rho_{C_0},a\rangle\psi_C(x,\l)=\psi_{R_a}(\tilde C_0,x,\l), $$
in which the sum runs over all chambers $C$ in $X$.
In particular, if $a^2\ne 1$, then the left-hand side of this equality
vanishes for regular $x$ and $R$-regular $\l$. \endproclaim
The proof is essentially the same as that of Lemma 1.4.
\heading 2. The function $\psi_R(C_0,x,\l)$ in case $-1 \in W$ \endheading
We continue with $X,X^*,R,R^\vee,W$ as in \S1. We still assume (for
convenience) that $R$ generates $X^*$, and we now add the assumption that
$-1_X \in W$. Let $\a$ be a root and define a root system $R_\a$ by
$$ \align
&R_\a^\vee=\{\b^\vee \in R^\vee \, | \, \langle\a,\b^\vee \rangle=0 \} \\ 
&R_\a=\{\b \in R \, | \, \b^\vee \in R_\a^\vee \}. \endalign $$
Let $Y$ denote the hyperplane $\{x \in X \, | \, \a(x)=0 \}$; then
$R_\a^\vee \subset Y$. Let $s_\a$ be the reflection in the root $\a$. Since
$-1_X$ belongs to $W$, so does $-s_\a$. But $-s_\a$ fixes $\a$, hence
belongs to $W(R_\a)$. Since $-s_\a$ acts by $-1$ on $Y$, we conclude that
$-1_Y \in W(R_\a)$. Therefore $(Y,Y^*,R_\a,R_\a^\vee)$ satisfies the same
conditions as $(X,X^*,R,R^\vee)$: $R_\a$ generates $Y^*$ and $-1_Y \in
W(R_\a)$. Note that $\a^\vee$ lies in the kernel of every root for $R_\a$.
Therefore $\a^\vee$ is a non-zero element in some 1-dimensional face of
some chamber in $X$, and $\a^\vee$ can serve as the element $\omega$
considered in \S1. Note that $R_\a=R_\omega$ with $R_\omega$ as in \S1.

There are two notions of chamber in $Y$. Of course we have the usual Weyl
chambers $D$ in $Y$ coming from the root system $R_\a$; these are
determined by the hyperplanes $\b=0$ ($\b \in R_\a$). There is a larger set
of hyperplanes in $Y$, namely those of the form $\b=0$ ($\b \in R \setminus
\{\pm\a\}$), and we will refer to the connected components $E$ of the
complement of this larger set of hyperplanes as \it chambers in \rm $Y$ \it
relative to \rm $R$. 

Fix a chamber $C_0$ in $X$ having $Y$ as a wall. 
As in \S1 we write $\tilde C$ for the unique chamber
for $X$ relative to $R_\a=R_\omega$ that contains $C$. 
It is easy to see that the
map $C \mapsto \tilde C 
 \cap Y$ is a bijection from the set of chambers $C$ in $X$
having $Y$ as a wall and lying on the same side of $Y$ as $C_0$ to the set
of closed 
chambers in $Y$ relative to $R$; note that the closure of $\tilde C \cap Y$
is equal to $\Cbar \cap Y$. 
Using the chamber $C_0$, we obtain a
function $\psi(C_0,\cdot,\cdot)$ on $X \times X^*$ as in \S1.
\proclaim{Lemma 2.1} Fix $\l \in X^*$. The function $\psi(C_0,\cdot,\l)$ on
$X$ is constant on the chambers in $X$. Suppose that $x,x'$ are regular
elements lying in adjacent chambers separated only by the hyperplane $Y$,
and assume that $x,C_0$ lie on the same side of $Y$ \rom(so that $x',C_0$ lie
on opposite sides of $Y$\rom). Then
$$
 \psi(C_0,x,\l)-\psi(C_0,x',\l)=2\psi_{R_\a}(D_0,y,\l_Y) $$
where $\l_Y \in Y^*$ denotes the restriction of $\l$ to $Y$, $y\in Y$ is
the unique point of $Y$ lying on the line segment joining $x$ and $x'$, and
$D_0$ is the chamber $\tilde C_0 \cap Y$ 
for $(Y,R_\a)$.
Moreover $y$ is regular in $Y$, and if $\l$ is $R$-regular then 
$\l_Y$ is $R_\a$-regular in $Y^*$.
\endproclaim 
It is clear that $\psi(C_0,\cdot,\l)$ is constant on chambers in $X$. By
Lemma A.2
$$
 \psi(C_0,x,\l)-\psi(C_0,x',\l) $$
is equal to 
$$
2\sum_C \e(C_0,C) \psi_{\Cbar \cap Y}(y,\l_Y), $$
where $C$ runs over the chambers in $X$ having $Y$ as a wall and lying on
the same side of $Y$ as $C_0$ (and $x$). The factor 2 arises since we have
combined the contributions of $C$ and the unique chamber adjacent to $C$
across the wall $Y$.  We denote by $C^\#$
the unique chamber for $(X,R)$ that is contained in $\tilde C$ and whose
closure contains $\omega$. Replacing $\a$ by $-\a$ if necessary, we may
assume without loss of generality that $\a$ is non-negative on $C_0$ and
$x$, and hence on any $C$ appearing in the sum above. Applying Lemma 1.1(d)
to both $C_0$ and any such $C$ (both satisfy conditions \therosteritem1 and
\therosteritem2), we see that
$$
\e(C_0,C)=\e(C_0^\#,C^\#), $$
and then from Lemma 1.1(a) we see further that
$$ \align
\e(C_0,C)&=\e(\tilde C_0,\tilde C) \\
&=\e(D_0,D) \endalign $$
where $D_0=\tilde C_0 \cap Y$ and $D=\tilde C \cap Y$.

We have now shown that the left-hand side of the equality we are trying to
prove is equal to
$$
2 \sum_D \e(D_0,D) \sum_E \psi_{\overline{E}}(y,\l_Y) $$
where $D$ runs over the chambers of $(Y,R_\a)$ and $E$ runs over the
chambers in $Y$ relative to $R$ such that $E \subset D$ (the function
$\psi_{\overline{E}}$ is the one attached in Appendix A to the closed
convex polyhedral cone $\overline{E}$ in $Y$). Each such $D$ is the
disjoint union of the corresponding $E$'s together with the relative
interiors of some closed convex polyhedral cones $F$ of lower dimension,
each of which is contained in some root hyperplane other than $Y$. It is
clear that $y$ lies on no root hyperplane of $R$ other than $Y$; therefore
$\psi_F(y,\l_Y)$ vanishes for all such $F$ (see (A.2)). By Proposition A.4
the inner sum is equal to $\psi_D(y,\l_Y)$, and therefore the whole
expression is equal to $2\psi_{R_\a}(D_0,y,\l_Y)$. This proves the lemma,
except for the last statement, which we leave to the reader.

Again let $C_0$ be a chamber in $X$ and let $R^+$ be the set of roots in
$R$ that are positive for $C_0$. Let $\e:W \to \{\pm1\}$ be the sign
homomorphism. The longest element of $W$ is $-1_X$. On the one hand
$$
    \e(-1_X)=(-1)^{|R^+|}. $$
On the other hand
$$
\e(-1_X)=\det(-1_X)=(-1)^{\dim(X)}. $$
Therefore $|R^+|$ and $\dim(X)$ have the same parity, and we can define an
integer $q(R)$ by
$$
  q(R):=(|R^+|+\dim(X))/2. $$
To understand the significance of the integer $q(R)$ one should note that
it is half the dimension of the symmetric space of the split semisimple
real group $G$ with root system $R$ (a number that is traditionally denoted
$q(G)$). 

Let $C_0^\vee$ be the Weyl chamber in $X^*$ corresponding to $C_0$. From
$C_0^\vee$ and the coroot system $R^\vee$ we get a function
$\psi_{R^\vee}(C_0^\vee,\cdot,\cdot)$ on $X^* \times X$ (the roles of $X$
and $X^*$ are now reversed). We have the notions of regularity (for $x \in
X$) and $R$-regularity (for $\l \in X^*$) from before. Applying these
definitions to $R^\vee$ rather than $R$, we have the notions of regularity
for $\l \in X^*$ and $R^\vee$-regularity for $x \in X$; note that the set
of regular elements in $X^*$ is the union of the Weyl chambers in $X^*$.
Since $-1_X \in W$ every coroot in $X$ is a non-zero element in a
1-dimensional face of some chamber in $X$ (we saw this during the
discussion at the beginning of this section). Therefore, if $\l \in X^*$ is
$R$-regular, it is automatically regular, and, similarly, if $x \in X$ is
$R^\vee$-regular, it is automatically regular. 
\proclaim{Lemma 2.2} For any $R^\vee$-regular $x \in X$ and any $R$-regular
$\l \in X^*$ there is an equality
$$
\psi_R(C_0,x,\l)=(-1)^{q(R)}\psi_{R^\vee}(C_0^\vee,\l,x). $$
\endproclaim
We prove this by induction on $\dim(X)$. It is certainly true when
$\dim(X)=0$ (the empty root system). Now assume that $\dim(X)>0$. 
Fix an $R^\vee$-regular element $x \in
X$. There exists $R$-regular $\l_0 \in X^*$ such that $\l_0(x) > 0$, and the
equality in the lemma holds for $x,\l_0$ since both sides of the equality
vanish by Proposition A.5. Therefore it is enough to show that
$$
\psi_R(C_0,x,\l)-\psi_R(C_0,x,\l')=(-1)^{q(R)}(\psi_{R^\vee}(C_0^\vee,\l,x)-\psi_{R^\vee}(C_0^\vee,\l',x))
\tag{2.1} $$
whenever $\l,\l'$ are $R$-regular elements of $X^*$ lying in adjacent
$R$-chambers. Let $Z$ denote the unique hyperplane separating these two
adjacent $R$-chambers. Thus $Z$ is of the form
$$
  Z=\{\l \in X^* \, | \, \l(\omega) =0 \} $$
for some non-zero $\omega$ lying in a 1-dimensional face of the closure of
some Weyl chamber in $X$; we may assume without loss of generality that
this Weyl chamber coincides with $C_0$ ( by property (1.1) changing $C_0$
changes both sides of (2.1) by the same sign). By switching $\l,\l'$ if
necessary we may also assume that $\l(\omega) > 0$ and $\l'(\omega) < 0$. 

First consider the case in which $\r\omega$ does not contain a coroot. Then
the left-hand side of (2.1) vanishes by Lemma 1.2, while the right-hand
side vanishes because $\psi_{R^\vee}(C_0^\vee,\cdot,x)$ is constant on Weyl
chambers in $X^*$, not just on $R$-chambers. We are left with the case in
which $\r\omega$ contains a coroot $\a^\vee$; replacing $\omega$ by a
positive scalar multiple we may as well assume that $\omega=\a^\vee$. Since
$\a^\vee \in \Cbar_0$, the root $\a$ is positive for $C_0$. By Lemma 1.2
the left-hand side of (2.1) is equal to
$$
   -2\psi_{R_\omega}(C_0^\omega,\tilde x,\tilde \l) $$
(with notation as in that lemma).

Of course we are going to use Lemma 2.1 to evaluate the right-hand side of
(2.1). However the hyperplane
$$
Z=\{ \l \in X^* \, | \, \l(\omega)=0 \}=\{\l \in X^* \, | \, \l(\a^\vee)=0
\} $$
need not be a wall of $C_0^\vee$, so that we need to introduce another
chamber $(C^\vee)''$ in $X^*$, better suited to our purposes. We take
$(C^\vee)''$ to be any chamber in $X^*$ satisfying the following three
conditions \roster \item $\a^\vee$ takes non-negative values on
$(C^\vee)''$, \item $Z$ is a wall of $(C^\vee)''$, \item
$((C^\vee)'')^\sim=(C_0^\vee)^\sim$. \endroster
The notation $(\cdot)^\sim$ used in \therosteritem3 has the following
meaning: for any chamber $C^\vee$ in $X^*$ we write $(C^\vee)^\sim$ for the
unique chamber in $X^*$ relative to  $R_\omega^\vee=(R^\vee)_{\a^\vee}$
that contains $C^\vee$. It is easy to see that $(C^\vee)''$ exists: pick
any chamber $E^\vee$ in $Z$ relative to $R^\vee$ contained in
$(C^\vee_0)^\sim \cap Z $ and take for $(C^\vee)''$ the unique chamber in
$X^*$ satisfying \therosteritem1 and \therosteritem2 and having the
property that the closure of $E^\vee$ is equal to the intersection of $Z$
with the closure of $(C^\vee)''$. 

By (1.1) and Lemma 2.1 the right-hand side of (2.1) is equal to
$$
    2(-1)^{q(R)}\e(C_0^\vee,(C^\vee)'')\psi_{R_\omega^\vee}(D_0^\vee,\tilde
\l,\tilde x), $$
with $\tilde \l,\tilde x$ as before and $D_0^\vee$ the 
chamber $((C^\vee)'')^\sim 
\cap Z$ in $Z$ for the root system $R_\omega^\vee$. 
It follows from Lemma 1.1(d) that
$$
 \e(C_0^\vee,(C^\vee)'')=(-1)^{(|R^+|-|R^+_\omega|-1)/2}; $$
since 
$$
   q(R)-q(R_\omega)=(|R^+|-|R^+_\omega|+1)/2, $$
we conclude that
$$
(-1)^{q(R)}\e(C_0^\vee,(C^\vee)'')=-(-1)^{q(R_\omega)}. $$
By our induction hypothesis
$$
  (-1)^{q(R_\omega)}\psi_{R^\vee_\omega}(D_0^\vee,\tilde \l,\tilde
x)=\psi_{R_\omega}(C_0^\omega,\tilde x,\tilde \l). $$ Of course we used
that $D^\vee_0=(C^\omega_0)^\vee$ and that $\tilde x,\tilde \l$ are
suitably regular. Therefore the right-hand side of (2.1) equals
$$
   -2\psi_{R_\omega}(C_0^\omega,\tilde x,\tilde \l), $$
which coincides with the expression we found for the left-hand side. This
concludes the proof of the lemma.\heading 3. Stable discrete series constants $\bar c_R$ \endheading
In the theory of stable discrete series characters on real groups, which we
will review briefly in \S4, there appear integer-valued functions
(see \cite{K}, \cite{He}, \cite{He2})
$$
\bar c_R: X_{\reg}\times X^*_{\reg} \to \Z $$
for every root system $(X,X^*,R,R^\vee)$ satisfying the two conditions of \S2
($R$ generates $X^*$ and $-1_X \in W$, where $W=W(R)$ denotes the Weyl
group of 
$R$). Here $X_{\reg}$ and $X^*_{\reg}$ denote the sets of regular elements
in $X$ and $X^*$ respectively (regular in the usual sense, so that
$X_{\reg}$, $X^*_{\reg}$ can also be described as the unions of the Weyl
chambers in $X$,$X^*$ respectively). The functions $\bar c_R$ satisfy the
following five properties:
\roster \item $\bar c_R(0,0)=1$ if $R$ is empty, \item $\bar c_R(x,\l)$
depends only on the chamber in $X$ in which $x$ lies and the chamber in
$X^*$ in which $\l$ lies,
\item $\bar c_R(x,\l)=0$ unless $\l(x)\le 0$, \item if $x,x' \in X$ lie in
adjacent chambers, separated only by the root hyperplane $Y$, then
$$
   \bar c_R(x,\l)+\bar c_R(x',\l)=2\bar c_{R_Y}(y,\l_Y), $$
where $R_Y \subset Y^*$ is the root system whose set of coroots is $R^\vee
\cap Y$, $\l_Y$ is the restriction of $\l$ to $Y$, and $y$ is the unique
point of $Y$ lying on the line segment joining $x$ and $x'$, \item
$\bar c_R(wx,w\l)=\bar c_R(x,\l)$ for all $w \in W(R)$. \endroster

It is well-known  that the collection of functions $\bar c_R$ is
characterized uniquely by properties \therosteritem1, \therosteritem3,
\therosteritem4 (this follows easily from an induction on $\dim(X)$ as in
the proofs of Lemma 2.2 and Corollary 1.3). Of course these properties are
reminiscent of ones enjoyed by the functions $\psi_R(C_0,x,\l)$ studied in
\S2. For $x \in X_{\reg}$ denote by $C_x$ the unique chamber in $X$
containing $x$. We now define an integer-valued function $m_R$ on
$X_{\reg}\times X^*_{\reg}$ by
$$
m_R(x,\l)=\psi_R(C_x,x,\l). $$

\proclaim{Theorem 3.1} The functions $\bar c_R$ and $m_R$ are equal.
\endproclaim

We need only show that $m_R$ satisfies properties \therosteritem1,
\therosteritem3, \therosteritem4 above. Property \therosteritem1 is
trivial. Property \therosteritem3 follows from Proposition A.5. Property
\therosteritem4 follows from (1.1) and Lemma 2.1.

There is a more efficient way to encode the information in the function
$\bar c_R$. Fix a Weyl chamber $C_0$ in $X$ and let $C_0^\vee$ be the
corresponding Weyl chamber in $X^*$. Then define an integer-valued function
$d$ on $W=W(R)$ by putting
$$
   d(w):=\bar c_R(x_0,w\l_0) \qquad (w \in W) $$
where $x_0$ (respectively, $\l_0$) is any point in $C_0$ (respectively,
$C^\vee_0$). Of course $d$ depends (in a simple way) on the choice of
$C_0$. Applying this construction to the root system $R^\vee$ (and the
chamber $C_0^\vee$) we get a function $d^\vee$ on $W=W(R^\vee)=W(R)$. 

\proclaim{Theorem 3.2} For all $w \in W$ there are equalities \roster
\item $d^\vee(w)=d(w)$ \item $d(w^{-1})=(-1)^{q(R)}\e(w)d(w)$ \endroster
where $\e(w)$ denotes the sign of $w$. \endproclaim

Pick a $W$-equivariant isomorphism $j:X \to X^*$. Let $x_0 \in C_0$; then
$\l_0:=j(x_0) \in C_0^\vee$. Since $\psi_R$ depends only on the root
hyperplanes and not on the roots themselves, it is clear that
$$
\psi_{R^\vee}(C_0^\vee,\l_0,wx_0)=\psi_R(C_0,x_0,w\l_0), \tag{3.1} $$
and by Theorem 3.1 this just says that 
$$
  d^\vee(w)=d(w). $$
The equality (2) follows from Theorem 3.1 and Lemma 2.2 (use (1.1) and
(3.1) as well). \heading 4. Background material on stable characters \endheading
In this section we review some of the theory of characters of irreducible
representations of real groups. Let $G$ be a connected reductive group over
$\r$. Let $E$ be an irreducible finite dimensional complex representation
of the algebraic group $G$. We are interested in irreducible
representations $\pi$ of $G(\r)$ (irreducible Harish-Chandra modules)
having the same infinitesimal character as $E$. Harish-Chandra 
associated to any such $\pi$ its character $\Theta_\pi$, a real-analytic
function on $G_{\reg}(\r)$, the set of regular semisimple elements in
$G(\r)$. 

Let $T$ be a maximal torus in $G$. Let $\Cal B(T)$ denote the set of Borel
subgroups of $G$ over $\C$ containing $T$. Let $R$ be the set of roots of
$T$ in $G$. For $B \in \Cal B(T)$ denote by $\l_B \in X^*(T)$ the highest
weight of $E$ relative to $B$, denote by $\rho_B \in X^*(T)_\r$ half the
sum of the roots in $R$ that are positive for $B$ and denote by $\Delta_B$
the Weyl denominator
$$
\Delta_B=\prod_{\a>0} (1-\alpha^{-1}) $$
for $T$ relative to $B$ (the index set is the subset of $R$ consisting of
roots that are positive for $B$). 

The character of $E$ on $T_{\reg}(\r):=T(\r)\cap G_{\reg}(\r)$ is given by
$$
     \tr(\g;E)=\sum_{B \in \Cal B(T)} \l_B(\g)\cdot \Delta_B(\g)^{-1}
\qquad (\g \in T_{\reg}(\r)). $$
The character $\Theta_\pi$ of $\pi$ on $T_{\reg}(\r)$ is given by a similar
expression
$$
   \Theta_\pi(\g)=\sum_{B \in \Cal B(T)} n(\g,B) \l_B(\g) \D_B(\g)^{-1}
\tag{4.1} $$
for certain integers $n(\g,B)$ depending on $(\g,B)$. Of course the
invariance of $\Theta_\pi$ under conjugation by $G(\r)$ implies that
$$
    n(\g,B)=n(w\g w^{-1},wBw^{-1}) \text{ for all $w \in
\Omega(T(\r),G(\r))$}, \tag{4.2} $$ 
where $\Omega(T(\r),G(\r))$ denotes the real Weyl group 
$$
     N_{G(\r)}(T)/T(\r). $$
For $\g \in T_{\reg}(\r)$ define subsets $R_\g$ and $R^+_\g$ of $R$ by
$$ \align
    &R_\g:=\{\a \in R \,|\, \text{$\a$ is real and $\a(\g) > 0$}\} \\
    &R^+_\g:=\{\a \in R \,|\, \text{$\a$ is real and $\a(\g) > 1$}\}.
\endalign $$
Note that $R_\g$ is a root system and that $R^+_\g$ is a positive system in
$R_\g$. Moreover $R_\g$ depends only on the connected component $\G$ of
$T(\r)$ in which $\g$ lies; thus we sometimes write $R_{\G}$ instead of
$R_\g$.  Harish-Chandra \cite{HC, Lemma 25} showed that 
$$
   n(\g_1,B)=n(\g_2,B) \text{ if $\Gamma_1 = \Gamma_2$ 
and $R^+_{\g_1}=R^+_{\g_2}$}, \tag{4.3} $$
where $\Gamma_i$ denotes the connected component of $T(\r)$ in which $\g_i$
lies. 

Of course any finite $\Z$-linear combination $\Theta$ of characters
$\Theta_\pi$ as above can also be expressed in the form (4.1) for integers
$n(\g,B)$ satisfying (4.2) and (4.3) (we refer to $\Theta$ as a virtual
character on $G(\r)$). We are particularly interested in virtual characters
$\Theta$ on $G(\r)$ that are \it stable \rm in the sense that
$$
     \Theta(\g)=\Theta(\g') $$
whenever $\g,\g' \in G_{\reg}(\r)$ are stably conjugate. A virtual character
$\Theta$ is stable if and only if the integers $n(\g,B)$ satisfy the
following strengthening of (4.2) (for all $T$):
$$
    n(\g,B)=n(w\g w^{-1},wBw^{-1})  \text{ for all $w \in W(\r)$,} \tag{4.4} $$
where $W$ is the Weyl group of $T_\C$ in $G_\C$ and $W(\r)$ is the subgroup
of $W$ consisting of all elements that are fixed by complex conjugation (of
course $W(\r)$ contains $\Omega(T(\r),G(\r)$).

Let $A$ be the maximal split subtorus of $T$ and let $M$ be the centralizer
of $A$ in $G$, a Levi subgroup of $G$. As usual for $\g \in M(\r)$ we
define a real number $D^G_M(\g)$ by 
$$
   D^G_M(\g)=\det(1-\Ad(\g);\Lie(G)/\Lie(M)). $$
We will need the following result of Arthur \cite{A} and Shelstad.
\proclaim{Lemma 4.1} For any stable virtual character $\Theta$ on $G(\r)$
the function 
$$
    \g \mapsto |D^G_M(\g)|^{1/2}\Theta(\g) $$
on $T_{\reg}(\r)$ extends continuously to $T(\r)$. \endproclaim

Let $\G$ be a connected component of $T(\r)$  and let $\G_{\reg}$ denote its
intersection with $T_{\reg}(\r)$. To prove the lemma we must show that
$$
    |D^G_M(\g)|^{1/2} \Theta(\g) $$
extends continuously from $\G_{\reg}$ to $\G$. Pick an element $a \in \G$
such that $a^2=1$ (it is easy to see that such an element exists). The root
system $R_\Gamma$ defined above is equal to the set of real roots $\a \in R$
such that $\a(a)=1$; thus $a$ lies in the center of the connected
reductive subgroup of $G$ containing $T$ with root system $R_\G$, and we
conclude that $a$ is fixed by the Weyl group $W(R_\G)$ of $R_\G$. Thus $\G$
is fixed by the subgroup $W(R_\G)$ of $\Omega(T(\r),G(\r))$, and since both
$|D^G_M(\g)|^{1/2}$, $\Theta(\g)$ are invariant under $\Omega(T(\r),G(\r))$
(which normalizes $M$), it follows that the function
$$
   |D^G_M(\g)|^{1/2}\Theta(\g) $$
on $\G_{\reg}$ is invariant under $W(R_\G)$. Let $T_c$ denote the maximal
anisotropic subtorus of $T$. Then 
$$
\G=a\cdot T_c(\r) \cdot \exp(\fA), $$
where
$$
   \fA=X_*(A)_\r=\Lie(A(\r)). $$
The Weyl group $W(R_\G)$ fixes $T_c(\r)$ as well as $a$. Fix a positive
system $R^+_\G$ in $R_\G$ and let $\Cbar$ be the corresponding closed
chamber in $\fA$. Then $\Cbar$ is a closed fundamental domain for the action
of $W(R_\G)$ on $\fA$, and therefore a $W(R_\G)$-invariant function on $\G$
is continuous if and only if its restriction to
$$
   \G^+:=a \cdot T_c(\r) \cdot \exp(\Cbar) $$
is continuous. Therefore it is enough to show that
$$
   |D^G_M(\g)|^{1/2}|\Theta(\g) $$
extends continuously to $\G^+$. For any regular element $\g \in \G^+$ we
have $R^+_\g=R^+_\G$, and thus there are integers $m(B)$ ($B \in \Cal
B(T)$) such that for all regular $\g \in \G^+$
$$
  \Theta(\g)=\sum_{B \in \Cal B(T)} m(B) \l_B(\g) \D_B(\g)^{-1}. $$

The Weyl group $W_M$ of $T_\C$ in $M_\C$ is a subgroup of $W(\r)$, and this
subgroup fixes $A$ pointwise and hence preserves $\G$ and $R^+_\G$.
Therefore it follows from (4.4) that 
$$
     m(B)=m(wBw^{-1}) \text{ for all $w \in W_M$.} \tag{4.5} $$
Choose a parabolic subgroup $P=MN$ having $M$ as Levi component
and having the property that every element of $R_\G$ that appears in
$\Lie(N)$ is non-negative on $\Cbar$ (here $N$ denotes the unipotent
radical of $P$). Put 
$$
  \D_P=\prod_\a (1-\a^{-1}) $$
where $\a$ runs through the roots of $T$ in $\Lie(N)$. We claim that
$\D_P(\g)$ is non-negative for all $\g \in \G^+$. Indeed, complex
conjugation preserves the set of roots of $T$ in $\Lie(N)$. If the complex
conjugate $\bar \a$ is different from $\a$, then the contribution of
$\a,\bar \a$ to $\D_P$ is $(1-\a(\g)^{-1})$ times its complex conjugate;
this contribution is certainly non-negative. If $\a$ is real, then
$\a(a)=\pm 1$. If $\a(a)=-1$, then $\a(\g)^{-1}$ is negative and therefore
$1-\a(\g)^{-1}$ is positive. If $\a(a)=1$, then $\a \in R_\G$ and by our
choice of $P$ we have $\a(\g)\ge 1$, so that $1-\a(\g)^{-1}\ge 0$. It
follows from the claim that
$$
   |D^G_M(\g)|^{1/2}=\D_P(\g) \cdot \d^{1/2}_P(\g), $$
where $\d_P$ denotes the modulus character
$$
    \d_P(x):=|\det(x;\Lie(N))| $$
on $M(\r)$. Therefore it is enough to show that
$$
    \sum_{B \in \Cal B(T)} m(B) \cdot \D_P(\g) \cdot \l_B(\g) \cdot
\D_B(\g)^{-1} $$ extends continuously to $\G^+$. But it follows immediately
from (4.5) that this last expression is a linear combination of characters
of irreducible finite dimensional representations of $M$, and of course
such a linear combination extends continuously to $\G^+$ (and even to 
all of $T(\r)$). This completes the proof of the lemma. 

For any stable virtual character $\Theta$ on $G(\r)$ we denote by
$\Phi^G_M(\g,\Theta)$ the (unique) continuous extension of
$$
   |D^G_M(\g)|^{1/2}\Theta(\g) $$
to $T(\r)$ whose existence is asserted in the lemma we just proved.
Sometimes it is convenient to extend $\Phi^G_M(\g,\Theta)$ to a function on
the set of all elliptic elements in $M(\r)$ (in other words, the set of
$M(\r)$-conjugates of elements in $T(\r)$) by taking the unique extension
that is invariant under conjugation by $M(\r)$.

The functions $\Phi^G_M(\cdot,\Theta)$ behave simply under induction. Let
$Q=LU$ be a parabolic subgroup of $G$ with Levi subgroup $L$ and unipotent
radical $U$. Let $\Theta_L$ be a stable virtual character on $L(\r)$ and
let $\Theta=i^G_L(\Theta_L)$ be the virtual character on $G(\r)$ obtained
from $\Theta_L$ by the usual normalized parabolic induction. Let $T \subset
M \subset G$ be as above. Then $\Theta$ is stable and for all $\g \in
T(\r)$ 
$$
   \Phi^G_M(\g,\Theta)=\sum_{gL(\r)}
\Phi^{gLg^{-1}}_M(\g,\Theta_{gLg^{-1}}), \tag{4.6} $$ 
where the sum runs over the set of cosets $gL(\r)$ of $L(\r)$ in $G(\r)$
such that $gLg^{-1} \supset M$, and where $\Theta_{gLg^{-1}}$ denotes the
virtual character on $gL(\r)g^{-1}$ obtained from $\Theta_L$ on $L(\r)$ via
the isomorphism 
$$
  \Int(g):L \to gLg^{-1} $$
($\Int(g)(x):=gxg^{-1}$). For regular $\g$ in $T(\r)$ the formula (4.6) is
just the usual formula for the character of a parabolically induced
representation, and by continuity the formula remains valid on all of
$T(\r)$. 

We finish this section by discussing stable discrete series characters.  We
now assume that there exists an elliptic maximal torus $T_e$ in $G$
(elliptic means that $T_e/Z$ is anisotropic, where $Z$ denotes the center
of $G$). As usual we let $q(G)$ denote half the dimension of the symmetric
space associated to the adjoint group of $G$ (our hypothesis on $G$
guarantees that this dimension is even). Let $\Pi$ be the L-packet
consisting of all (isomorphism classes of) discrete series representations of
$G(\r)$ having the same infinitesimal and central characters as the finite
dimensional representation $E$. Put
$$ 
\Theta^E=(-1)^{q(G)}\sum_{\pi \in \Pi} \Theta_\pi, $$ where $\Theta_\pi$
denotes the character of $\pi$; then the virtual character $\Theta^E$ is
stable \cite{HC,Lemma 61},\cite{S}.  
Any discrete series representation $\pi$ of
$G(\r)$ is obtained by induction from a discrete series representation of
the normal subgroup $$ Z(\r)\im[G_{\sc}(\r) \to G(\r)], $$ where $G_{\sc}$
denotes the simply connected cover of the derived group $G_{\der}$ of $G$.
Therefore $\Theta^E$ is supported on this normal subgroup (of finite index).

Let $A \subset T \subset M \subset G$ be as above. Let $\g \in T_{\reg}(\r)$
and define $R_\g$ as above. The character value $\Theta^E(\g)$ is given by
(4.1) for certain integers $n(\g,B)$. We will now review how these integers
are related to the stable discrete series constants discussed in \S3.
Let $T_c$ denote the
maximal anisotropic subtorus in $T$; note that $T=AT_c$. Let $L$
denote the centralizer of $T_c$ in $G$; then $L_\C$ is a Levi subgroup of
$G_\C$ and $L$ contains $T$. Note that the roots of $T$ in $L$ are
precisely the real roots of $T$.

Let $T(\r)_1$ denote the maximal compact subgroup of $T(\r)$. Then
$T_c(\r)$ is the identity component of $T(\r)_1$, and there is a direct
product decomposition 
$$
  T(\r)=A(\r)^0 \times T(\r)_1. \tag{4.7} $$

We decompose our regular element $\g \in T(\r)$ according to the
decomposition (4.7):
$$
   \g=\exp(x)\cdot \g_1 $$
for uniquely determined elements $x$ in $X_*(A)_\r=\Lie(A)$ and $\g_1$ in
$T(\r)_1$. Let $J$ denote the identity component of the 
centralizer of $\g_1$ in $L$. The root system of $T$ in $J$ is precisely
$R_\g$. 

Of course we may as well assume that $\g$ belongs to
$$
   Z(\r)\im[G_{\sc}(\r) \to G(\r)]; $$
otherwise $n(\g,B)=0$ for all $B \in \Cal B(T)$. In this case we claim that
$-1$ belongs to the Weyl group of $R_\g$. In proving the claim we may as
well assume that $\g$ lies in the image of $G_{\sc}(\r)$, and therefore we
may as well assume that $G_{\sc}=G$.
 Replacing $T_e$ by a conjugate, we may assume
that $T_c$ is contained in $T_e$; then $T_e$ is contained in $L$. Therefore
the connected center of $L$ is equal to $T_c$ (it contains $T_c$ and is
contained in both $T$ and $T_e$).
Moreover the maximal compact subgroups of $L(\r)$ are connected since the
derived group 
$L_{\der}$ of $L$ is 
simply connected and $L/L_{\der}$ is anisotropic. It follows
that by conjugating $T$ in $L$ we may assume that $\g_1$ belongs to $T_e$
and hence that $T_e$ is contained in $J$. Therefore the connected center of
$J$ is also equal to $T_c$. The maximal torus $T$ in $J$ is split modulo
the connected center $T_c$ of $J$, and therefore its split component $A$ is
a split maximal torus in $J_{\der}$. But $J$ contains an anisotropic
maximal torus, namely $T_e$, and therefore $-1 \in W(R_\g).$

Thus the root system $R_\g$ in $X^*(A/A_G)_\r$ is of the type considered in
\S3, and from this root system we obtain an integer-valued function $\bar
c$ on 
$$
   (X_*(A/A_G)_\r)_{\reg} \times (X^*(A/A_G)_\r)_{\reg} $$
(see \S3). The integer $n(\g,B)$ is given by
$$
n(\g,B)=\bar c(x,p(\l_B+\rho_B-\l_0)) $$
where 
$$
p:X^*(T)_\r \to X^*(A)_\r$$ 
is the natural restriction map and $\l_0 \in X^*(T)_\r$ is obtained from
the character $\l_0 \in X^*(A_G)$ by which $A_G$ acts on $E$ by viewing
$X^*(A_G)_\r$ as a direct summand of $X^*(T)_\r$ in the usual way. 
\heading 5. The stable virtual characters $\Theta_\nu$ \endheading
Let $F$ be a subfield of $\r$ (the two examples we have in mind are $\Q$
and $\r$). Let $G$ be a connected reductive group over $F$. By a parabolic
subgroup of $G$ we mean a parabolic subgroup of $G$ defined over $F$, and
by a Levi subgroup of $G$ we mean a Levi component defined over $F$ of some
parabolic subgroup of $G$. 

Let $M$ be a Levi subgroup of $G$. We write $\Cal F^G(M)$ for the set of
parabolic subgroups of $G$ containing $M$, and we write $\Cal P^G(M)$ for
the subset of $\Cal F^G(M)$ consisting of those parabolic subgroups for
which $M$ is a Levi component; we often abbreviate $\Cal F^G(M)$,$\Cal
P^G(M)$ to $\Cal F(M)$,$\Cal P(M)$. Let $A_M$ denote the maximal $F$-split
torus in the center of $M$, and write $\fA_M$ for the real vector space
$$
   \fA_M:=X_*(A_M)_\r, $$
where the subscript $\r$ indicates that we have tensored $X_*(A_M)$ over
$\Z$ with $\r$. By a root of $A_M$ we mean a non-zero weight of $A_M$ in 
$\Lie(G)$. Any $P \in \Cal P(M)$ determines a chamber $C_P$ in $\fA_M$,
consisting of the points $x \in \fA_M$ such that
$$
    \langle x,\a \rangle > 0 $$
for every root of $A_M$ in $\Lie(N)$, where $N$ denotes the unipotent
radical of $P$. The map $P \mapsto C_P$ is a bijection from $\Cal P(M)$ to
the set of chambers in $\fA_M$ (a chamber in $\fA_M$ is a connected
component of the complement in $\fA_M$ of the union of the root hyperplanes
in $\fA_M$). Let $Q \in \Cal F(M)$ and let $L$ denote the unique Levi
component of $Q$ containing $M$. Then $\fA_L$ is a subspace of $\fA_M$, so
that $C_Q$ can be regarded as a cone in $\fA_M$. Moreover $\fA_M$ is equal
to the disjoint union
$$
    \fA_M = \coprod _{Q \in \Cal F(M)} C_Q. $$

We fix a minimal parabolic subgroup $P_0$ of $G$ and fix a Levi component
$M_0$ of $P_0$ over $F$. We say that a parabolic subgroup $P$ of $G$ is \it
standard \rm if it contains $P_0$, and we say that $P$ is \it semistandard
\rm if it contains $M_0$. Thus $\Cal F(M_0)$ is the set of semistandard
parabolic subgroups. Given a semistandard parabolic subgroup $P$ of $G$, we
write $N_P$ for the unipotent radical of $P$ and $M_P$ for the unique Levi
component of $P$ containing $M_0$. We will often write $A_P$,$\fA_P$ rather
than $A_{M_P}$,$\fA_{M_P}$. In fact we will often abbreviate $M_P$,$N_P$ to
$M$,$N$, so that
$$P=MN.$$
When we use $Q$ to denote a semistandard parabolic subgroup, we will often
write $L$,$U$ instead of $M_Q$,$N_Q$, so that
$$ Q=LU.$$

Let $E$ be an irreducible representation of the algebraic group $G$ on a
finite dimensional complex vector space, and let $\nu$ be an element in
$(\fA_{P_0})^*$, the real vector space dual to $\fA_{P_0}$, such that the
restriction of $\nu$ to $\fA_G$ coincides with the character by which $A_G$
acts on $E$ (this character is an element of $X^*(A_G)$, a lattice in
$\fA^*_G$). 

We are going to use $E,\nu$ to define a virtual representation of the real
group $G(\r)$, the character of which we will denote by $\Theta_\nu$. The
first step is to define an element $\nu_P \in \fA_P^*$ for each semistandard
parabolic subgroup $P=MN$. There is a unique standard parabolic subgroup
$P'=M'N'$ conjugate to $P$ under $G(F)$. There exists $g \in G(F)$, unique
up to right multiplication by $M(F)$, such that $gPg^{-1}=P'$ and
$gMg^{-1}=M'$, and the inner automorphism $x \mapsto gxg^{-1}$ of $G$
induces isomorphisms
$$ \align
&A_P\simeq A_{P'} \\
&\fA_P \simeq \fA_{P'} \endalign
$$
independent of the choice of $g$. Let $\nu' \in \fA^*_{P'}$ be the
restriction of the linear form $\nu$ to the subspace $\fA_{P'}$ of
$\fA_{P_0}$, and then use the isomorphism $\fA_P\simeq \fA_{P'}$ to
transport $\nu'$ over to an element $\nu_P \in \fA^*_P$. This completes the
definition of $\nu_P$. It is easy to see that if $P,Q \in \Cal F^G(M_0)$
and $P \subset Q$, then $\nu_Q$ is the image of $\nu_P$ under the natural
restriction map
$$
   \fA^*_P \to \fA^*_Q. $$

The next step is to use $E,\nu$ to define a virtual finite dimensional
complex representation $E^\nu_P$ of $M$ for any semistandard parabolic
subgroup $P=MN$. We begin by considering the Lie algebra cohomology groups
$$
   H^i(\Lie(N),E); $$
these are finite dimensional complex representations of $M$ (we use the
usual \it left \rm action of $M$). Of course the action of the split torus
$A_P$ on $H^i(\Lie(N),E)$ decomposes this space as a direct sum of weight
subspaces
$$
   H^i(\Lie(N),E)_\mu,$$
where $\mu$ runs through $X^*(A_P)$ (a lattice in $\fA^*_P$). We write
$C^*_P$ for the closed convex cone in $\fA^*_P$ dual to $C_P$; thus $C^*_P$
consists of all $\mu \in \fA^*_P$ such that $\langle x,\mu \rangle \ge 0 $
for all $x \in C_P$. We write
$$
   H^i(\Lie(N),E)_{\ge \nu_P} $$
for the subspace of $H^i(\Lie(N),E)$ obtained by taking the direct sum of
all the weight spaces $H^i(\Lie(N),E)_\mu$ for $\mu \in X^*(A_P)$ such that
$\mu -\nu_P \in C^*_P$; of course $$H^i(\Lie(N),E)_{\ge \nu_P}$$ is stable
under the action of $M$. We write $E^\nu_P$ for the virtual $M$-module
$$
  \sum_i (-1)^i H^i(\Lie(N),E)_{\ge \nu_P}. $$

We now use a theorem of Kostant \cite{Ko} to express $E^\nu_P$ in terms of
irreducible representations of $M$. Let $T$ be a maximal torus of $M$ over
$\C$ and let $B_M$ be a Borel subgroup of $M$ over $\C$ containing $T$. For
any $B_M$-dominant weight $\mu \in X^*(T)$ we let $V^M_\mu$ be an
irreducible finite dimensional complex representation of $M$ with highest
weight $\mu$. The set of Borel subgroups $B$ of $G$ over $\C$ containing
$T$ and contained in $P$ is in natural bijection with the set of Borel
subgroups of $M$ over $\C$ containing $T$. Thus our choice of $B_M$
determines a Borel subgroup $B$ of $G$ over $\C$ containing $T$ and
contained in $P$, characterized by the equality
$$ 
   B_M=B \cap M. $$
Let $R$ (respectively, $R_M$) be the set of roots of $T$ in $G$
(respectively, $M$). The Borel subgroups $B,B_M$ determine positive systems
$R^+,R^+_M$ in $R,R_M$, and of course
$$
   R^+_M=R_M \cap R^+. $$

Let $\l_B \in X^*(T)$ denote the highest weight (with respect to $B$) of
the irreducible representation $E$ of $G$, and let $\rho_B \in X^*(T)_{\r}$
denote half the sum of the roots in $R^+$. Let $W$ (respectively, $W_M$)
denote the Weyl group of $T_\C$ in $G_\C$ (respectively, $M_\C$). Let $W'$
denote the set of Kostant representatives for the cosets $W_M\backslash W$;
thus $W'$ consists of the elements $ w \in W$ such that
$$
  w^{-1}(R^+_M) \subset R^+ $$
(obviously $W'$ depends on the choice of $B_M$). Kostant's theorem on
$\Lie(N)$-cohomology states that as an $M$-module
$H^i(\Lie(N),E)$ is isomorphic to
$$
   \bigoplus_w V^M_{w(\l_B+\rho_B)-\rho_B} $$
where $w$ runs through the set of Kostant representatives of length $i$ (we
use the length function on $W$ determined by $B$). Note that the weight
$w(\l_B+\rho_B)-\rho_B$ is indeed $B_M$-dominant for any Kostant
representative $w$.

Let 
$$
 \epsilon :W \to \{\pm 1\} $$
be the usual sign function on $W$ ($\epsilon(w)$ is $(-1)^{l(w)}$, where
$l(w)$ denotes the length of $w$). We see from Kostant's theorem that the
virtual representation $E^\nu_P$ is given by
$$
     \sum_{w \in W'} \epsilon(w) \cdot
\xi_{C^*_P}(p_M(w(\l_B+\rho_B)-\rho_B)-\nu_P) \cdot
V^M_{w(\l_B+\rho_B)-\rho_B} $$
where $\xi_{C^*_P}$ denotes the characteristic function of the subset
$C^*_P$ of $\fA^*_M$ and $p_M$ denotes the restriction map 
$$
  X^*(T)_\r=(X_*(T)_\r)^* \to \fA^*_M $$
induced by the inclusion of $A_M$ in $T$.

Now we are ready to define the virtual character $\Theta_\nu$ on the real
group $G(\r)$. For any semistandard parabolic subgroup $P=MN$ we write
$\d_P$ for the modulus quasicharacter on $M(\r)$ given by
$$
    \d_P(x)=|\det(x;\Lie(N))| $$
for $x \in M(\r)$. We write $(E^\nu_P)^*$ for the contragredient of the
virtual representation $E^\nu_P$. Then
$$
   \d_P^{-1/2} \otimes (E^\nu_P)^* $$
is a virtual representation of the real group $M(\r)$, which we may induce
from $P(\r)$ to $G(\r)$ to obtain a virtual representation
$$
   i^G_P(\d_P^{-1/2} \otimes (E^\nu_P)^*) $$
of $G(\r)$. We are using the usual \it normalized \rm parabolic induction,
which builds in a factor of $\d_P^{1/2}$; if we used unnormalized
induction we would simply be inducing $(E^\nu_P)^*$ from $P(\r)$ to
$G(\r)$. We write $\Theta^\nu_P$ for the (Harish-Chandra) character of
$$
   i^G_P(\d^{-1/2}_P \otimes (E^\nu_P)^*) $$
and define a virtual character $\Theta_\nu$ on $G(\r)$ by putting
$$
\Theta_\nu:=\sum_P (-1)^{\dim(A_P/A_G)} \Theta^\nu_P, $$
where $P$ runs over the set of standard parabolic subgroups of $G$. Note
that $\Theta_\nu$ is stable. Indeed, the character of $E^\nu_P$ is
obviously stable on $M(\r)$, and stability is preserved by parabolic
induction.

Now we fix a Levi subgroup $M$ of $G$ containing $M_0$, and we assume that
$M_\r$ contains a maximal torus $T$ over $\r$ such that $T/A_M$ is
anisotropic over $\r$. It follows that $A_M$ coincides with the maximal
$\r$-split torus in the center of $M$, and this in turn implies that any
parabolic subgroup of $G$ over $\r$ containing $M$ is automatically defined
over $F$. Note that $A_M$ is the maximal $\r$-split torus in $T$ and that
$T$ is elliptic in $M_\r$.
The discussion following Lemma 4.1  
 applies to the stable character $\Theta_\nu$, and
thus we obtain a continuous function $\Phi^G_M(\g,\Theta_\nu)$ on $T(\r)$.
Sometimes we abbreviate $\Phi^G_M(\g,\Theta_\nu)$ to
$\Phi_M(\g,\Theta_\nu)$. 

We are now going to use $E,\nu$ to define another function $L^\nu_M(\g)$ on
$T(\r)$ (with $M$ and $T$ as above); we will see in \S7 that this function
arises naturally in the Lefschetz trace formula for Hecke operators. 
Once the definition is complete our goal will be to 	
show that $\Phi_M(\g,\Theta_\nu)$ is in fact equal to $L^\nu_M(\g)$.
Let $\g \in T(\r)$. There is a direct product decomposition
$$
   T(\r)=A_M(\r)^0 \times T(\r)_1, $$
where $T(\r)_1$ denotes the maximal compact subgroup of $T(\r)$. Therefore
we can write $\g$ as
$$
     \g = \exp(x) \cdot \g_1 $$
for unique elements $x \in \fA_M$ and $\g_1 \in T(\r)_1$. 
The complex number $L^\nu_M(\g)$ that we are in the process of defining has
the form
$$
   L^\nu_M(\g):=\sum_Q (-1)^{\dim(A_L/A_G)} \cdot |D^L_M(\g)|^{1/2} \cdot
\d_Q^{-1/2}(\g) \cdot L^\nu_Q(\g) $$
where the sum runs over $Q=LU$ in $\Cal F(M)$ such that $x$ is contained in
the subspace $\fA_L$ of $\fA_M$ and where $L^\nu_Q(\g)$ is a complex number
we have yet to define. The factor $D^L_M(\g)$ was defined in \S4, just
before Lemma 4.1.  

In order to define $L^\nu_Q(\g)$ we choose a Borel subgroup $B$ of $G$ over
$\C$ containing $T$ and contained in $Q$, and we put $B_L:=B \cap L$, a
Borel subgroup of $L$ over $\C$ containing $T$; it turns out that
$L^\nu_Q(\g)$ is independent of this choice. We now use the same notational
system as we used when discussing Kostant's theorem (though we are
now using $Q,L$ instead of $P,M$). In particular we have the set $W'$ of
Kostant representatives for the cosets $W_L\backslash W$, the irreducible
representations 
$$
    V^L_{w(\l_B+\rho_B)-\rho_B} $$
of $L$, and the restriction map
$$ 
   p_L : X^*(T)_\r \to \fA_L^*. $$
The \it open \rm convex polyhedral cone $C_Q$ in $\fA_L$ determines a
function
$$
   \varphi_{C_Q}(\cdot,\cdot) $$
on $\fA_L \times \fA_L^*$, as in the last part of Appendix A, and we will
denote this function simply by $\varphi_Q(\cdot,\cdot)$. We define
$L^\nu_Q(\g)$ by
$$ \align
    L^\nu_Q(\g):=(-1)^{\dim(A_L)} \sum_{w \in W'} \epsilon(w) &\cdot
\tr(\g^{-1};V^L_{w(\l_B+\rho_B)-\rho_B}) \\ &\cdot 
\varphi_Q(-x,p_L(w(\l_B+\rho_B)-\rho_B)-\nu_Q). \endalign $$
\proclaim{Theorem 5.1} The two functions $\Phi_M(\g,\Theta_\nu)$ and
$L^\nu_M(\g)$ on $T(\r)$ are equal. \endproclaim
By definition $\Phi_M(\g,\Theta_\nu)$ is given by 
$$
   \sum_Q (-1)^{\dim(A_L/A_G)} \cdot \Phi_M(\g,\Theta^\nu_Q) \tag{5.1} $$
where $Q=LU$ runs over the set of standard parabolic subgroups of $G$.
Applying equation (4.6) to the induced character $\Theta^\nu_Q$ of $G(\r)$,
we see that $\Phi_M(\g,\Theta_\nu)$ is equal to 
$$
   \sum_Q (-1)^{\dim(A_L/A_G)} \sum_{Q'}
\Phi^{L'}_M(\g,\d_{Q'}^{-1/2}\otimes(E^\nu_{Q'})^*) \tag{5.2} $$
where the index set for the first sum is the same as before and the index
set for the second sum is the set of parabolic subgroups $Q'$ over $\r$
containing $M$ such that $Q'$ is conjugate under $G(\r)$ to $Q$. Since, as
we remarked earlier, every parabolic subgroup of $G$ over $\r$ containing
$M$ is automatically defined over $F$, we see that $\Phi_M(\g,\Theta_\nu)$
is equal to
$$
   \sum_{Q \in \Cal F(M)} (-1)^{\dim(A_L/A_G)}
\Phi^L_M(\g,\d_Q^{-1/2}\otimes(E^\nu_Q)^*) \tag{5.3} $$
(as usual $Q=LU$). Recalling the expression for $E^\nu_Q$ that we found
using Kostant's theorem, we see that $\Phi_M(\g,\Theta_\nu)$ is equal to 
$$ \align
   \sum_{Q \in \Cal F(M)} &(-1)^{\dim(A_L/A_G)} \cdot |D^L_M(\g)|^{1/2}
\cdot \d_Q^{-1/2}(\g) \cdot \sum_{w \in W'} \epsilon(w) \\ &\cdot
\tr(\g^{-1};V^L_{w(\l_B+\rho_B)-\rho_B}) \cdot
\xi_{C_Q^*}(p_L(w(\l_B+\rho_B)-\rho_B)-\nu_Q). \tag{5.4} \endalign $$
The notation here is the same as that used during our discussion of
Kostant's theorem. In particular, given $Q \in \Cal F(M)$ we must choose a
Borel subgroup $B$ of $G$ over $\C$ containing $T$ and contained in $Q$ in
order to define $W'$, the set of Kostant representatives. 

Let $Q=LU$ be a parabolic subgroup in $\Cal F(M)$. As usual $\fA_M$ is a
disjoint union of convex cones $C_{Q'}$, one for each $Q' \in \Cal F(M)$.
But $M$ is also a Levi subgroup of $L$, and therefore $\fA_M$ is also a
disjoint union of convex cones $C_{Q''}$, one for each $Q'' \in \Cal
F^L(M)$. For any $Q' \in \Cal F(M)$ such that $Q' \subset Q$ we put
$Q'':=Q' \cap L$, an element of $\Cal F^L(M)$. The map $Q' \mapsto Q''$
sets up a bijection 
$$
  \{Q' \in \Cal F(M) \,|\, Q' \subset Q \} \simeq \Cal F^L(M), $$
and the convex cones $C_{Q'}$,$C_{Q''}$ in $\fA_M$ are related by the
equality 
$$
   C_{Q''}=C_{Q'}+\fA_L. $$

Recall that we have written $\g$ as
$$
    \g=\exp(x) \cdot \g_1 $$
for uniquely determined $x \in \fA_M$ and $\g_1 \in T(\r)_1$. 
For each parabolic subgroup $Q \in \Cal F(M)$ we denote by $Q'=L'U'$ the
unique element of $\Cal F(M)$ such that
\roster \item $Q' \subset Q$, and \item $-x \in C_{Q'}+\fA_L.$
\endroster 
It follows from the second condition that $x$ belongs to $\fA_{L'}$.

Now let $Q_1=L_1U_1$ be an element of $\Cal F(M)$ such that $x \in
\fA_{L_1}$. Pick a Borel subgroup $B$ in $G$ over $\C$ containing $T$ and
contained in $Q_1$. We are interested in the terms in (5.4) indexed by
parabolic subgroups $Q \in \Cal F(M)$ such that $Q' = Q_1$. We have
inclusions
$$
   T \subset B \subset Q_1 \subset Q, $$
so that we can (and do) use $B$ to define the set $W'$ of Kostant
representatives for the cosets $W_L\backslash W_G$. As before we write
$Q''$ for the element $Q_1 \cap L$ of $\Cal F^L(M)$. Define a function
$\D^L_{Q''}$ on $T(\r)$ (a partial Weyl denominator for the group $L$) by 
$$
    \D^L_{Q''}(\g)=\prod_\a (1-\a(\g)^{-1}) $$
where $\a$ runs through the set of roots of $T$ in $\Lie(N'')$ ($N''$
denotes the unipotent radical of $Q''$). We claim that $\D^L_{Q''}(\g^{-1})$
is a non-negative real number (for $\g$,$Q$,$Q_1$,$Q''$ as in the
discussion preceding the definition of $\D^L_{Q''}$). Since $Q''$ is
defined over $\r$, the set of roots $\a$ of $T$ in $\Lie(N'')$ is stable
under complex conjugation. Complex conjugate pairs $\bar \a$,$\a$ with
$\bar \a\ne\a$ make a non-negative contribution to $\D^L_{Q''}$, since
$1-\bar\a(\g)$ is complex conjugate to $1-\a(\g)$. Let $\a$ be a root of
$T$ in $\Lie(N'')$ such that $\bar\a=\a$. It is enough to show that
$1-\a(\g)$ is non-negative. Since $Q'=Q_1$, the element $-x$ belongs to
$C_{Q''}$, which implies that
$$
    \langle x,\a \rangle \le 0. $$
Since $\a=\bar\a$, the value of $\a$ on any element of $T(\r)_1$ is $\pm1$.
Therefore $\a(\g_1)=\pm 1$ and 
$$
   \a(\g)=\exp(\langle x,\a \rangle) \cdot \a(\g_1) \le 1, $$
as desired.

It follows from the claim that
$$
   |D^L_{L_1}(\g)|^{1/2} \cdot \d_Q^{-1/2}(\g) = \d_{Q_1}^{-1/2}(\g) \cdot
\D^L_{Q''}(\g^{-1}).$$
Moreover, applying Kostant's theorem to $L$ and its parabolic subgroup
$Q''$, it is easy to see that for $w \in W'$ 
$$
   \D^L_{Q''}(\g^{-1}) \cdot \tr(\g^{-1};V^L_{w(\l_B+\rho_B)-\rho_B}) $$
is equal to
$$
\sum _{u \in W'_L} \epsilon(u) \cdot
\tr(\g^{-1};V^{L_1}_{uw(\l_B+\rho_B)-\rho_B}) $$ where $W'_L$ is the set of
Kostant representatives for the cosets $W_{L_1}\backslash W_L$ (relative to
the Borel subgroup $B \cap L$ of $L$). It is also easy to see that
$$W'_LW'$$ is the set $W''$ of Kostant representatives for the cosets
$W_{L_1}\backslash W$ (relative to the Borel subgroup $B$ of $G$) and that
for $w \in W'$, $u \in W'_L$ 
$$
   p_L(uw(\l_B+\rho_B))=p_L(w(\l_B+\rho_B)). $$
Therefore the contribution of such a $Q$ to (5.4) is
$$ \align
   (-1)^{\dim(A_L/A_G)} \cdot |D^{L_1}_M(\g)|^{1/2} &\cdot
\d_{Q_1}^{-1/2}(\g) \cdot \sum_{w \in W''} \epsilon(w) \cdot
\tr(\g^{-1};V^{L_1}_{w(\l_B+\rho_B)-\rho_B}) \\ &\cdot
\xi_{C^*_Q}(p_L(w(\l_B+\rho_B)-\rho_B)-\nu_Q). \endalign $$
Since all factors in this expression except for the first and last depend
on $Q$ only through $Q_1$, we see that $\Phi_M(\g,\Theta_\nu)$ is equal to
$$ \align
   \sum_{Q_1}(-1)^{\dim(A_{L_1}/A_G)} &\cdot |D^{L_1}_M(\g)|^{1/2} \cdot
\d_{Q_1}^{-1/2}(\g) \cdot \sum_ {w \in W''} \epsilon(w) \cdot
\tr(\g^{-1};V^{L_1}_{w(\l_B+\rho_B)-\rho_B}) \\ &\cdot
\sum_Q(-1)^{\dim(A_{L_1}/A_L)} \cdot
\xi_{C^*_Q}(p_L(w(\l_B+\rho_B)-\rho_B)-\nu_Q), \tag{5.5} \endalign $$
where the index set for the first sum is the set of $Q_1=L_1U_1 \in \Cal
F(M)$ such that $x \in \fA_{L_1}$ and the index set for the second sum is
the set of $Q \in \Cal F(M)$ such that $Q'=Q_1$.

Comparing (5.5) with the definition of $L^\nu_M(\g)$, we see that in order
to prove that $\Phi_M(\g,\Theta_\nu)$ is equal to $L^\nu_M(\g)$, it is
enough to prove  the equality 
$$ \align
\sum_Q
(-1)^{\dim(A_{L_1}/A_L)} &\cdot
\xi_{C^*_Q}(p_L(w(\l_B+\rho_B)-\rho_B)-\nu_Q)
\\ &=(-1)^{\dim(A_{L_1})}\cdot
\varphi_{Q_1}(-x,p_{L_1}(w(\l_B+\rho_B)-\rho_B)-\nu_{Q_1}). \tag{5.6}
\endalign $$
The sum in (5.6) is taken over the set of $Q \in \Cal F(M)$ such that $Q' =
Q_1$, or, in other words, the set of $Q \in \Cal F(M)$ such that $ Q
\supset Q_1$ and $-x \in C_{Q_1}+\fA_L$. Therefore the left-hand side of
(5.6) is equal to
$$
   \sum_Q (-1)^{\dim(A_{L_1}/A_L)} \cdot \xi_{C_{Q_1}+\fA_L}(-x) \cdot
\xi_{C_Q^*}(p_L(w(\l_B+\rho_B)-\rho_B)-\nu_Q), $$
where the sum is now taken over all $Q \in \Cal F(M)$ such that $Q \supset
Q_1$, 
and this is indeed equal to the right-hand side of (5.6), as one easily
sees from Lemma A.6 and the fact that the restriction of
$$
   p_{L_1}(w(\l_B+\rho_B)-\rho_B)-\nu_{Q_1} $$
to the subspace $\fA_L$ of $\fA_{L_1}$ is equal to
$$
   p_L(w(\l_B+\rho_B)-\rho_B)-\nu_Q. $$
This concludes the proof of Theorem 5.1.

We will now make a particular choice for $\nu$, and we will show (Theorem
5.2) that with this choice the virtual character $\Theta_\nu$  becomes
especially simple. Denote by $N_0$ the unipotent radical of our chosen
minimal parabolic subgroup $P_0$. Let $\rho_0 \in \fA^*_{P_0}$ denote half
the sum of the roots (counted with multiplicity) of $A_{P_0}$ in
$\Lie(N_0)$. As usual we regard $\fA^*_G$ as a direct summand of
$\fA^*_{P_0}$. Define $\nu_m \in \fA^*_{P_0}$ by
$$
    \nu_m:=-\rho_0 + \l_0, $$
where $\l_0 \in X^*(A_G) \subset \fA^*_G$ is the character by which $A_G$
acts in the representation $E$. From $\nu_m$ we obtain the virtual
character $\Theta_{\nu_m}$ on $G(\r)$.

Suppose first that there exists an elliptic maximal torus $T_e$ in $G$ over
$\r$, so that $G(\r)$ has a discrete series. As in \S4 we put
$$
    \Theta^E= (-1)^{q(G)} \sum_{\pi \in \Pi} \Theta_\pi, $$
where $\Pi$ is the L-packet of discrete series representations of $G(\r)$
having the same infinitesimal and central characters as the finite
dimensional representation $E$. Note that the contragredient of $\Theta^E$
is equal to $\Theta^{E^*}$, where $E^*$ denotes the contragredient of $E$.
\proclaim{Theorem 5.2} The virtual characters $\Theta_{\nu_m}$ and
$\Theta^{E^*}$ agree on $T_{\reg}(\r)$ for any maximal torus $T$ in $G$
over $\r$ whose $\r$-split component is both defined and split 
over $F$. \endproclaim
It is enough to prove the theorem when $F=\r$, in which case we must
show that $\Theta_{\nu_m}$ is equal to $\Theta^{E^*}$. We appeal to the
characterization of $\Theta^{E^*}$ provided by Theorem 3 of \cite{HC}.
Clearly $\Theta_{\nu_m}$ and $\Theta^{E^*}$ are invariant distributions
with the same infinitesimal and central characters. Moreover it is obvious
from the definition of $\Theta_{\nu_m}$ that $\Theta_{\nu_m}$ agrees with
$\Theta^{E^*}$ on $T_e(\r)\cap G_{\reg}(\r)$ (on $T_e(\r) \cap
G_{\reg}(\r)$ the virtual character $\Theta^{E^*}$ coincides with that of
the finite dimensional representation $E^*$).
Thus the only non-trivial point is to check the validity of the second
condition in Harish-Chandra's theorem:
$$
\sup_{\g \in G_{\reg}(\r)} |D(\g)|^{1/2}|\Theta_{\nu_m}(\g)|\omega(\g) <
\infty , $$
where $\omega$ is the unique homomorphism from $G(\r)$ to the group of
positive real numbers whose restriction to $A_G(\r)$ is equal to the
absolute value of the quasi-character by which $A_G(\r)$ acts in $E$, and
where
$$
   D(\g)=\det(1-\Ad(\g);\Lie(G)/\Lie(T)), $$
where $T$ denotes the unique maximal torus of $G$ containing $\g$. 

It is enough to check that for every maximal torus $T$ of $G$
$$
    \Phi_M(\g,\Theta_{\nu_m}) \cdot \omega(\g) \tag{5.7} $$
is bounded on $T_{\reg}(\r)$. Here $M$ is (as usual) the centralizer of the
split component of $T$; we used that roots of $T$ in $M$ are imaginary and
hence that
$$
    |\det(1-\Ad(\g);\Lie(M)/\Lie(T))| $$
is bounded on $T(\r)$. But the boundedness of (5.7) follows directly from
Theorem 5.1. Indeed, it is enough to check that for all $P=MN \in \Cal
P(M)$ 
$$
    \tr(\g^{-1};V^M_{w(\l_B+\rho_B)-\rho_B}) \cdot \omega(\g) \cdot
\d_P^{-1/2}(\g) \tag{5.8} 
$$ 
is bounded on $T(\r)$ whenever $w \in W'$ is such that 
$$
   \varphi_P(-x,p_M(w(\l_B+\rho_B)-\l_0)) \ne 0. \tag{5.9} $$
By Proposition A.5 the condition (5.9) implies that
$$
 \langle x, p_M(w(\l_B+\rho_B)-\l_0) \rangle \ge 0. \tag{5.10} $$
Since
$$
   T(\r)=A_M(\r)^0 \times T(\r)_1$$ 
and $T(\r)_1$ is compact, only the character by which $A_M(\r)^0$ acts in
$V^M_{w(\l_B+\rho_B)-\rho_B} $ is relevant to the boundedness of (5.8). In
fact the function (5.8) of $\g$ transforms under $A_M(\r)^0$ by the element
$$
   p_M(-w(\l_B+\rho_B)+\rho_B+\l_0-\rho_N)=p_M(-w(\l_B+\rho_B)+\l_0) \in
\fA^*_M $$
where $\rho_N$ is half the sum of the roots of $T$ in $\Lie(N)$, and thus
(5.10) does imply that (5.8) is bounded on $T(\r)$.  This completes the
proof.

Now we drop the assumption that $G$ has an elliptic maximal torus over
$\r$. First, for arbitrary $G$ and suitably regular $\g$ 
we will rewrite $\Phi_M(\g,\Theta_{\nu_m})$
in terms of the functions $\psi_R$ of \S1. Second, for $G$ having no
elliptic maximal torus we will use this expression for
$\Phi_M(\g,\Theta_{\nu_m})$ to show that it vanishes under a certain
regularity hypothesis on the highest weight of $E$.

Let $T \subset M \subset G$ be as usual. Let $\g \in T(\r)$ and write
$$
   \g = \exp(x) \cdot \g_1 $$
with $x \in \fA_M$ and $\g_1 \in T(\r)_1$. Assume that $x$ is regular in
$\fA_M$, in the sense that no root of $A_M$ vanishes on $x$. Then by
(A.2) (or rather its analog for the function $\varphi_Q$)
$L^\nu_M(\g)$ is given by 
$$ \align
    L^\nu_M(\g)=(-1)^{\dim(A_G)}\sum_{P \in \Cal P(M)} \d_P^{-1/2}(\g)
\cdot \sum_{w\in W'} 
&\epsilon(w) \cdot \tr(\g^{-1};V^M_{w(\l_B+\rho_B)-\rho_B}) 
\\ &\cdot \varphi_P(-x,p_M(w(\l_B+\rho_B)-\rho_B)-\nu_P). \endalign $$
Fix a Borel subgroup $B_M$ of $M$ over $\C$ containing $T$. For each $P \in
\Cal P(M)$ we let $B(P)$ denote the unique Borel subgroup of $G$ over $\C$
such that
$$
   T \subset B(P) \subset P $$
and $B(P) \cap M = B_M$. Write $P_x=MN_x$ for the unique element of $\Cal
P(M)$ whose chamber in $\fA_M$ contains $-x$, and write $B(x)$ for
$B(P_x)$. Let $\Cal B(T)$ denote the set of Borel subgroups of $G$ over
$\C$ containing $T$, and let $\Cal B(T)'$ denote the subset of $\Cal B(T)$
consisting of Borel subgroups $B$ such that $B \cap M=B_M$. Then
$$ \align
   L^\nu_M(\g)=(-1)^{\dim(A_G)} \sum_{B \in \Cal B(T)'}\sum_{P=MN \in \Cal
P(M)} \epsilon(B,B(P)) &\cdot \d_P^{-1/2}(\g) \cdot
\tr(\g^{-1};V^M_{\l_B+\rho_B-\rho_{B(P)}}) 
\\ &\cdot \varphi_P(-x,p_M(\l_B+\rho_B-\rho_N)-\nu_P) \endalign $$
where $\epsilon(B,B(P))=\epsilon(w)$ for the unique element $w \in W$ such
that $wB(P)w^{-1}=B$ and where $\rho_N$ is half the sum of the roots of $T$
in $\Lie(N)$. Of course
$$
    \rho_{B(P)}=\rho_M+\rho_N, $$
where $\rho_M$ denotes half the sum of the roots of $T$ in $M$ that are
positive for $B_M$. Thus
$$
     \rho_{B(P)}-\rho_{B(x)}=\rho_N-\rho_{N_x} \in X^*(T) $$
is trivial on the intersection of $T$ with the derived group of $M$, and
therefore defines a homomorphism 
$$
    M \to   \Bbb G_m. $$
It follows that
$$
    \tr(\g^{-1};V^M_{\l_B+\rho_B-\rho_{B(P)}})=\tr(\g^{-1};V^M_{\l_B+\rho_B-\rho_{B(x)}})\cdot \langle \g^{-1},\rho_{B(x)}-\rho_{B(P)} \rangle $$
and this shows that
$$
   L^\nu_M(\g)=(-1)^{\dim(A_G)} \sum _{B \in \Cal B(T)'} \tr(\g^{-1};
V^M_{\l_B+\rho_B-\rho_{B(x)}}) \cdot \d_{P_x}^{-1/2}(\g) \cdot \epsilon(B,B(x))
\cdot a_B $$
where
$$ \align
   a_B= \sum_{P=MN \in \Cal P(M)} \langle \g^{-1}, \rho_{B(x)}-\rho_{B(P)}
\rangle &\cdot 
\d_{P_x}^{1/2}(\g) \cdot \d_P^{-1/2}(\g) \cdot \epsilon(B(P),B(x))
\\ &\cdot \varphi_P(-x,p_M(\l_B+\rho_B-\rho_N)-\nu_P). \endalign $$
Let $\overline P=M \overline N$ be the element of $\Cal P(M)$ opposite $P$.
Note that
$$
    \langle \g^{-1},\rho_{B(x)}-\rho_{B(P)} \rangle \d^{1/2}_{P_x}(\g)
\d_P^{-1/2}(\g) $$ 
is equal to the sign of the real number
$$
    \langle \g^{-1},\rho_{B(x)}-\rho_{B(P)} \rangle =\langle
\g^{-1},\rho_{N_x}-\rho_N \rangle. $$
But $\rho_{N_x}-\rho_N$ is the sum of all roots $\a$ of $T$ in $N_x \cap
\overline N$, and since this set of roots is preserved by complex
conjugation, the sign in question is
$$
    \prod_\a \sgn \a(\g^{-1}) $$
where $\a$ runs through the real roots of $T$ in $N_x \cap \overline N$.
For such a real root
$$
  \sgn \a(\g^{-1})=\sgn \a(\g_1). $$
Of course the sign $\epsilon(B(P),B(x))$ is $-1$ raised to the number of
roots of $T$ in $N_x \cap \overline N$, and again since this set of roots
is stable under complex conjugation, this sign is $-1$ raised to the number
of real roots of $T$ in $N_x \cap \overline N$.
Let $R_\g$ be the set of real roots $\a$ of $T$ in $G$ such that
$\a(\g_1)=1 $ (or, equivalently, such that $\a(\g) > 0$). Of course $R_\g$
is a root system in $(\fA_M/\fA_G)^*$ (though it need not span that space).
Let $\Cal C$ be the set of Weyl
chambers in $\fA_M$ for the root system $R_\g$. For $C_1$,$C_2 \in \Cal C$
let $\epsilon(C_1,C_2)$ be $-1$ raised to the number of root hyperplanes
(for roots in $R_\g$) separating $C_1$ and $C_2$. Let $C_0$ be the
unique element in $\Cal C$ that contains $-x$. For $P \in \Cal P(M)$ let
$C_P$ denote the unique element of $\Cal C$ that contains the chamber in
$\fA_M$ determined by $P$. Then it follows from the discussion above that
$$
  a_B = \sum_{P=MN \in \Cal P(M)} \epsilon(C_0,C_P)
\varphi_P(-x,p_M(\l_B+\rho_B-\rho_N)-\nu_P). $$

Now suppose that $\nu=\nu_m$. Then 
$$
   p_M(\l_B+\rho_B-\rho_N)-\nu_P=p_M(\l_B+\rho_B-\l_0) $$
is independent of $P$. Since $x$ is regular, it follows from (A.13) and 
Lemma A.4 that
$a_B$ is equal to 
$$
    \sum_{C \in \Cal C} \epsilon(C_0,C) \psi_{\overline
C}(-x,p_M(\l_B+\rho_B-\l_0)). \tag{5.11} $$
This is nothing but
$$
   (-1)^{\dim(A_G)} \psi_{R_\g}(C_0,-x,p_M(\l_B+\rho_B-\l_0)), $$
(strictly speaking we only defined the functions $\psi_R$ for root systems
spanning the vector space in which they lie, but of course the definition
extends immediately to the general case). Note that by (A.2) each term of
(5.11) vanishes unless $p_M(\l_B+\rho_B-\l_0)$ belongs to the span of
$R_\g$. 

Now suppose that the highest weight of $E$ satisfies the following
property: for every proper Levi subgroup $M$ of $G$, for every maximal
torus $T$ of $M$ over $\C$ and for every $B \in \Cal B(T)$ the element 
$$
   \l_B+\rho_B-\l_0 $$
of $X^*(T)$ is non-trivial on $A_M$.

\proclaim{Theorem 5.3} Assume that $G$ contains no elliptic maximal torus
over $\r$.
Then, under the hypothesis above on the highest weight of $E$, the complex
number $\Phi_M(\g,\Theta_{\nu_m})$ is $0$, and consequently, if $F$ is
$\r$, the virtual character $\Theta_{\nu_m}$ is $0$. \endproclaim

By Theorem 5.1 we must show that $L^\nu_M$ is 0. By continuity it is enough
to show this for $\g \in T(\r)$ such that $x$ is regular in $\fA_M$. Then,
by the discussion above, it is enough to show that (5.11) vanishes. As in
\S4 let $T_c$ denote the maximal anisotropic subtorus in $T$ and let $J$
denote the centralizer in $G$ of $T_c$ and $\g_1$. Then $R_\g$ is the root
system $R_J$ of $T$ in $J$. Clearly $T_c$ is central in $J$ and $J/T_c$ is
a split group with split maximal torus $T/T_c$.

Since by hypothesis $G/A_G$ has no anisotropic maximal torus over $\r$, the
same is true of $J/A_G$ and $J/T_c A_G$. Since $J/T_c A_G$ is a split
group, either $A_J$ is strictly bigger than $A_G$ or $A_J=A_G$ (in which
case $R_J$ spans $(\fA_M/\fA_G)^*$) and $-1_{\fA_M/\fA_G}$ does not belong
to the Weyl group of $R_J$. In the first case $\fA_J$ is of the form
$\fA_L$ for some proper Levi subgroup $L$ of $G$ containing $M$, and
therefore our hypothesis on the highest weight of $E$ implies that every
term of (5.11) vanishes. In the second case
Corollary 1.3 (applied to the root system $R_\g$ in $\fA_M/A_G$) implies
that (5.11) is 0, since our hypothesis on the highest weight of $E$ ensures
that $p_M(\l_B+\rho_B-\l_0)$ is $R_\g$-regular (any intersection of root
hyperplanes in $\fA_M/\fA_G$ for roots in $R_\g$ is of the form
$\fA_L/\fA_G$ for some Levi subgroup $L$ of $G$ containing $M$, and
$R_\g$-regularity is equivalent to non-vanishing on every such non-zero
intersection). 
\heading 6. Discrete series constants \endheading
In the beginning of \S4 (see (4.1)--(4.3)) we reviewed the form taken by
the character of an irreducible representation of a real reductive
group. In this section we are concerned with the case of discrete series
representations. We then refer to the integers $n(\g,B)$ appearing in (4.1)
as {\it discrete series constants\/} (however we will no longer use the
notation $n(\g,B)$). In this section we give a simple formula for the
discrete series constants. Because of a descent property satisfied by the
constants (see \cite{K,13.4}) it is enough to give the formula in the
following special case.

Let $G$ be a split semisimple simply connected group over $\r$, and assume
that $G$ contains an anisotropic maximal torus~$T_e$. Let $A$ be a split
maximal torus in~$G$. We choose an isomorphism $A \simeq T_e$ over~$\C$
that is induced by an inner automorphism of~$G$ over~$\C$ and use it to
identify the character groups of~$A$ and~$T_e$. We put
$$ X^*:=X^*(T_e)_\r \simeq X^*(A)_\r $$ and
$$ X:= X_*(T_e)_\r\simeq X_*(A)_\r.$$
The roots and coroots of~$T_e$ in~$G$ give us a root system
$(X,X^*,R,R^\vee)$. Of course the set~$R$ spans $X^*$, and $-1$ belongs to
the Weyl group $W=W(R)$. 

Let $\tau$ be a regular element in $X^*(T_e)$. Associated to $\tau$ is a
discrete series representation $\pi(\tau)$ of $G(\r)$ having infinitesimal
character~$\tau$ and having the same central character as the finite
dimensional representation having infinitesimal character~$\tau$. 
We are interested in the constants needed to express the
values of the character of~$\pi(\tau)$ at regular elements in the identity
component $A(\r)^0$ of $A(\r)$; this is the special case alluded to
above. 

We need a little preparation before we can state our formula for these
constants. There is a unique maximal compact subgroup $K$ of $G(\r)$
containing $T_e(\r)$, and the roots of~$T_e$ in~$K$ form a subset
$R_c$ of~$R$ (such roots are said to be {\it compact\/}). We write $W_c$ for
the Weyl group of $R_c$ and identify it with a subgroup of~$W$. It is not
hard to see that the normalizer $\widetilde W_c$ of $W_c$ in~$W$ is given by
$$
\widetilde W_c=\{ w \in W \,|\, w(R_c)=R_c \}. $$

Let $C$
be a chamber in $X$. The chamber $C$ determines a subset $R_C$ of $R$ in
the following way. Let $\d_C \in X$ denote the half-sum of the coroots that
are positive for $C$, and put
$$
   R_C = \{ \a \in R \,|\, \a(\d_C) \in 2\Z \}; $$
note that no simple root (for $C$) belongs to $R_C$. 
We denote by $W_C$ the Weyl group of $R_C$. 
We identify $W_C$ with a subgroup of~$W$, and let $\widetilde 
W_C$ denote its normalizer in~$W$.  As before we write $C^\vee$
for the Weyl chamber in $X^*$ corresponding to $C$.

It is known (see \cite{AV,6.24(f)}) that there exists a chamber $C$ such
that $R_C$ equals $R_c$. For such a chamber $C$ we have $W_C=W_c$ and
$\widetilde W_C=\widetilde W_c$. The $\widetilde W_c$-orbit of $C$ is
uniquely determined by the condition that $R_C$ equal $R_c$.

As before let $\tau$ be a regular element of $X^*(T_e)$ (actually, in the
definition we are about to make we could just as well let $\tau$ be any
regular element in $X^*(T_e)_\r=X^*$). Let $C$ be a chamber in $X$ (for the
time being we do {\it not\/} assume that $R_C=R_c$). Let $x$ be an 
$R^\vee$-regular element of $X$, and let $\l$ be an element of $X^*$
lying in the $W$-orbit $W\cdot \tau$ of $\tau$ (see \S1 for the definition
of $R^\vee$-regularity). 
We define an integer
$b_R(\tau,C;x,\l)$ by 
$$
b_R(\tau,C;x,\l)=(-1)^{q(R)}\sum_{w \in W(\tau,C,\l)} \epsilon(x,wC)
\psi_{wC^\vee}(\l,x). \tag{6.1} $$
In case $R_C=R_c$ these constants (for $\l \in W\cdot \tau$) are the ones
needed to express the value of the character of $\pi(\tau)$ at the point $a
\in A(\r)^0 $ obtained from $x \in X$ via the exponential map (we have
identified $X_*(A)$ with $X_*(T_e)$ and thus we may view $X$ as the Lie
algebra of $A(\r)$). However for technical reasons it is best to define
$b_R(\tau,C;x,\l)$ for {\it any\/} chamber~$C$. 

The expression (6.1) requires some explanation. The integer
$$
   q(R)=[|R^+|+\dim(X)]/2 $$
was used already in \S2, and its interpretation in terms of $G$ was also
given there. The index set for the sum is the coset
$$
    W(\tau,C,\l):=\{ w \in W \, | \, w^{-1}\l \in W_C\cdot \tau \} $$
of $W_C$ in $W$ (of course $W_C \cdot \tau$ denotes the orbit of $\tau$
under $W_C$). As in \S1, for any two chambers $C_1,C_2$ in $X$ we write
$\epsilon(C_1,C_2)$ for the sign of the Weyl group element $w$ such that
$wC_1=C_2$. The sign $\epsilon(x,wC)$ appearing in (6.1) is by definition
$\epsilon(C_x,wC)$, where $C_x$ denotes the unique chamber in $X$
containing the (regular) element $x \in X$. Finally $\psi_{wC^\vee}(\l,x)$
is the function of $(\l,x) \in X^*\times X$ defined in \S1; it is obtained
from the coroot system $R^\vee$ and the Weyl chamber $wC^\vee$ in $X^*$.

Of course  the integer $b_R(\tau,C;x,\l)$ depends
only on the $R^\vee$-chamber of $X$ in which $x$ lies. But in fact we claim
that $b_R(\tau,C;x,\l)$ depends only on the Weyl chamber of $X$ in which
$x$ lies ($x$ is still assumed to be $R^\vee$-regular). 
Indeed it follows from Lemma 1.5 (applied to $R^\vee$) that for
all $s \in \hat A_{\sc}$ such that $s^2 \in Z^\vee$
$$
\sum_D \epsilon(D_0^\vee,D^\vee)\langle \d_D-\d_{D_0},s\rangle
\psi_{D^\vee}(\l,x) \tag{6.2} $$
depends only on the Weyl chamber in which $x$ lies. Here we are using the
notation $\d_D$,$\hat A_{\sc}$,$s$,$Z^\vee$
of Lemma 1.4 (since we are applying
Lemma 1.5 to $R^\vee$ rather than $R$, we need the notation of Lemma 1.4),
and we have applied Lemma 1.2 to the root system $R^\vee_s$ of Lemma 1.4.
Summing (6.2) over all $s \in \hat A_{\sc}$ such that $s^2 \in Z^\vee$ we
find that
$$
\sum_D \epsilon(D_0^\vee,D^\vee)\psi_{D^\vee}(\l,x) \tag{6.3} $$
depends only on the Weyl chamber in which $x$ lies, where the sum is now
taken over all chambers $D$ such that the element
$\d_D-\d_{D_0} \in Q$ lies in $2Q$ (here $Q \subset X$ denotes the lattice
generated by $R^\vee$). This set of chambers can also be described as the
set of chambers $wD_0$ where $w$ ranges through the stabilizer in $W$ of
the element $\d_{D_0} \in (\frac12 Q)/2Q$. We can think of $(\frac12 Q)/2Q$
as consisting of 4-torsion elements in a maximal torus $A_{\sc}$ in the
semisimple simply connected complex group $G_{\sc}$ with root system $R$,
and by a theorem of Steinberg this stabilizer is the Weyl group of the root
system of the centralizer of $\d_{D_0} \in A_{\sc}$ in $G_{\sc}$, namely 
$$
\{ \a \in R \, | \, \langle \d_{D_0},\a \rangle \in 2\Z \} = R_{D_0}. $$
Therefore the sum in (6.3) is over $wD_0$ ($w \in W_{D_0}$) and our claim
has been proved (take $D_0=w_0C$ for any $w_0 \in W(\tau,C,\l)$).

In order to prove that the integers $b_R(\tau,C;x,\l)$ are the ones
appearing in the character formula for $\pi(\tau)$ on $A(\r)^0$, we
must show that they satisfy various properties.
 The first is that
$$
   b_R(\tau,C;x,\l)=b_R(\tau,wC;x,\l) \quad \text{for all $w \in \tilde
W_C$} \tag{6.4} $$
(in other words $b_R(\tau,C;x,\l)$ only depends on the subset $R_C$
determined by~$C$).
It is trivial that (6.4) holds for all $w \in W_C$, but
to prove it for $w \in \widetilde W_C$ we need to use Lemma 1.5 (applied to
$R^\vee$), just as in the previous proof. Indeed, by Lemma 1.5 the
expression (6.2) vanishes unless $s^2=1$. Therefore summing (6.2) over $\{
s \in \hat A_{\sc} \, | \, s^2 =1 \}$ yields the same result as summing
over $\{s \in \hat A_{\sc} \, | \, s^2 \in Z^\vee \}$. It follows that
(drop the subscript 0 from $D$)
$$
   \sum_{w \in W_D} \epsilon(w)\psi_{wD^\vee}(\l,x) $$
is equal to
$$
|Z^\vee|^{-1}\sum_{w \in \widetilde W_D} \epsilon(w) \psi_{wD^\vee}(\l,x), $$
and it is clear that this expression is multiplied by $\epsilon(w)$ if $D$
is replaced by $wD$ for $w \in \widetilde W_D$. This proves (6.4) (take
$D=w_0C$ for any $w_0 \in W(\tau,C,\l)$). 

The next two properties of $b_R(\tau,C;x,\l)$ are obvious.
$$
b_R(w\tau,wC;x,\l)=b_R(\tau,C;x,\l) \quad \text{for all $w \in W$}.
\tag{6.5} $$
$$
b_R(\tau,C;wx,w\l)=b_R(\tau,C;x,\l) \quad \text{for all $w \in W$}.
\tag{6.6} $$
Moreover it follows from Proposition A.5 that
$$
  b_R(\tau,C;x,\l)=0 \quad \text{unless $\l(x)\le 0$,} $$
and since $b_R(\tau,C;x,\l)$ depends only on the chamber $C_x$ containing
$x$ we find that
$$
b_R(\tau,C;x,\l)=0 \quad \text{unless $\l \le 0$ on $C_x$}. \tag{6.7} 
$$
There is one more elementary property of the constants:
$$
b_R(\tau,C;x,\l)=1 \quad \text{if $R$ is empty}. \tag{6.8}
$$

The last property we need requires a bit more work. Suppose that $\a \in R$
and put $Y=\ker(\a) \subset X$. Define $R_\a$,$R_\a^\vee$ as in the
discussion at the beginning of \S2. Recall that $R_\a$ generates $Y^*$ and
that $-1_Y \in W(R_\a)$. Let $s=s_\a \in W$ be the reflection in the root
$\a$. Assume further that $C$ is a chamber in $X$ such that $\a$ belongs to
the closure of $C^\vee$. Let $x,x'$ be $R^\vee$-regular elements 
in $X$ that lie in
adjacent chambers separated by the wall $Y$. We are going to derive a
formula for
$$
    b_R(\tau,C;x,\l)+b_R(\tau,C;x',\l) $$
in terms of the constants $b_{R_{\a}}$ associated to the root system $R_\a$.

To get a clean formula we need to use the constants for the root system
$R_\a$ to define constants $b^R_{R_\a}(\tau,C;y,\l)$ for
$R_\a^\vee$-regular $y \in Y$
and $\tau \in X^*$, $\l \in W \cdot \tau$, $C$ as before (subject to the
requirement that $\a$ belongs to the closure of $C^\vee$). Write $W_\a$ for
the Weyl group of $R_\a$. Then we define
$$
  b^R_{R_\a}(\tau,C;y,\l)=0 \quad \text{unless $\l \in W_\a W_C\cdot \tau$}. $$
If $\l$ does belong to $W_\a W_C\cdot \tau$, choose $\tau' \in W_C\cdot
\tau$ such that $\l \in W_\a \cdot \tau'$ and put
$$
    b^R_{R_\a}(\tau,C;y,\l)=b_{R_\a}(\tilde\tau',C_Y;y,\tilde\l). $$
Here $\tilde\tau'$,$\tilde\l \in Y^*$ denote the restrictions of $\tau',\l
\in X^*$ to $Y$ and $C_Y$ is the chamber in $Y$ determined by $C$ (thus
$C_Y=Y\cap \tilde C$, where $\tilde C$ is the unique chamber in $X$
relative to $R_\a$ that contains $C$). Since $\tau'$ is well-determined up
to an element of $W_\a\cap W_C=W_{C_Y}$, we see from (6.4) that
$b_{R_\a}(\tilde\tau',C_Y;y,\tilde\l)$ is independent of the choice of
$\tau'$ (the equality $W_\a\cap W_C=W_{C_Y}$ is a consequence of our
assumption that $\a$ belongs to the closure of $C^\vee$). It is easy to see
that for any $\l \in W\cdot \tau$ we have the formula
$$
  b^R_{R_\a}(\tau,C;y,\l)=(-1)^{q(R_\a)}\sum_{w \in W_\a\cap W(\tau,C,\l)}
\epsilon(y,wC_Y) \psi_{wC^\vee_Y}(\tilde \l,y). \tag{6.9} $$

Now we are ready to formulate the last property of our constants: for
$\a$,$s$,$C$,$x$,$x'$ as above
$$
b_R(\tau,C;x,\l)+b_R(\tau,C;x',\l)=b^R_{R_\a}(\tau,C;y,\l)+b^R_{R_\a}(\tau,C;y,s\l),
\tag{6.10} $$
where $y$ is the unique point of $Y$ lying on the line segment joining $x$
and $x'$. Note that by (6.6) the left-hand side of (6.10) can also be
written as
$$
b_R(\tau,C;x,\l)+b_R(\tau,C;x,s\l). $$

Let us now prove (6.10). Assume without loss of generality that $\a(x) > 0$
and $\a(x') < 0$. Then by Corollary A.3 the left-hand side of (6.10) is
equal to 
$$
(-1)^{q(R)}\sum_w \epsilon(x,wC)\cdot \psi_{(wC)^\vee_Y} (\tilde \l,y) \cdot
\eta(w) $$
where the sum is taken over the set of all $w \in W(\tau,C,\l)$ such that
the closure of $wC^\vee$ contains either $\a$ or $-\a$, and where
$\eta(w)=\pm 1$ is defined by
$$ 
  \eta(w)=\cases -1 &\text{if $w\bar C^\vee$ contains $\a$}, \\
1 &\text{if $w\bar C^\vee$ contains $-\a$}. \endcases $$
Note that $w\bar C^\vee$ contains $\a$ if and only if $w \in W_\a$, and
$w\bar C^\vee$ contains $-\a$ if and only if $sw \in W_\a$ (since
$s\a=-\a$). Looking back at the proof of Lemma 2.2, we see that if $w \in
W_\a$, then
$$
-(-1)^{q(R)}\epsilon(x,wC)=(-1)^{q(R_\a)}\epsilon(y,wC_Y). $$
Therefore the left-hand side of (6.10) is equal to the difference of
$$
(-1)^{q(R_\a)}\sum_{w \in W_\a \cap W(\tau,C,\l)}
\epsilon(y,wC_Y)\psi_{wC^\vee_Y}(\tilde \l,y)
$$
and 
$$
(-1)^{q(R_\a)}\sum_{w \in
sW_\a \cap W(\tau,C,\l)} \epsilon(y,wC_Y)\psi_{(wC)^\vee_Y}(\tilde \l,y).
$$
Replacing $w$ by $sw$ in the second sum (and using that $(swC)_Y=(wC)_Y$
and $\widetilde{s\l}=\tilde\l$), we see that the left-hand side of (6.10)
is equal to the right-hand side of (6.10), as we wished to show.

The constants $b_R$ are determined uniquely by
properties (6.4)-(6.8) and (6.10) (and the property that $b_R(\tau,C;x,\l)$
depends only on the chamber in which $x$ lies). To see this fix
$\tau$,$C$,$\l$ and regard $b_R(\tau,C;x,\l)$ as a function of $x$. If $R$
is empty, $b_R$ is given by (6.8). If it is non-empty, the value of
$b_R(\tau,C;x,\l)$ is given by (6.7) for $x$ in at least one chamber in
$X$. Therefore it is enough to know
$$
   b_R(\tau,C;x,\l)+b_R(\tau,C;x',\l) $$
whenever $x,x'$ lie in adjacent chambers. But by (6.5) (which we use to put
$C$ in good position relative to the wall separating $x,x'$) and (6.10) the
sum above can be written in terms of the constants for a root system of
lower rank, which we may assume have already been determined. 

It is known (see \cite{K,13.4}) that the discrete series constants
satisfy these same properties; therefore they are equal to the constants
$b(\tau,C;x,\l)$. Before making this statement more precise, we need to
change the indexing of our constants in order to facilitate comparison with
\cite{K}. We now fix a chamber $C$ such that $R_C=R_c$ and define 
constants $c(w,\l,\Delta^+)$ as follows. For a regular element $\l \in
X^*$, 
$w \in W$, and a system $\Delta^+$ of positive roots for $R$, we put
$$
    c(w,\l,\Delta^+):=b(\l,C;x,w\l), $$
where $x \in X$ is any $R^\vee$-regular 
element in the (positive) Weyl chamber in $X$ determined by
$\Delta^+$.  It follows from (6.4) that the right side in this definition
is independent of the choice of chamber $C$ such that $R_C=R_c$. 
The constants $c(w,\l,\Delta^+)$ are those in
\cite{K} (see properties (13.32)-(13.34) in \cite{K}).
In other words our $\tau$ corresponds to Knapp's $\l$, and our $\l$
corresponds to Knapp's $w\l$. \heading 7. Lefschetz formula on reductive Borel-Serre compactifications
\endheading 
\subheading{(7.1) The group $G$} Let $G$ be a connected
reductive group over $\Q$ and let $A_G$ denote the maximal $\Q$-split torus
in the center of $G$.  Choose a maximal $\Q$-split torus $A_0$ in $G$ and
let $M_0$ denote its centralizer, a Levi subgroup of $G$. Fix a parabolic
subgroup $P_0$ of $G$ over $\Q$ having $M_0$ as Levi component; then $P_0$
is a minimal parabolic subgroup of $G$ over $\Q$. See \S5 for notation and
terminology concerning parabolic and Levi subgroups. In particular for any
standard parabolic subgroup $P$ we write $M$ for the unique Levi component
of $P$ containing $M_0$ and $N$ for the unipotent radical of $P$; thus
$P=MN$.

\subheading{(7.2) The locally symmetric spaces $S_K$}
Let $K$ be a suitably small compact open subgroup of $G(\A_f)$.
Choose a maximal compact subgroup $K_G$ of $G(\r)$ in good
position relative to $M_0$, in the sense that the Cartan involution on $G$
associated to $K_G$ preserves $M_0$. For each standard parabolic subgroup
$P=MN$ we denote by $K_M$ the intersection of $K_G$ with $M(\r)$, a
maximal compact subgroup in $M(\r)$. We denote by $A_G(\r)^0$ the
identity component of the topological group $A_G(\r)$, and we denote by
$X_G$ the homogeneous space 
$$
G(\r)/(K_G\cdot A_G(\r)^0)$$ 
for $G(\r)$. We then denote by $S_K$ the space
$$ 
G(\Q)\backslash [(G(\A_f)/K)\times X_G]. $$

\subheading{(7.3) The local system $\bold E_K$ on $S_K$}
Let $E$ be an irreducible representation of the algebraic group $G$ on a
finite dimensional complex vector space. Then $E$ gives rise to a local
system $\bold E_K$ on $S_K$. By definition $\bold E_K$ is the sheaf of flat
sections of the flat vector bundle
$$
     G(\Q)\backslash [(G(\A_f)/K) \times X_G \times E]
$$ over $S_K$.

\subheading{(7.4) The Hecke correspondence $(c_1,c_2)$ on $S_K$}
Now fix an element $g$ in $G(\A_f)$, and let $K'$ be any
compact open subgroup of $G(\A_f)$ that is contained in 
$K\cap g^{-1}Kg$. We use $g$,$K'$ to form a Hecke
correspondence on $S_K$, as follows. 
The inclusion $K' \subset K$ induces a surjection 
$$c_1:S_{K'} \to S_{K}.$$
The inclusion $K' \subset g^{-1}Kg $ induces a surjection
$$S_{K'} \to S_{g^{-1}Kg}$$ 
which we compose with the canonical isomorphism (use the element $g$)
$$S_{g^{-1}Kg} \simeq S_K$$ to get a second surjection
$$
c_2:S_{K'} \to S_K. $$
There are canonical isomorphisms
$$
   c_1^*\bold E_K \simeq \bold E_{K'} \simeq c_2^*\bold E_K.
\tag{7.4.1} $$

\subheading{(7.5) The reductive Borel-Serre compactification $\Sbar_K$ of
$S_K$}
For any standard parabolic subgroup $P=MN$ we denote by
$S^P_K$ the space
$$ S^P_K:=
M(\Q)\backslash [(N(\A_f)\backslash G(\A_f)/K) \times X_M ], $$
where $X_M$ denotes the analog for $M$ of $X_G$, namely
$M(\r)/(K_M\cdot A_M(\r)^0)$.
Now we can make precise what it means for $K$ to be suitably small: we
require that for each standard $P$ the group $M(\Q)$ act freely on 
$$ 
    (N(\A_f)\backslash G(\A_f)/K)\times X_M. $$ 
The reductive Borel-Serre compactification (see \cite{GHM}) $\Sbar_K$ of
$S_K$ is a stratified space whose statra are indexed by standard
parabolic subgroups $P$ of $G$, the stratum indexed by $P$ being the
manifold $S^P_K$ described above. 

\subheading{(7.6) The weighted cohomology complex $\Ebar_K$ on $\Sbar_K$} 
Let $p$ be a weight profile (see \S1.1 of \cite{GHM}). Associated to the
representation $E$ and the weight profile $p$ is a constructible complex of
sheaves $\bold W^p\bold C^\bullet(\bold E)$ of complex vector spaces on
$\Sbar_K$ (see \S1.3 of \cite{GHM}). In this paper we will denote this
complex of sheaves by $\Ebar_K$; as the notation suggests, the restriction
of $\Ebar_K$ to $S_K$ may be identified with $\bold E_K$. 

\subheading{(7.7) The Hecke correspondence $(\bar c_1,\bar c_2)$ on
$\Sbar_K$} The maps
$$c_1,c_2:S_{K'} \to S_K $$ have unique continuous extensions
$$\bar c_1,\bar c_2:\Sbar_{K'} \to \Sbar_K.$$ These maps carry $S^P_{K'}$ 
onto $S^P_K$; in fact, representing points of $S^P_{K'}$ by
pairs $(x,x_\infty)$ where $x \in G(\A_f)$
and $x_\infty \in X_M$, we have that the image of the pair
$(x,x_\infty)$ under $\bar c_1$ (respectively, $\bar c_2$) is the point of
$S^P_K$ represented by $(x,x_\infty)$ (respectively,
$(xg^{-1},x_\infty)$). 

It follows from
the definition of weighted cohomology complexes that there are
canonical isomorphisms
$$
\bar c_1^*\Ebar_K \simeq \Ebar_{K'} \simeq \bar c_2^*\Ebar_K. 
\tag{7.7.1} $$ 
The Verdier dual of the weighted cohomology complex $\Ebar_K$ 
is (a shift of) the
weighted cohomology complex obtained from the contragredient of the
representation $E$ and the weight profile $\bar p$ dual to $p$ (see \S1.3
of \cite{GHM}). Thus, applying Verdier duality to (7.7.1), we find that
there are canonical isomorphisms 
$$
\bar c_1^!\Ebar_K \simeq \Ebar_{K'} \simeq \bar c_2^!\Ebar_K. 
\tag{7.7.2} $$ It follows that there is a canonical isomorphism 
$$
   \bar c_2^*\Ebar_K \to \bar c_1^!\Ebar_K, \tag{7.7.3} $$
obtained as the composition of the isomorphism (7.7.1) from $\bar
c_2^*\Ebar_K$ to $\Ebar_{K'}$ and the isomorphism (7.7.2) from $\Ebar_{K'}$
to $\bar c_1^!\Ebar_K$. Thus there is a canonical extension, namely the
morphism (7.7.3), of the Hecke correspondence $(\bar c_1,\bar c_2)$ to the
weighted cohomology complex $\Ebar_K$.

\subheading{(7.8) The goal} The canonical morphism (7.7.3) induces 
self-maps on hypercohomology groups 
$$
    H^i(\Sbar_K,\Ebar_K) \to  H^i(\Sbar_K,\Ebar_K).  \tag{7.8.1} $$
These maps are obtained as the composition of the canonical pullback map
$$
    H^i(\Sbar_K,\Ebar_K) \to H^i(\Sbar_{K'},\bar c_2^*\Ebar_K), $$
the map
$$    H^i(\Sbar_{K'},\bar c_2^*\Ebar_K) \to H^i(\Sbar_{K'},\bar
c_1^!\Ebar_K)
$$ induced by (7.7.3), and the canonical proper pushforward map
$$   H^i(\Sbar_{K'},\bar c_1^!\Ebar_K) \to H^i(\Sbar_K,\Ebar_K). $$
The Lefschetz fixed point formula is a formula for the alternating sum of
the traces of the self-maps (7.8.1). An explicit version of the Lefschetz 
formula
(for the case at hand) is given in the theorem on page 474 
of \cite{GM}; our goal here
is to rewrite that formula in terms of stable virtual
characters on the group $G(\r)$, using the results in \S5 of this paper.

\subheading{(7.9) Fixed points} 
First we need to determine the fixed points of the correspondence. Of
course a fixed point is an element $x$ of $\Sbar_{K'}$ such that
$\bar c_1(x)=\bar c_2(x)$. Let us fix a standard parabolic subgroup $P=MN$ and
determine the fixed points of the correspondence that lie in the
subset $S^P_{K'}$ of $\Sbar_{K'}$. The group $P(\A_f)$ acts on
$G(\A_f)/K'$ with finitely many orbits. Choose a set of
representatives $x_0 \in G(\A_f)$ for these orbits and put
$$ \align
K'_P(x_0)&=P(\A_f) \cap x_0 K'x_0^{-1} \\
K'_M(x_0)&=\text{image of $K'_P(x_0)$ in $M(\A_f)$.} \endalign $$
Then $S^P_{K'}$ is the disjoint union of the subsets
$$
   S^P_{K'}(x_0):= M(\Q) \backslash [(M(\A_f)/K'_M(x_0)) \times X_M],
$$
the disjoint union being indexed by the set of representatives $x_0$
chosen above. 

A pair $(y,y_\infty) \in M(\A_f) \times X_M$ represents a fixed point
in $S^P_{K'}(x_0)$ of our Hecke correspondence if and only if there
exists $\gamma \in M(\Q)$ such that \roster \item $\gamma y_\infty =
y_\infty$, and \item there exists $n \in N(\A_f)$ such that
$y^{-1}n\gamma y \in x_0Kgx_0^{-1}$. \endroster The conjugacy class of
$\gamma$ in $M(\Q)$ depends only on the fixed point we started with.
Now fix an element $\gamma \in M(\Q)$ and denote by
$\Fix(P,x_0,\gamma)$ the subset of $S^P_{K'}(x_0)$ consisting of all fixed
points of our correspondence for which the associated conjugacy class
in $M(\Q)$ is equal to that of $\gamma$. The discussion above shows
that $\Fix(P,x_0,\gamma)$ is equal to $$ M_\gamma(\Q) \backslash
(Y^\infty \times Y_\infty) $$ where $M_\gamma$ denotes the centralizer
of $\gamma$ in $M$, $Y^\infty$ denotes the subset of
$M(\A_f)/K'_M(x_0)$ consisting of elements in that set represented by
elements $y \in M(\A_f)$ such that $y^{-1}\gamma y $ belongs to the
image in $M(\A_f)$ of $P(\A_f) \cap x_0Kgx_0^{-1}$, and $Y_\infty$
denotes the set of fixed points of $\gamma$ in $X_M$.

The group $M_\gamma(\Q)$ acts freely on $Y^\infty \times Y_\infty$.
We write $I$ for the identity component of $M_\gamma$. The
group $I(\A_f)$ acts on $Y^\infty$ with finitely many orbits.
The space $Y_\infty$ is empty unless $\gamma$ is conjugate in $M(\r)$
to an element of $K_M\cdot A_M(\r)^0$, in which case $I(\r)$ acts
transitively on $Y_\infty$ and in fact 
$$
  Y_\infty = I(\r)/(K_I \cdot A_M(\r)^0) $$
for some maximal compact subgroup $K_I$ of $I(\r)$. 
In line with our usual notational conventions we write $X_I$ for the
homogeneous space
$$ I(\r)/(K_I \cdot A_I(\r)^0); $$
note that $Y_\infty$ maps onto $X_I$ and is in fact a principal fiber
bundle over that space for the vector group 
$$ A_I(\r)^0/A_M(\r)^0. $$

\subheading{(7.10) Euler characteristic of $S_K$}
We need to recall Harder's formula (see [H]) for
the Euler characteristic of the space $S_K$ (we should also note that this
Euler characteristic coincides with the Euler characteristic with compact
support of $S_K$). Harder's formula
involves several ingredients, which we now explain. Let us choose a
Haar measure $dg_f$ on $G(\A_f)$. Then the Euler characteristic of
$S_K$ has the form 
$$
    \chi(G)\cdot \vl(K)^{-1}. $$
Of course $\vl(K)$ denotes the measure of $K$ with respect to $dg_f$. 
The quantity $\chi(G)$ depends on $G$ and the Haar measure $dg_f$, but
not on $K$. Moreover $\chi(G)$ is 0 unless the group $G$ has a
maximal torus $T$ over $\r$ such that $T/A_G$ is anisotropic over $\r$. 
Assume now that this
condition is satisfied. Let $\Cal D(G)$ denote the finite set
$$
    \Cal D(G):=\ker[H^1(\r,T)\to H^1(\r,G)] $$
(as usual we write $H^1(\r,G)$ as an abbreviation for
$H^1(\Gal(\C/\r),G(\C))$). Since $G/A_G$ has an anisotropic maximal
torus over $\r$, there is an inner form $\overline G$ of $G$ over $\r$
such that $\overline G/A_G$ is anisotropic over $\r$. We pick a Haar
measure $dg_\infty$ on $G(\r)$ and transport it to the inner form
$\overline G(\r)$ in the usual way, by identifying the space of invariant
top degree differential forms on $G$ with the analogous space for
$\overline G$ (this identification is defined over $\r$ since $\overline G$
is an \it inner \rm form of $G$). 
We define $q(G)$ to be half the
real dimension of the symmetric space associated to the real points of
the adjoint group of $G$. Then $\chi(G)$ is equal to
$$
   (-1)^{q(G)}\vl(G(\Q)A_G(\r)^0\backslash
G(\A))\vl(A_G(\r)^0\backslash \overline G(\r))^{-1}|\Cal D(G)|. $$
Of course we use $dg_f$ and $dg_\infty$ to get a Haar measure on
$G(\A)$; note that $\chi(G)$ is independent of the choice of Haar
measures on $G(\r)$ and $A_G(\r)^0$. 

\subheading{(7.11) Euler characteristic with compact support of $\Fix(P,x_0,\gamma)$}
Assume that $Y_\infty$ is non-empty. 
The Euler characteristic with compact support 
of the space $\Fix(P,x_0,\gamma)$ (which is one
of the ingredients in the Lefschetz formula) is equal to
$$
|M_\gamma(\Q)/I(\Q)|^{-1} $$
times the Euler characteristic with compact support of
$$
   I(\Q)\backslash (Y^{\infty} \times Y_{\infty}), $$
and this latter Euler characteristic with compact support is equal to
$(-1)^{\dim(A_I/A_M)}$ times that of
$$
   I(\Q)\backslash (Y^{\infty} \times X_I), $$
since the natural surjection
$$
    I(\Q)\backslash (Y^{\infty} \times Y_{\infty}) \to I(\Q)\backslash
(Y^{\infty} \times X_I) $$
is a principal fiber bundle under the vector group 
$$ 
    A_I(\r)^0/A_M(\r)^0. $$
It follows from Harder's theorem (see (7.10)) that the Euler characteristic
with compact support of $\Fix(P,x_0,\gamma)$ is equal to
$$
(-1)^{\dim(A_I/A_M)} |M_\gamma(\Q)/I(\Q)|^{-1} \cdot \chi(I) \cdot \sum_y \vl(I(\A_f) \cap
yK'_M(x_0)y^{-1})^{-1}, $$
where the index set for the sum is the subset of
$$
   I(\A_f)\backslash M(\A_f)/K'_M(x_0) $$
consisting of elements that can be represented by an element $y \in
M(\A_f)$ such that $y^{-1}\gamma y$ belongs to the image in $M(\A_f)$ of
$P(\A_f) \cap x_0Kgx_0^{-1}$. Of course we have chosen a Haar measure
$di_f$ on $I(\A_f)$. Let us fix a Haar measure $dm$ on $M(\A_f)$ as well.
Define a locally constant compactly supported function $f_{P,x_0}$ on
$M(\A_f)$ as follows: $f_{P,x_0}$ is $\vl(K'_M(x_0))^{-1}$ times the
characteristic function of the image in $M(\A_f)$ of $P(\A_f) \cap
x_0Kgx_0^{-1}$. For any locally constant compactly supported function $f$
on $M(\A_f)$ write $O_\gamma(f)$ for the orbital integral
$$
    \int_{I(\A_f)\backslash M(\A_f)} f(m^{-1}\gamma m) \, dm/di_f. $$
Then
$$
  O_\gamma(f_{P,x_0})=\sum_y \vl(I(\A_f) \cap yK'_M(x_0)y^{-1})^{-1}
$$
with the same index set as above. Therefore the Euler characteristic with
compact support of
$\Fix(P,x_0,\gamma)$ is equal to
$$
(-1)^{\dim(A_I/A_M)} |M_\gamma(\Q)/I(\Q)|^{-1} \cdot \chi(I) \cdot O_\gamma(f_{P,x_0}). 
\tag{7.11.1} $$

\subheading{(7.12) Lefschetz formula (qualitative version)} We now need to
recall the general form taken by the Lefschetz formula in \cite{GM}. The
formula is a sum of contributions, one for each connected component $C$ of the
fixed point set of the Hecke correspondence. We further decompose each such
connected component into locally closed pieces
$$
    C_P:=C \cap S^P_{K'}. $$
There are two natural ways to break up the contribution of $C$ to the
Lefschetz formula as a sum of contributions from the pieces $C_P$. In
\cite{GM} one of these two ways was chosen; it leads to the version of 
the Lefschetz formula given in that paper. However it is the other version 
that we are using here. 

This alternative version differs in two respects from the one chosen in
\cite{GM}. The first is that it involves the Euler characteristic with
compact support of $\Fix(P,x_0,\g)$ (rather than its Euler characteristic).
The second is that neutral directions are treated as being contracting
(rather than expanding); this change affects the definition of the set
$I(\g)$ appearing in (7.14), as we explain in more detail when we make the
definition. 

The
subset $\Fix(P,x_0,\gamma)$ of the fixed point set is a disjoint union of
certain sets of the form $C_P$, and from \cite{GM} we see
that the total contribution of $\Fix(P,x_0,\gamma)$ to the Lefschetz
formula is given by the product of three factors:
\roster \item the Euler characteristic with compact support 
of $\Fix(P,x_0,\gamma)$, 
\item the ramification index
$$
   r(x_0) := [N(\A_f) \cap x_0Kx_0^{-1}:N(\A_f) \cap x_0K'x_0^{-1}] $$
of the map $\bar c_1$ at any point in $S^P_{K'}(x_0)$,
\item a factor $L_P(\gamma)$ that depends only on the $G(\r)$-conjugacy
class of the pair $(P,\gamma)$ (and, of course, the representation $E$ and
the weight profile $p$ as well). \endroster

We will review the precise form of the factor $L_P(\gamma)$ later. All that
matters for the moment is the property stated in \therosteritem3. The
discussion above shows that the Lefschetz formula (for the alternating sum
of the traces of the self-maps (7.8.1)) is given by the following sum
$$
   \sum_P \sum_\gamma (-1)^{\dim(A_I/A_M)} \cdot |M_\gamma(\Q)/I(\Q)|^{-1} \cdot \chi(I) \cdot
L_P(\gamma) \cdot O_\gamma(f_P), \tag{7.12.1} $$
where $f_P$ is the locally constant compactly supported function on
$M(\A_f)$ defined by
$$
f_P:=\sum_{x_0} r(x_0)f_{P,x_0}. \tag{7.12.2} $$
In the sum defining $f_P$, the index $x_0$ runs over a set of
representatives for the orbits of $P(\A_f)$ on $G(\A_f)/K'$, as before.
In the first sum in (7.12.1) $P$ runs through the standard parabolic
subgroups of $G$, and in the second sum $\gamma$ runs through the set of
$M(\Q)$-conjugacy classes of elements $\gamma \in M(\Q)$ such that the fixed
point set of $\gamma$ in $X_M$ is non-empty.

\subheading{(7.13) Some familiar harmonic analysis} Let $P=MN$ be a
standard parabolic subgroup of $G$. In (7.11) we fixed a Haar measure $dm$
on $M(\A_f)$. Now we fix a Haar measure $dg$ on $G(\A_f)$ as well. Pick a
compact open subgroup $K_0$ of $G(\A_f)$ such that
$$
    G(\A_f)=P(\A_f)K_0. $$
Choose Haar measures $dn$ on $N(\A_f)$ and $dk$ on $K_0$ so that the usual
integration formula holds:
$$
 \int_{G(\A_f)} f(g) \, dg= \int_{M(\A_f)}\int_{N(\A_f)}\int_{K_0} f(mnk)\,
dk\,dn\,dm  \tag{7.13.1} $$
for any $f$ in $C^\infty_c(G(\A_f))$, the space of all locally constant
compactly supported functions on $G(\A_f)$. Let $\delta_{P(\A_f)}$ denote
the modulus function on $P(\A_f)$; thus, for $x \in P(\A_f)$ we have
$$
   \delta_{P(\A_f)}(x):=|\det(\Ad(x);\Lie(N)\otimes\A_f)|_{\A_f}, $$
where $|\cdot|_{\A_f}$ is the normalized absolute value on $\A_f^\times$.
Given $f \in C^\infty_c(G(\A_f))$ we define a function $f_M \in
C^\infty_c(M(\A_f))$ in the usual way, by putting
$$
   f_M(m):=\delta_{P(\A_f)}^{-1/2}(m)\int_{N(\A_f)}\int_{K_0}
f(k^{-1}nmk)\,dk\,dn. \tag{7.13.2} $$
The function $f_M$ depends on $P$ and even on $K_0$, but its orbital
integrals do not. It is worth noting that its orbital integrals are one of
the ingredients in Arthur's trace formula \cite{A}.

We are interested in a particular function $f^\infty \in
C^\infty_c(G(\A_f))$, namely the Hecke operator associated to the Hecke
correspondence in (7.4). Explicitly, $f^\infty$ is by definition
$\vl_{dg}(K')^{-1}$ times the characteristic function of the coset $Kg$
(with $K$,$g$,$K'$ as in (7.4)). Applying the discussion above to
$f^\infty$, we get $f^\infty_M \in C^\infty_c(M(\A_f))$ (defined by
(7.13.2), with $f$ replaced by $f^\infty$). In (7.12.2) we defined a
function $f_P \in C^\infty_c(M(\A_f))$. 
\proclaim{Lemma 7.13.A} The functions $f^\infty_M$ and
$\delta^{-1/2}_{P(\A_f)} \cdot f_P$ 
have the same orbital integrals. \endproclaim
It is equivalent to prove that $\delta^{1/2}_{P(\A_f)}\cdot f^\infty_M$ and
$f_P$ have the same orbital integrals. Note that although $f_P$ depends on
the choice of representatives $x_0$ for the double cosets
$$  
P(\A_f)\backslash G(\A_f) /K', $$
its orbital integrals do not. Suppose that we replace $K'$ by a compact
open subgroup $K''$ of $G(\A_f)$ contained in $K'$. Then $f^\infty$ is
multiplied by the index $[K':K'']$, as is $f^\infty_M$. An easy calculation
shows that $f_P$ is also multiplied by $[K':K'']$, as long as we take
representatives for
$$
P(\A_f)\backslash G(\A_f) /K'' $$
of the form $x_0k'$, where $x_0$ is one of our previous representatives and
$k' \in K'$. Therefore, by shrinking $K'$ if necessary, it is enough to
prove the lemma in the case that $K'$ is contained in $K_0$. Then we have
$$
   \delta^{1/2}_{P(\A_f)}(m)f^\infty_M(m)=\sum_{x \in K_0/K'}
\vl_{dk}(K')\int_{N(\A_f)} f^\infty(x^{-1}nmx) \, dn. $$ Thus
$\delta^{1/2}_{P(\A_f)} \cdot f^\infty_M$ has the same orbital integrals as
the function of $m \in M(\A_f)$ given by
$$
\sum_{x \in (P(\A_f)\cap K_0)\backslash K_0/K'} a(x) \int_{N(\A_f)}
f^\infty (x^{-1}nmx) \, dn, $$ where
$$
 \align    a(x)=&\vl_{dk}(K')\cdot [P(\A_f)\cap K_0:P(\A_f) \cap xK'x^{-1}]
\\             =&\vl_{dg}(K')\cdot \vl_{dm\,dn}(P(\A_f)\cap xK'x^{-1})^{-1}. 
\endalign $$
Here we used that
$$
\align
   \vl_{dm\,dn}(P(\A_f)\cap K_0)=&\vl_{dg}(K_0)\cdot \vl_{dk}(K_0)^{-1} \\
                                =&\vl_{dg}(K')\cdot \vl_{dk}(K')^{-1}, 
\endalign $$
which follows from (7.13.1), applied to the characteristic function of
$K_0$. From the equality 
$$
      G(\A_f)=P(\A_f)K_0 $$
it follows that
$$
   (P(\A_f)\cap K_0) \backslash K_0 /K' \simeq P(\A_f) \backslash G(\A_f)
/K'. $$ Thus the elements $x$ used here can serve as the elements $x_0$
used to define $f_P$. Moreover it is easy to see that
$$
   \int_{N(\A_f)} f^\infty (x^{-1}nmx) \, dn $$
is 0 unless $m$ lies in the image in $M(\A_f)$ of 
$$
    P(\A_f) \cap xKgx^{-1}, $$
in which case it equals
$$
   \vl_{dn}(N(\A_f) \cap xKx^{-1}) \cdot \vl_{dg}(K')^{-1}; $$
this shows that (with $x_0 = x$)
$$
   a(x)\int_{N(\A_f)} f^\infty(x^{-1}nmx) \, dn = r(x_0)f_{P,x_0}(m). $$
The proof of the lemma is now complete.

\subheading{(7.14) Manipulation of the Lefschetz formula} 
From Lemma 7.13.A we see
that the orbital integral $O_\gamma(f_P)$ appearing in (7.12.1) can be
rewritten as
$$
   O_\gamma(f_P)=\delta^{1/2}_{P(\A_f)}(\gamma)\cdot O_\gamma(f^\infty_M),
$$ with $f^\infty_M$ as in (7.13). Moreover, since $\gamma \in M(\Q)$ the
product formula shows that
$$
  \delta_{P(\A_f)}(\gamma)=\delta^{-1}_{P(\r)}(\gamma), $$
where for $x \in P(\r)$ we put
$$
   \delta_{P(\r)}(x)=|\det(\Ad(x);\Lie(N)\otimes \r)|. $$
Therefore the Lefschetz formula (7.12.1) can be rewritten as
$$
   \sum_P \sum_{\gamma} |M_\gamma(\Q)/I(\Q)|^{-1} \cdot \chi(I) \cdot
O_{\gamma} (f^\infty_M)\cdot (-1)^{\dim(A_I/A_M)}
\cdot \delta^{-1/2}_{P(\r)}(\gamma) \cdot
L_P(\gamma), \tag{7.14.1} $$
with the two index sets for the sum the same as in (7.12.1).

The factors $|M_\gamma(\Q)/I(\Q)|^{-1}$ and $\chi(I)$ depend only on
$M_\gamma$. The factor $O_\gamma(f^\infty_M)$ depends only on the
$G(\Q)$-conjugacy class of the pair $(M,\gamma)$. Therefore we  could
rewrite (7.14.1) by grouping together terms according to the
$G(\Q)$-conjugacy class of $(M,\gamma)$. However, there are still more
terms that can be grouped together, and it is this phenomenon that we must
study next.

The elements $\gamma$ appearing in (7.14.1) are semisimple, since they are
required to fix some point in $X_M$. We now fix a semisimple element
$\gamma \in G(\Q)$ and consider the set $\Cal M_\gamma$ of Levi subgroups
of $G$ such that $\gamma \in M(\Q)$. Recall that $M$ coincides with the
centralizer $\Cent_G(A_M)$ of $A_M$ in $G$ ($A_M$ denotes the maximal
$\Q$-split torus in the center of $M$, in line with the notational
conventions established in (7.1)). Thus a Levi subgroup $M$ of $G$ contains
$\gamma$ if and only if $A_M$ is contained in $G_\gamma^0$, the identity
component of the centralizer of $\gamma$ in $G$. We are now going to define
a map $M \mapsto M^*$ from $\Cal M_\gamma$ to itself. Let $M \in \Cal
M_\gamma$. Put $I:=M^0_\gamma$ and note that $A_I$ contains $A_M$. Then put
$$
   M^*:=\Cent_G(A_I)=\Cent_M(A_I). $$
Clearly $M^*$ is a member of $\Cal M_\gamma$ and $M^*$ is contained in $M$.

\proclaim{Lemma 7.14.A} The following statements hold.
\roster \item $\gamma \in I \subset M^*$. \item $I = (M^*)^0_\gamma$. \item
$A_{M^*}=A_I$. \item $M^{**}=M^*$. \item Suppose that $M^*=M$ and that
$M_1$ is a Levi subgroup of $G$ containing $M$. Then $M^*_1=M$ if and only
if $(M_1)^0_\gamma=M^0_\gamma$. \endroster \endproclaim
The first statement is obvious. The second statement follows from the first
and the fact that $M^*$ is contained in $M$. As for the third statement,
the inclusion $A_{M^*} \supset A_I$ is trivial and the opposite inclusion
follows from the obvious inclusion $A_{M^*} \subset A_{(M^*)^0_\gamma}$,
since $(M^*)^0_\gamma=I$ by the second statement. The fourth statement
follows immediately from the second. Finally we verify the fifth statement.
If $M^*_1=M$, apply the second statement to $M_1$ to see that
$(M_1)^0_\gamma=M^0_\gamma$. Conversely, suppose that
$(M_1)^0_\gamma=M^0_\gamma$. Then, applying the third statement to both
$M_1$ and $M$, we see that $A_{M^*_1}=A_M$, which in turn implies that
$M^*_1=M$.   

The first and second statements of the lemma together allow us to apply the
usual descent theory for orbital integrals. Put (as usual) 
$$
   D^M_{M^*}(\gamma):=\det(1-\Ad(\gamma);\Lie(M)/\Lie(M^*)) \in \Q. $$
Since $M^0_\gamma=(M^*)^0_\gamma$, it follows that
$$
    D^M_{M^*}(\gamma)\ne 0, $$
and from descent theory we have the equality
$$
    O_\gamma(f^\infty_{M^*})=|D^M_{M^*}(\gamma)|^{1/2}_{\A_f}\cdot
O_\gamma(f^\infty_M). \tag{7.14.2} $$
Furthermore the product formula implies that
$$
   |D^M_{M^*}(\gamma)|^{1/2}_{\A_f}=|D^M_{M^*}(\gamma)|^{-1/2}_\r. $$
 
The Lefschetz formula (7.14.1) can now be rewritten as
$$ \align
   \sum_{(P,M,\gamma)} |M_\gamma(\Q)/I(\Q)|^{-1} &\cdot \chi(I) \cdot
O_\gamma(f^\infty_{M^*}) \cdot (-1)^{\dim(A_I/A_M)} \\ &\cdot
|D^M_{M^*}(\gamma)|^{1/2}_\r \cdot
\delta^{-1/2}_{P(\r)}(\gamma) 	
\cdot L_P(\gamma), \tag{7.14.3} \endalign $$
where the sum runs through a set of representatives for the
$G(\Q)$-conjugacy classes of triples $(P,M,\gamma)$ consisting of a
parabolic subgroup $P$ of $G$, a Levi factor $M$ of $P$ and an element
$\gamma \in M(\Q)$, satisfying the condition that the fixed point set
$X_M^\gamma$ of $\gamma$ on $X_M$ be nonempty. As usual we write $I$ for
$M^0_\gamma$. We may impose the additional condition that the real group
$(I/A_I)_\r$ contain some anisotropic maximal $\r$-torus, since otherwise
$\chi(I)=0$ (see (7.10)).

For any triple $(P,M,\gamma)$ satisfying these conditions the real group
$(M^*/A_{M^*})_\r$ contains some anisotropic maximal $\r$-torus, and
moreover $\gamma$ is elliptic in $M^*(\r)$ (in other words, is contained in
some maximal $\r$-torus in $M^*$ that is anisotropic modulo $A_{M^*}$).
Indeed, by assumption there exists a maximal torus $T$ in $I$ over $\r$
such that $T/A_I$ is anisotropic over $\r$. But $T$ is also a maximal torus
in $M^*$ (use the second statement in Lemma 7.14.A), and $A_I=A_{M^*}$ (use
the third statement in Lemma 7.14.A), so that $T/A_{M^*}$ is an anisotropic
maximal torus in $(M^*/A_{M^*})_\r$ containing (the image of) $\gamma$. 

Therefore the Lefschetz formula (7.14.3) can be rewritten as
$$ \align
\sum_{(M,\gamma)} \chi(I) \cdot O_\gamma(f^\infty_M) \cdot \sum_{(P_1,M_1)}
&|(M_1)_\gamma(\Q)/I(\Q)|^{-1} \cdot |D^{M_1}_M(\gamma)|^{1/2}_\r \\ &\cdot
(-1)^{\dim(A_M/A_{M_1})} \cdot
\delta^{-1/2}_{P_1(\r)}(\gamma) \cdot L_{P_1}(\gamma). \tag{7.14.4}
\endalign $$
In (7.14.4) $I$ again denotes $M^0_\gamma$, and the index sets for the sums
are as follows. The first sum runs over a set of representatives for the
$G(\Q)$-conjugacy classes of pairs $(M,\gamma)$, where $M$ is a Levi
subgroup of $G$ and $\gamma \in M(\Q)$, satisfying the conditions that
$(M/A_M)_\r$ contain some anisotropic maximal $\r$-torus and that $\gamma$
be elliptic in $M(\r)$. Write $N_G(M)$ for the normalizer of $M$ in $G$.
Then the second sum in (7.14.4) runs over a set of representatives for the
$N_G(M)(\Q)\cap G_\gamma(\Q)$-conjugacy classes of pairs $(P_1,M_1)$
consisting of a parabolic subgroup $P_1$ in $G$ and a Levi factor $M_1$ of
$P_1$, satisfying the three conditions \roster \item $M_1 \supset M$ (which
implies that $\gamma \in M_1(\Q)$), \item $M^*_1 = M$, \item $X^\gamma_{M_1}$
is non-empty. \endroster

Consider two pairs $(M,\gamma)$ and $(P_1,M_1)$ appearing in (7.14.4), and
let $I=M^0_\gamma$, as before. Our assumptions on $(M,\gamma)$ imply that
$(I/A_I)_\r$ contains some anisotropic maximal $\r$-torus and that
$A_I=A_M$ (note in passing that this implies that $M^*=M$). Let
$n(P_1,M_1)$ denote the number of elements in the conjugacy class of
$(P_1,M_1)$ under $N_G(M)(\Q)\cap G_\gamma(\Q)$. We are now going to check
that
$$
    |(M_1)_\gamma(\Q)/I(\Q)| \cdot n(P_1,M_1)=|(N_G(M)(\Q)\cap
G_\gamma(\Q))/I(\Q)|. \tag{7.14.5} $$
Indeed, it follows from the equality $A_M=A_I$ that
$$ \align
N_G(M)(\Q)&=N_G(A_M)(\Q) \\
&=N_G(A_I)(\Q). \endalign $$
The stabilizer of $(P_1,M_1)$ in $G(\Q)$ is $M_1(\Q)$, and hence its
stabilizer in $N_G(M)(\Q) \cap G_\gamma(\Q)$ is $N_{M_1}(M)(\Q) \cap
(M_1)_\gamma(\Q)$. Since any element of $(M_1)_\gamma(\Q)$ normalizes the
split component $A_I$ of $(M_1)^0_\gamma=I$, we see that
$$
     N_{M_1}(M)(\Q) \supset (M_1)_\gamma(\Q), $$
and therefore the stablizer of $(P_1,M_1)$ in $N_G(M)(\Q)\cap G_\gamma(\Q)$
is simply $(M_1)_\gamma(\Q)$, which proves (7.14.5). Therefore the
Lefschetz formula (7.14.4) can be rewritten as
$$
  \sum_{(M,\gamma)}|(N_G(M)(\Q) \cap G_\gamma(\Q))/I(\Q)|^{-1} \cdot
\chi(I) \cdot O_\gamma(f^\infty_M) \cdot (-1)^{\dim A_M/A_G} \cdot
L_M(\gamma) \tag{7.14.6} $$
where the index set for the sum is the same as that for the first sum in
(7.14.4), and where
$$
   L_M(\gamma):= \sum_{(P_1,M_1)} (-1)^{\dim A_{M_1}/A_G}
|D^{M_1}_M(\gamma)|^{1/2}_\r \cdot \delta^{-1/2}_{P_1(\r)}(\gamma) \cdot
L_{P_1}(\gamma), \tag{7.14.7}  $$
with the index set equal to the set of \it all \rm pairs $(P_1,M_1)$
satisfying \therosteritem1, \therosteritem2, \therosteritem3 (in other
words, we no longer divide by the action of $N_G(M)(\Q) \cap G_\gamma(\Q)$
on the set of pairs). 

Note that condition \therosteritem2 is equivalent to the condition that the
connected centralizers of $\g$ in $M$ and $M_1$ coincide (recall that
$M^*=M$ and apply Lemma 7.14.A(5)), and this in turn is equivalent to the
condition that $D^{M_1}_M(\g)$ be non-zero. Therefore condition
\therosteritem2 can be dropped without changing (7.14.7).

Our assumptions on $(M,\g)$ imply that there exists a maximal torus $T$ of
$G$ over $\r$ containing $\g$ such that $T/A_M$ is anisotropic over $\r$.
The element $\g$ can be written as
$$
\g = \exp(x) \cdot \g_1 $$
for unique elements $x \in \fA_M$ and $\g_1 \in T(\r)_1$, where $T(\r)_1$
denotes the maximal compact subgroup of $T(\r)$ (see \S5 for the definition
of $\fA_M$). It is easy to see that $M_1$ satisfies condition
\therosteritem3 if and only if $x$ belongs to the subspace $\fA_{M_1}$ of
$\fA_M$. We conclude that
$$
  L_M(\g)=\sum_Q (-1)^{\dim(A_L/A_G)} \cdot |D^L_M(\g)|^{1/2}_\r \cdot
\d^{-1/2}_{Q(\r)}(\g) \cdot L_Q(\g) \tag{7.14.8} $$
where the sum ranges over all parabolic subgroups $Q=LU$ containing $M$
such that $x$ lies in $\fA_L$ (here, as usual, $L$ denotes the unique Levi
component of $Q$ containing $M$ and $U$ denotes the unipotent radical of
$Q$).

At this point we need to be more precise about the complex numbers
$L_Q(\g)$, and in order to do so we must first discuss weight profiles.
As in \S5 let $\nu$ be an element in $(\fA_{P_0})^*$ whose restriction to
$\fA_G$ coincides with the character by which $A_G$ acts on $E$. Then $\nu$
determines a weight profile (see \S1.1 of \cite{GHM}) in the following way.
As in \S5 $\nu$ determines elements 
$$
  \nu_P \in \fA_P^* $$
for every $P \in \Cal F(M_0)$ and in particular for every standard parabolic
subgroup $P$. The weight profile we associate to $\nu$ is the one such that
the restriction to the stratum $S^P_K$ (for standard $P=MN$) of the $i$-th
cohomology sheaf of the weighted cohomology complex $\Ebar_K$ is the local
system associated to the finite dimensional representation
$$
    H^i(\Lie(N),E)_{\ge\nu_P} $$
of $M$ defined in \S5. See Proposition 17.2 of \cite{GHM} in order to see
how to describe this weight profile in the language of that paper; note
that in \cite{GHM} the subspace 
$$
   H^i(\Lie(N),E)_{\ge\nu_P} $$
is denoted by
$$
  H^i({\frak n}_P,E)_+. $$
 
Now we are in a position to give a precise formula for $L_Q(\g)$. Here, as
above, $Q=LU$ is a parabolic subgroup containing $M$ such that $x$ lies in
the subspace $\fA_L$ of $\fA_M$. This formula is essentially the one given
in \cite{GM,\S11}, although that paper only discusses two particular weight
profiles and uses the other natural way of breaking up the Lefschetz
formula as a sum over strata of fixed point components. Let
$\a_1,\dots,\a_n \in \fA_L^*$ be the simple roots of $A_L$ in $\Lie(U)$.
Let $I=\{1,\dots,n\}$ (we temporarily reuse the notation $I$ in this way
for the sake of compatibility with \cite{GM}).
Put
$$
     I(\g):=\{i \in I \,|\, \langle \a_i,x \rangle < 0 \} $$
(the set of ``expanding directions").
In \cite{GM} the set $I(\g)$ is defined instead by the condition 
$$
   \langle \a_i,x \rangle \le 0; $$
this, together with the use of Euler characteristics rather than Euler
characteristics with compact support, is the other natural way of breaking
up the Lefschetz formula.

Choose a Borel subgroup $B$ of $G$ over $\C$ containing $T$ and contained
in $Q$, and let $W'$ be the corresponding Kostant representatives for the
cosets $W_L\backslash W$. Let $\l_B,\rho_B \in X^*(T)_{\r}$ be as usual
(see \S5). 
Let $t_1,\dots,t_n \in \fA_L/\fA_G$ be the basis of $\fA_L/\fA_G$ dual to
the basis $\a_1,\dots,\a_n$ of $(\fA_L/\fA_G)^*$. For $w \in W'$ put
$$
  I(w)=\{ i \in I \,|\, \langle p_L(w(-\l_B-\rho_B)+\rho_B)+\nu_Q,t_i
\rangle > 0 \}. $$
In \cite{GM} the set $I(w)$ was defined using $\le$ rather than $>$, but this
was a misprint and should have been $\ge$. If we take $\nu=\nu_m$, our set
$I(w)$ still differs ($>$ rather than $\ge$) from the (corrected) one in
\cite{GM}: there are two middle weighted cohomology complexes, upper and
lower, and the formula in \cite{GM} arises from one of them while the
formula here arises from the other.

The complex number $L_Q(\g)$ is defined by 
$$
   (-1)^{|I(\g)|} \sum_w \epsilon(w) \cdot
\tr(\g^{-1};V^L_{w(\l_B+\rho_B)-\rho_B}) $$
where the index set for the sum is the set of $w \in W'$ such that
$I(w)=I(\g)$. Perhaps the following remarks will help the reader understand
the number $L_Q(\g)$. The local contribution to the Lefschetz fixed point
formula from a connected component of the fixed point set is given (see
\cite{GM2}) by the (weighted) cohomology (of a regular neighborhood of the
fixed point component) with supports which are compact in the expanding
directions and which are closed in the contracting directions. The piece of
weighted cohomology which is indexed by $w \in W'$ (via Kostant's theorem)
contributes to stalk cohomology with compact supports in the directions
$I(w)$ and with degree shift by $|I(w)|$. The Hecke correspondence is
expanding away from the stratum $S^Q_K$ in the directions $I(\g)$. Thus the
contributions to the Lefschetz formula occur when $I(w)=I(\g)$. 

As in \S5 let $C_Q$ be the (open) chamber in $\fA_L$
corresponding to $Q \in \Cal P(L)$; of course the image of the 
cone $C_Q$ in $\fA_L/\fA_G$ is generated
by $t_1,\dots,t_n$. Again as in \S5 let $\varphi_Q$ denote the function
$$
  \varphi_{C_Q}(\cdot,\cdot) $$
on $\fA_L \times \fA_L^*$ determined by the \it open \rm cone $C_Q$ (see
the last part of Appendix A). By Lemma A.7 and the analog of (A.2)
mentioned just before Lemma A.7 the value of
$\varphi_Q(-x,p_L(w(\l_B+\rho_B)-\rho_B)-\nu_Q)$ is given by
$$
   (-1)^{\dim(A_G)}\cdot (-1)^{\dim(A_L/A_G)-|I(\g)|} $$
if $I(w)=I(\g)$ and is 0 otherwise. Therefore
$$ \align
   L_Q(\g)=\sum_{w \in W'} \epsilon(w) \cdot (-1)^{\dim(A_L)} \cdot
&\varphi_Q(-x,p_L(w(\l_B+\rho_B)-\rho_B)-\nu_Q) \\ \cdot
&\tr(\g^{-1};V^L_{w(\l_B+\rho_B)-\rho_B}). \endalign $$

This expression for $L_Q(\g)$ coincides exactly with the one used to
define $L^\nu_Q(\g)$ in \S5. Moreover, comparing the definition of
$L^\nu_M(\g)$ given in \S5 with (7.14.8), we now see that the numbers
$L^\nu_M(\g)$ from \S5 and $L_M(\g)$ from this section are equal. Applying
Theorem 5.1 yields the following result (we now let $I$ once again denote
$M^0_\g$). 
\proclaim{Theorem 7.14.B} The Lefschetz formula for the alternating sum of
the traces of the self-maps (7.8.1) is given by
$$
   \sum_{(M,\g)} |(N_G(M)(\Q) \cap G_\g(\Q))/I(\Q)|^{-1} \cdot \chi(I)
\cdot O_\g(f^\infty_M) \cdot (-1)^{\dim(A_M/A_G)} \cdot
\Phi_M(\g,\Theta_{\nu}) $$ 
with $\Theta_{\nu}$ and $\Phi_M(\g,\Theta_{\nu})$ as in \S5, the index set
for the sum being the same as that for the first sum in (7.14.4).
\endproclaim 
\subheading{(7.15) $\Q$-equivalence} There are three special cases in which
the statement of Theorem 7.14.B can be simplified, but before we can do
this we need to make a couple of definitions. Let $\Theta$, $\Theta'$ be
stable virtual characters on~$G(\r)$. We say that $\Theta$ and $\Theta'$
are $\Q$-{\it equivalent\/} if they agree on $T_{\reg}(\r)$ for every
maximal $\r$-torus $T$ in~$G$ whose $\r$-split component is both defined and
split over~$\Q$. Note that $\Theta$ and $\Theta'$ are $\Q$-equivalent if
and only if the functions $\Phi_M(\cdot,\Theta)$ and
$\Phi_M(\cdot,\Theta')$ coincide for every Levi subgroup $M$ of~$G$
over~$\Q$ such that $(M/A_M)_\r$ contains an anisotropic maximal
$\r$-torus. Clearly the expression for the Lefschetz formula given in
Theorem 7.14.B remains valid when $\Theta_\nu$ is replaced by any stable
virtual character $\Theta'$ that is $\Q$-equivalent to $\Theta_\nu$. 

\subheading{(7.16) The orientation character $\chi_G$} We will also need
the following sign character
$$ \chi_G:G(\r) \to \{ \pm 1 \}. $$
Recall from 7.2 the real manifold
$$
X_G=G(\r)/(K_G\cdot A_G(\r)^0), $$
on which $G(\r)$ acts by left translations. Of course
$X_G$ is diffeomorphic to a Euclidean space and hence is orientable. 
For $g \in G(\r)$ we
define $\chi_G(g)$ to be $-1$ if $g$ reverses the orientation of $X_G$ and
$+1$ if $g$ preserves the orientation of $X_G$.

Let $P=MN$ be a parabolic subgroup of $G$ whose Levi component $M$ contains
$M_0$. Suppose that $M$ contains a maximal $\r$-torus $T$ such that $T/A_M$
is anisotropic over~$\r$. We claim that
$$
    \chi_G(\g)=\sgn(\det(\g;\Lie(N))) \tag{7.16.1} $$
for all $\g \in T(\r)$.

Indeed, since $P(\r)$ acts transitively on $X_G$, we see that as an
$M(\r)$-space $X_G$ is given by
$$
    X_G=\Lie(N) \times X_M \times (A_M(\r)^0/A_G(\r)^0). \tag{7.16.2} $$
It follows that
$$
   \chi_G(m)=\chi_M(m) \cdot \sgn(\det(m;\Lie(N))) \tag{7.16.3} $$
for all $m \in M(\r)$. Since $T(\r)/A_M(\r)$ is connected, and since
$A_M(\r)$ acts trivially on $X_M$, we see that $\chi_M$ is trivial on
$T(\r)$, so that (7.16.1) follows from (7.16.3).

\subheading{(7.17) Very positive $\nu$} Now suppose that $\nu$ is
sufficiently positive. Then the weighted cohomology complex $\Ebar_K$ is
equal to $Rj_!\bold E_K$, where $j$ denotes the inclusion of $S_K$ in
$\Sbar_K$, and the global weighted cohomology
is the cohomology with compact support of $S_K$ with coefficients in
$\bold E_K$. Moreover it is obvious from the definition of $\Theta_\nu$ that
$\Theta_\nu$ is equal to (the character of) the representation $E^*$
contragredient to~$E$ (for sufficiently positive $\nu$ the virtual modules
$E^\nu_P$ used to define $\Theta_\nu$ are trivial except when
$P=G$). Therefore in this case the expression for the Lefschetz formula
given by Theorem 7.14.B involves the character of the representation $E^*$,
as was also observed by J.~Franke \cite{F} and G.~Harder. 

\subheading{(7.18) Very negative $\nu$} Now suppose that $\nu$ is
sufficiently negative. Then the weighted cohomology complex $\Ebar_K$ is
equal to the full direct image $Rj_*\bold E_K$, and the global weighted
cohomology is the ordinary cohomology of $S_K$ with coefficients in $\bold
E_K$. 
Moreover we claim that $\Theta_\nu$
is $\Q$-equivalent to (the character of) the virtual
representation $(-1)^{\dim(X_G)}\chi_G \otimes E^*$ of $G(\r)$, 
where $X_G$ is the space defined in 7.2 and $\chi_G$ is the
orientation character defined in 7.16.

Let $T$ be a maximal $\r$-torus in~$G$ whose $\r$-split component $A$ is
defined and split over~$\Q$. Let $M$ be the centralizer of~$A$ in~$G$, a
Levi subgroup of~$G$ over~$\Q$. Replacing $T$ by a conjugate under $G(\Q)$,
we may assume that $M$ contains $M_0$. We must show that for all $\g \in
T(\r)$
$$
  \Phi_M(\g,\Theta_\nu)=(-1)^{\dim(X_G)} \cdot 
\chi_G(\g) \cdot \Phi_M(\g,E^*). \tag{7.18.1} $$
As usual we write 
$$
   \g =\exp(x)\cdot \g_1 $$
for unique elements $x \in \fA_M$ and $\g_1 \in T(\r)_1$.
By continuity it is enough to prove (7.18.1) when $\g$ is regular in
$T(\r)$ and $x$ is regular in $\fA_M$.

We use Theorem 5.1 to evaluate $\Phi_M(\g,\Theta_\nu)$. Since we are
assuming that $\nu$ is sufficiently negative, the factor
$$
 \varphi_Q(-x,p_L(w(\l_B+\rho_B)-\rho_B)-\nu_Q) $$
entering into the definition of $L^\nu_Q(\g)$ is $0$ unless $L=M$ and $x$
lies in the chamber $C_Q$ in $\fA_M$ determined by $Q$, in which case the
factor is $(-1)^{\dim(A_M)}$ (use Lemma A.7 and the analog for $\varphi_Q$ of
(A.2) to evaluate $\varphi_Q(-x,\mu)$ for very positive $\mu$). Therefore
Theorem 5.1 tells us that $\Phi_M(\g,\Theta_\nu)$ is equal to
$$ (-1)^{\dim(A_M/A_G)} \cdot \d_Q^{-1/2}(\g) $$
times
$$
    \sum_{w \in W'} \epsilon(w) \cdot
\tr(\g^{-1};V^M_{w(\l_B+\rho_B)-\rho_B}) \tag{7.18.2} $$
where $Q=MN$ is the unique parabolic subgroup with Levi component~$M$
such that $x$ lies in the chamber~$C_Q$. The number (7.18.2) is equal to
$$
    \tr(\g^{-1};E) \cdot \tr(\g;\wedge^\bullet(\Lie(N)), $$
where $\wedge^\bullet(\Lie(N))$ denotes the virtual $M$-module
$$
    \sum_{i \ge 0} (-1)^i \bigwedge^i(\Lie(N)). $$
Therefore, in order to prove the equality (7.18.1) it is enough to prove
that 
$$
(-1)^{\dim(A_M/A_G)} \cdot \d_Q^{-1/2}(\g) \cdot \prod_{\a} (1 - \a(\g))=
(-1)^{\dim(X_G)} \cdot \chi_G(\g) \cdot |D^G_M(\g)|^{1/2}, \tag{7.18.3} $$
where the product is taken over roots $\a$ of $T$ in $\Lie(N)$. The two
sides of (7.18.3) have the same absolute value, and therefore it is enough
to prove that the number
$$
  (-1)^{\dim(X_G)} \cdot
  (-1)^{\dim(A_M/A_G)} \cdot \chi_G(\g) \cdot \prod_{\a} (1-\a(\g))
\tag{7.18.4} $$
is positive. Using (7.16.1) we see that the number (7.18.4) has the same
sign as
$$
(-1)^{\dim(X_G)} \cdot 
(-1)^{\dim(A_M/A_G)} \cdot \prod_{\a} (\a^{-1}(\g)-1). $$
Complex conjugation preserves the set of roots $\a$ of $T$ in $\Lie(N)$. If
$\a$ is not real, then the product of $\a^{-1}(\g)-1$ and
$\bar\alpha^{-1}(\g)-1$ is positive. If $\a$ is real, then $\a^{-1}(\g)-1$
is negative (use that $x$ lies in the chamber $C_Q$). Therefore the sign of
the number (7.18.4) is equal to
$$
  (-1)^{\dim(X_G)} \cdot   (-1)^{\dim(A_M/A_G)} \cdot (-1)^{\dim(N)}. $$
It follows from (7.16.2) that
$$
\dim(X_G)=\dim(A_M/A_G) + \dim(N) + \dim(X_M). $$
Moreover $\dim(X_M)$ is even since $(M/A_M)_\r$ contains an anisotropic
maximal $\r$-torus.
Therefore the sign of the number (7.18.4) is indeed $+1$, as we wished to
show. 

\subheading{(7.19) Middle $\nu$} Now suppose that $\nu=\nu_m$. Then
$\Ebar_K$ is the upper middle weighted cohomology complex. Assume further
that $(G/A_G)_\r$ has an anisotropic maximal $\r$-torus. Then by Theorem
5.2 the virtual character $\Theta_\nu$ is $\Q$-equivalent to
$\Theta^{E^*}$. Thus the expression for the Lefschetz 
formula given in Theorem 7.14.B
essentially agrees with Arthur's formula \cite{A,Theorem 6.1} for the
alternating sum of the traces of a Hecke operator on the $L^2$-cohomology
groups of $S_K$ with coefficients in $\bold E_K$. 
Actually
Theorem 7.14.B appears at first to disagree with Arthur's formula, which
contains the factor $\Phi_M(\g,\Theta^E)$ rather than
$\Phi_M(\g,\Theta^{E^*})$, but by replacing $\g$ by $\g^{-1}$ in one of the
two formulas we come closer to 
agreement, since Arthur uses the usual \it left \rm
action of the Hecke algebra on the $L^2$-cohomology of $S_K$, while in this
paper we have (implicitly) used the \it right \rm action of the Hecke
algebra on the upper 
middle weighted cohomology of $\Sbar_K$. In other words it
is the \it right \rm action of $f^\infty$ on $H^i(\Sbar_K,\Ebar_K)$ that
coincides with the self-map defined by our Hecke correspondence $(\bar
c_1,\bar c_2)$; of course the right action of $f^\infty$ coincides with the
left action of the Hecke operator $f^\infty_r$ defined by
$$
    f^\infty_r(x)=f^\infty(x^{-1})  \qquad (x \in G(\A_f)). $$

We should also note that the index set for the sum in our formula is
somewhat different from the index set for the double sum in Arthur's
formula, and that our formula has the factor
$$
    |(N_G(M)(\Q) \cap G_\g(\Q) )/I(\Q) |^{-1} $$
while Arthur's has the factor
$$
    |W^M_0||W^G_0|^{-1}|\iota^M(\g)|^{-1}; $$
however it is easy to see that these are just two different ways of writing
the same thing. 

However, our formula disagrees with Arthur's in that we sum only over Levi
subgroups $M$ such that $M/A_M$ contains an \it anisotropic \rm maximal
torus over $\r$. In \cite{A} all Levi subgroups $M$ are allowed, although
the term indexed by $M$ vanishes (due to the factor $\Phi_M$) unless
$M/A_M$ contains an \it elliptic \rm maximal torus over $\r$. Thus Arthur's
formula may have more non-zero terms than ours (the split component of the
center of $M$ over $\r$ may be strictly bigger than $A_M$); however these
extra terms in Arthur's formula should not actually be there (the error
occurs in equation (4.1) of \cite{A}, which is only valid if the split
components of the center of $M$ over $\Q$ and $\r$ are the same).

\heading A. Functions associated to convex polyhedral cones \endheading
Let $X$ be a finite dimensional real vector space and let $C$ be a closed
convex polyhedral cone in $X$. Recall that this means that there exists a
finite subset $S$ of $X$ such that $C$ is equal to the set of non-negative
linear combinations of elements of $S$; equivalently, it means that there
exists a finite subset $T$ of the dual vector space $X^*$ such that $C$ is
equal to the intersection of the sets
$$
   \{x \in X\,|\, \l (x) \ge 0 \} $$
where $\l$ ranges through $T$. We denote by $C^*$ the closed convex
polyhedral cone in $X^*$ dual to $C$. Recall that $C^*$ is defined by
$$
   C^*:=\{\l \in X^* \, | \, \text{$\l(x) \ge 0$ for all $x \in C$}\} $$
and that the map
$$
   F \mapsto F^\perp:=C^* \cap \{ \l \in X^* \, | \, \text{$\l(x)=0$ for all
$x \in F$}\} $$
sets up a bijection between the set $\Cal F$ of (closed) faces of $C$ to the
set $\Cal F^*$ of (closed) faces of $C^*$. Moreover $C^{**}=C$ and
$F^{\perp\perp}=F$. The map $F \mapsto F^{\perp}$ is order-reversing (we
order faces by inclusion) and
$$
   \dim(F)+\dim(F^\perp)=\dim(X), $$
where $\dim(F)$ denotes the dimension of the linear span of $F$.

     Define an integer-valued function $\psi_C$ on $X\times X^*$ by 
$$
   \psi_C=\sum_{F \in \Cal F} (-1)^{\dim(F)} \xi_{(F^\perp)^* \times
F^*}, $$
where we have denoted by $\xi_{(F^\perp)^* \times
F^*} $ the characteristic function of the subset ${(F^\perp)^* \times
F^*}$ of $X \times X^*$ (since $F^\perp$ and $F$ are themselves closed
convex polyhedral cones, it makes sense to consider their dual cones
$(F^\perp)^*$ and $F^*$). Note that $(F^\perp)^*$ is equal to $C+\spn(F)$,
where $\spn(F)$ denotes the linear span of $F$. Clearly 
$$ \psi_{C^*}(\l,x)=(-1)^{\dim(X)}\psi_C(x,\l). \tag{A.1}$$

Let $X_1$ denote the linear span of $C$ and let $X_2$ denote the largest
linear subspace of $X$ contained in $C$.  In $X^*$ we have the
perpendicular subspaces
$$
  X_i^\perp=\{\l \in X^* \, |\, \text{$\l(x)=0$ for all $x \in X_i$ }\} \qquad
(i=1,2). $$  The cone $C$ gives rise to a cone $\tilde C$ in $X_1/X_2$ and
it is easy to see that $\psi_C(x,\l)$ is 0 unless $x \in X_1$ and $\l \in
X_2^{\perp}$, in which case
$$
   \psi_C(x,\l)=(-1)^{\dim(X_2)}\psi_{\tilde C}(\tilde x,\tilde \l)
\tag{A.2} $$ 
where $\tilde x,\tilde \l$ denote the images of $x,\l$ under the natural
surjections
$$
 \align
        &X_1 \to X_1/X_2 \\
        &X_2^\perp \to X_2^{\perp}/X_1^{\perp}=(X_1/X_2)^*
\endalign $$ respectively.

Suppose that $C$ is a simplicial cone in $X$. Then the pair $(X,C)$ is
isomorphic to the pair $(\r^n,(\r_{\ge0})^n)$, where $n=\dim X$ and
$\r_{\ge0}$ denotes the set of non-negative real numbers. We will now
calculate $\psi_C$ for the pair $(\r^n,(\r_{\ge0})^n)$. Of course we
identify $X^*$ with $\r^n$ as well. Let $x=(x_1,\dots,x_n) \in \r^n$ and
$\l = (\l_1,\dots,\l_n) \in \r^n$. Write $I$ for $\{1,\dots,n\}$ and then
put
$$ \align
   &I_x=\{i \in I \, | \, x_i \ge 0 \} \\
   &I_\l=\{i \in I \, | \, \l_i \ge 0 \}. \endalign $$

\proclaim{Lemma A.1} The number $\psi_C(x,\l)$ is $0$ unless the subsets
$I_x$ and $I_\l$ of $I$ are complementary, in which case
$$
    \psi_C(x,\l)=(-1)^{|I_\l|}, $$
where $|I_\l|$ denotes the cardinality of $I_\l$. \endproclaim

Indeed, the faces of $C$ are indexed by the subsets $J \subset I$, the face
corresponding to $J$ being 
$$ \{ x\in \r^n \, | \, \text{$x_i \ge 0$ for all $i \in J$ and $x_i=0$ for
all $i \in I\setminus J$ }\}. $$ Thus
$$
   \psi_C(x,\l)=\sum_J (-1)^{|J|} $$
where $J$ ranges through all subsets of $I$ such that
$$
     (I\setminus I_x) \subset J \subset I_\l, $$
and the lemma follows from the elementary fact that for any finite set $S$
the sum 
$$
     \sum_{T \subset S} (-1)^{|T|} $$
is 0 unless $S$ is empty, in which case it is 1.

In the next lemma and corollary it will be convenient to assume that
$$\dim(C)=\dim(C^*)=\dim(X). $$ 
We refer to the faces of $C$ having codimension 1 in $C$ as \it facets \rm
of $C$, and we refer to the linear spans of the facets of $C$ as the \it
walls \rm of $C$. The walls of $C$ form a finite set of hyperplanes in $X$.
We say that a point in $X$ is $C$-regular if it lies on no wall of $C$. We
refer to the connected components of the set of $C$-regular points in $X$
as \it chambers \rm in $X$. Of course the interior of $C$ is one of these
chambers. Similarly the walls of $C^*$ divide $X^*$ into chambers. The
value of $\psi_C(x,\l)$ on $C$-regular $x \in X$ and $C^*$-regular $\l \in
X^*$ depends only on the chambers in which $x,\l$ lie.

\proclaim{Lemma A.2} Let $x,x'$ be $C$-regular elements of $X$ and suppose
that there is exactly one wall $Y$ of $C$ separating $x$ and $x'$ \rom(in 
other
words $x,x'$ lie in adjacent chambers\rom).  Put $C_Y:=C \cap Y$, a facet of
$C$ which can also be regarded as a closed convex polyhedral cone in $Y$,
so that the function $\psi_{C_Y}$ on $Y \times Y^*$ is defined. Assume that
$x,C$ lie on the same side of $Y$. Then
$$
   \psi_C(x,\l)-\psi_C(x',\l)=\psi_{C_Y}(y,\l_Y), $$
where $\l_Y \in Y^*$ is the restriction of $\l$ to $Y$ and $y \in Y$ is the
unique point of $Y$ lying on the line segment joining $x$ and $x'$.
\endproclaim
 
Let $\a \in C^*$ be a non-zero element of the 1-dimensional face
$(C_Y)^\perp$ of $C^*$. Thus $Y$ is the kernel of the linear form $\a$, and
$\a$ is positive on $x$, negative on $x'$. Since $x,x'$ lie in adjacent
chambers, $\b(x)$ and $\b(x')$ have the same sign for any non-zero element
$\b$ of any 1-dimensional face of $C^*$ other than $(C_Y)^\perp$, and this
sign is the same as that of $\b(y)$. Therefore
$$
   \xi_{(F^\perp)^* \times
F^*}(x,\l) - \xi_{(F^\perp)^* \times
F^*} (x',\l) $$ 
is 0 unless $F$ is contained in $Y$ (equivalently, unless $F^\perp$
contains $\a$), in which case it is equal to
$$
   \xi_{G \times H}(y,\l_Y) $$
where $G=(F^\perp)^* \cap Y$ and $H$ is the image of $F^*$ under the
canonical surjection from $X^*$ to $Y^*$. But the faces of $C_Y$ are
precisely the faces $F$ of $C$ contained in $Y$, and, writing $F_Y$ to
indicate that we are regarding such a face $F$ as a face of $C_Y$, we have
$F_Y^*=H$ and 
$$
  \align
          (F_Y^\perp)^*&=C_Y+\spn(F) \\
                       &=(C+\spn(F)) \cap Y \\
                       &=(F^{\perp})^* \cap Y \\
                       &=G. 
\endalign $$
This proves the lemma.

Applying the lemma to $C^*$ and using (A.1), we obtain the following
result.
\proclaim{Corollary A.3} Let $\l,\l'$ be $C^*$-regular elements of $X$ and
suppose that there is exactly one wall $Z$ of $C^*$ separating $\l$ and
$\l'$. Let $\omega \in C$ be a non-zero element of the 1-dimensional face
$(C^*\cap Z)^\perp$ of $C$ corresponding
 to the facet $C^* \cap Z$ of $C^*$. Put
$\tilde X:=X/\r \omega$ and let $\tilde C$ denote the image of $C$ under
the canonical surjection $X \to \tilde X$; note that $Z$ is the dual vector
space to $\tilde X$ and that $\tilde C$ is the closed convex polyhedral
cone $(C^* \cap Z)^*$ in $\tilde X$ dual to the cone $C^* \cap Z$ in $Z$. In
particular the function $\psi_{\tilde C}$ on $\tilde X \times Z$ is
defined. Assume that $\l,C^*$ lie on the same side of $Z$ \rom(equivalently,
assume that $\l(\omega) > 0$ and $\l'(\omega) < 0$\rom). Then
$$
    \psi_C(x,\l)-\psi_C(x,\l')=-\psi_{\tilde C}(\tilde x,\tilde \l) $$
where $\tilde x$ denotes the image of $x$ under the canonical surjection
$X \to \tilde X$ and $\tilde\l$ denotes the unique point of $Z$ lying on
the line segment joining $\l$ and $\l'$. \endproclaim

Now we come to the main point of the appendix, which is to show that if $C$
decomposes as a union of cones, then $\psi_C$ decomposes accordingly. For
this we need a little preparation.  For any closed convex polyhedral cone
$C$ in $X$ we write $\inC$ for its relative interior (the interior of $C$
in $\spn(C)$). We say that a function $f:X \to \Z $ is \it conic \rm if it
can be written as a finite $\Z$-linear combination of characteristic
functions $\xi_C$ of closed convex polyhedral cones $C$ in $X$. Let $\l$
be a non-zero element of $X^*$, let $C$ be a closed convex polyhedral cone
in $X$, and put
$$\align
C_+&=\{x \in C \, | \, \l(x) \ge 0 \} \\
C_-&=\{x \in C \, | \, \l(x) \le 0 \} \\
C_0&=\{x \in C \, | \, \l(x) = 0 \}. \endalign $$
Then $C_+,C_-,C_0$ are closed convex polyhedral cones in $X$, and we have
the relation
$$
\xi_C + \xi_{C_0} - \xi_{C_+} -\xi_{C_-} = 0.  \tag{A.3} $$
It is not difficult to see that the abelian group $C(X)$ of conic functions
on $X$ is presented by the generators $\xi_C$ and the relations (A.3). 

The characteristic function $\xi_{\inC}$ of the relative interior $\inC$
of $C$ is a conic function on $X$. Indeed it can be expressed in the
following way in terms of our generators
$$
\xi_{\inC}=(-1)^{\dim(C)}\sum_F (-1)^{\dim(F)}\xi_F, \tag{A.4} $$
where the sum is taken over the set of faces $F$ of $C$. To prove this note
that the value of the right-hand side of (A.4) at a point $x \in X$ is 0
unless $x \in C$, in which case it is
$$
(-1)^{\dim(C)} \sum_F (-1)^{\dim(F)} $$
where $F$ now ranges through all faces of $C$ containing $F(x)$, the
smallest face of $C$ containing $x$. These faces correspond bijectively to
the faces of a new cone, namely $C+\spn(F(x))$. Therefore (A.4) follows
from the following well-known fact: for any closed convex polyhedral cone
$C$  $$
\sum_F (-1)^{\dim(F)}= 
\cases (-1)^{\dim(C)} &\text{if $C$ is a linear subspace of $X$},\\
0 &\text{otherwise}, \endcases   \tag{A.5}
$$
where the sum is taken over all faces $F$ of $C$ (to prove this fact reduce
to the case in which $C \ne \{0\}$ and $C$ contains no non-zero linear
subspace and note that in this case $C\setminus \{0\}$ is contractible,
hence has Euler characteristic 1; on the other hand this Euler
characteristic is equal to 
$$
-\sum_F (-1)^{\dim(F)} $$
where $F$ ranges through all faces of $C$ other than $\{0\}$).

We now define a homomorphism $f \mapsto f^*$ from $C(X)$ to itself by
putting
$$
   (\xi_C)^*:=(-1)^{\dim(C)}\xi_{\inC}. $$
Of course we must check that if this definition is applied to the left-hand
side of (A.3), we get 0.  We may as well assume that $\l$ takes both
strictly positive and strictly negative values on $C$ (otherwise the
relation (A.3) is itself trivial). Then
$$
\align
&\dim(C_+)=\dim(C_-)=\dim(C) \\
&\dim(C_0)=\dim(C)-1 \endalign $$
and $\inC$ is the disjoint union of $\inC_+$,$\inC_-$ and $\inC_0$,
which gives what we want.

Next we define a homomorphism $f \mapsto \hat f$ from $C(X)$ to $C(X^*)$ by
putting 
$$
(\xi_C)^\wedge := (-1)^{\dim(X)-\dim(C^*)}\xi_{\inC^*}. $$
Here $\inC^*$ denotes the relative interior of the dual cone 
$C^*$. Again we must check
that when we apply this definition to the left-hand side of (A.3), we get
0. We will see in a moment that the operation $*$ is an isomorphism;
therefore it is enough to consider instead the operation obtained by
following $\wedge$ by $*$ and multiplying by $(-1)^{\dim(X)}$; this
operation sends $\xi_C$ to $\xi_{C^*}$. We must show that
$$
\xi_{C^*} + \xi_{C^*_0} - \xi_{C^*_+} - \xi_{C^*_-} = 0. $$
This is an immediate consequence of the following well-known fact: if
$C_1,C_2$ are closed convex cones whose union is convex, then $C^*_1 \cup
C^*_2$ is also convex and
$$
   \align
(C_1 \cup C_2)^* &= C_1^* \cap C_2^* \\
(C_1 \cap C_2)^* &= C_1^* \cup C_2^*. \endalign $$

It is easy to see from the definitions of our two operations that
$f^{**}=f$ and $f^{\wedge * \wedge *}=f$ for all $f \in C(X)$ (for the first
equality use (A.4) and for the second use that $C^{**}=C$). It follows that
both of our operations are isomorphisms. These operations are best
understood by placing them in a more general context (see
\cite{KS},\cite{M}), in which $*$ comes from Verdier duality and $\wedge$
from the Fourier-Sato transformation, but this point of view is 
not needed for what follows.

It is interesting to calculate $(\xi_{\inC})^\wedge$. We claim that
$$
    (\xi_{\inC})^\wedge = (-1)^{\dim(C)}\xi_{-C^*}. \tag{A.6} $$
We will now give an elementary proof of this consequence of Lemma
3.7.10(ii) of \cite{KS}. Using (A.4), we see that (A.6) is equivalent to
the equality
$$
    \sum_{F \in \Cal F(C)} (-1)^{\dim(F)} \xi_F^{\wedge *}= (\xi_{-C^*})^*,
$$
which is in turn equivalent to 
$$
   \sum _{F \in \Cal F(C)} (-1)^{\dim(F)}
\xi_{F^*}=(-1)^{\dim(X)-\dim(C^*)} \xi_{-\inC^*}. \tag{A.7} $$
Let $\l \in X^*$. The value of the left-hand side of (A.7) at $\l$ is
$$
    \sum_{F \in \Cal F_\l} (-1)^{\dim(F)}, \tag{A.8} $$
where $\Cal F_\l$ denotes the set of faces $F$ of $C$ such that $\l$ is
non-negative on $F$. Let $C_+$ be the set
$$
  C_+:=\{ x \in C \,|\, \l(x) \ge 0 \}. $$
We claim that (A.8) is equal to
$$
   \sum_{F \in \Cal F(C_+)} (-1)^{\dim(F)}. \tag{A.9} $$
If $\l =0$, the claim is obvious, so we assume that $\l \ne 0$. Then $\Cal
F(C_+)$ contains $\Cal F_\l$. Elements of $\Cal F(C_+)$ that are not in
$\Cal F_\l$ arise in pairs $F \cap C_+$ and $F \cap \ker(\l)$, one pair for
each face $F$ of $C$ such that $\ker(\l)$ meets $\intF$ but does not
contain $\intF$. For such $F$
$$
    \dim(F\cap \ker(\l)) = \dim(F \cap C_+)-1; $$
and therefore
each such pair of faces contributes 0 to (A.9), and the claim follows.

It follows from (A.5) (applied to $C_+$) that (A.9) is equal to
$$ \cases   (-1)^{\dim(C_+)} &\text{if $C_+$ is a linear subspace,} \\
            0                &\text{otherwise}. \endcases $$ 
But $C_+$ is a linear subspace if and only if $\l$ is 
$\le 0$ on $C$ and $C \cap
\ker(\l)$ is a linear subspace. If $\l$ is $\le 0$ on $C$, then $C \cap
\ker(\l)$ is a face of $C$, and therefore $C \cap \ker(\l)$ is a linear
subspace if and only if it is equal to the unique minimal face of $C$. Thus
$C_+$ is a linear subspace if and only if $\l \in -\inC^*$, in which case
$C_+=C \cap \ker(\l)$ is the minimal face of $C$ and thus has dimension
equal to 
$$\dim(X)-\dim(C^*).$$
This completes the proof of the equality (A.7).

We say that an integer-valued 
function on $X \times X^*$ is \it biconic \rm if it is
a finite $\Z$-linear combination of characteristic functions $\xi_{C
\times D}$, where $C$ (respectively, $D$) is a closed convex polyhedral
cone in $X$ (respectively, $X^*$). Of course the functions $\psi_C$
considered earlier are biconic functions on $X \times X^*$. In fact we are
now going to generalize $\psi_C$ by associating to any conic function $f$
on $X$ a biconic function $\psi_f$ on $X \times X^*$. If $f$ is the
characteristic function $\xi_C$ of a closed convex polyhedral cone $C$ we
define $\psi_f$ to be the function $\psi_C$ defined earlier. Since any
conic function $f$ can be written as a $\Z$-linear combination 
$$
f= \sum_C n_C \xi_C $$
for a finite set of closed convex polyhedral cones $C$ and integers $n_C$,
we are then forced to define $\psi_f$ by
$$
   \psi_f:=\sum_C n_C \psi_C. $$

Of course we must show that $\psi_f$ depends only on $f$ and not on the way in
which we write it as a $\Z$-linear combination of characteristic functions
of cones. This can be done directly, but instead we will use the two
operations introduced earlier to give a shorter proof. 
For any biconic function $\psi$ and any $\l \in X^*$ the function
$\psi(\cdot,\l)$ on $X$ is conic, and therefore we can apply the operation
$*$ to $\psi$ in the variable $x \in X$ (holding $\l$ fixed); this
operation sends $\xi_{C \times D}$ to $(-1)^{\dim(C)}\xi_{\inC \times
D}$, which is again biconic. It follows that applying $*$ to $\psi$ in the
variable $x \in X$ yields a biconic function.
Similarly, applying $*$ to $\psi$ in the variable $\l \in X^*$ yields a
biconic function. 
Any biconic function is conic on $X \times X^*$ and so we may apply the
operation $\wedge$ on $X \times X^*$ to any biconic function $\psi$ on $X
\times X^*$. Since the dual of $X \times X^*$ is $X^* \times X$, which we
may identify with $X \times X^*$ by switching the order of the two factors,
we may regard $\wedge$ as an operation taking biconic functions on $X
\times X^*$ to biconic functions on $X \times X^*$ (it is clear that $\hat
\psi$ is biconic, not just conic). 

For a biconic function $\psi$ on $X \times X^*$ we denote by $\psi'$ the
biconic function on $X \times X^*$ obtained by applying $*$ in the variable
$\l \in X^*$ to the biconic function $\hat \psi$. We will use the operation
$\psi \mapsto \psi'$ to show that $f \mapsto \psi_f$ is well-defined. It is
enough to show that
$$
   \sum_C n_C (\psi_C)' $$
depends only on $f$. We need to find a simple expression for $(\psi_C)'$;
by definition it is 
$$
\sum_F \xi_{\intF \times F^\perp}. $$
But for any $x \in \intF$ we have
$$
F^\perp=C^* \cap \{ \l \in X^* \, | \, \l(x)=0 \}. $$
Therefore $(\psi_C)'(x,\l)$ is 0 unless $\l(x)=0$, in which case it is
equal to
$$
\sum_F \xi_{\intF \times C^*}(x,\l)=\xi_{C \times C^*}(x,\l). $$

Our problem has been reduced to the following: for $(x,\l) \in X \times
X^*$ such that $\l(x)=0$ we must show that
$$
\sum_C n_C \xi_{C \times C^*} (x,\l) \tag{A.10} $$
depends only on $f$. In fact we will show more: for such $(x,\l)$ the
expression (A.10) is equal to
$$
[f^*\xi_{X_+}]^*(x) \tag{A.11} $$
where $X_+$ is the set
$$
X_+=\{x \in X \, | \, \l(x) \ge 0 \}. $$
Since this statement is linear in $f$, we may assume without loss of
generality that $f=\xi_C$. Then (A.10) is equal to $\xi_{C \times
C^*}(x,\l)$, which is in turn equal to
$$
\cases \xi_C(x) &\text{if $C \subset X_+$,}\\ 0 &\text{otherwise.}
\endcases 
$$ 
By definition $f^*$ equals $(-1)^{\dim(C)}\xi_{\inC}$ and therefore the
product $f^*\xi_{X_+}$ equals $$(-1)^{\dim(C)}\xi_{\inC \cap X_+}.$$

There are three cases. If $C$ is contained in $X_+$, 
then $\inC \cap X_+=\inC$ and
$[f^*\xi_{X_+}]^*=\xi_C$, so that (A.11) agrees with (A.10) in this case.
If $C$ is not contained in $X_+$ but is contained in $X_-$, where 
$$
   X_-=\{x \in X \, | \, \l(x) \le 0 \}, $$
then $\inC \cap X_+$ is empty and therefore (A.11) is 0, which again agrees
with (A.10). If $C$ is contained in neither $X_+$ nor $X_-$, then we are in
the situation considered during the proof that $*$ is well-defined, and
$\inC$ is the disjoint union of $\inC_+$, $\inC_-$
and $\inC_0$, where $C_+=C \cap X_+$, $C_-=C\cap X_-$ and 
$C_0=C \cap \ker(\l)$. Therefore $\inC \cap X_+$ is
the disjoint union of $\inC_+$ and $\inC_0$, and it follows that
$$ \align
[f^*\xi_{X_+}]^*&=(-1)^{\dim(C)}[\xi_{\inC_+}+\xi_{\inC_0}]^* \\
&=\xi_{C_+}-\xi_{C_0}. \endalign $$
For $x$ such that $\l(x)=0$, we have 
$$
(\xi_{C_+}-\xi_{C_0})(x)=\xi_C(x)-\xi_C(x)=0, $$
which shows that (A.11) again agrees with (A.10). This concludes the proof
that $\psi_f$ is well-defined.

The following proposition is an easy consequence of the fact that $\psi_f$
is well-defined.
\proclaim{Proposition A.4} Let $C$ be a closed convex polyhedral cone and
suppose that $\inC$ is the disjoint union of the relative interiors
$\inC_1,\dots,\inC_r$ of $r$ closed convex polyhedral cones
$C_1,\dots,C_r$. Then
$$
\psi_C=\sum_{i=1}^r (-1)^{\dim(C)-\dim(C_i)}\psi_{C_i}. $$
\endproclaim

Indeed our hypothesis is that 
$$
\xi_{\inC}=\sum_{i=1}^r \xi_{\inC_i}. $$
Applying the operation $*$, we find that
$$
(-1)^{\dim(C)}\xi_C=\sum_{i=1}^r (-1)^{\dim(C_i)} \xi_{C_i}. $$
Applying the map $f \mapsto \psi_f$, we then get the desired
equality.
\proclaim{Proposition A.5} For any conic function $f$ on $X$ the quantity
$\psi_f(x,\l)$ vanishes unless $\l(x)\le 0$. \endproclaim

The conic function $f$ can be written as a $\Z$-linear combination of
characteristic functions of closed convex polyhedral cones $C$ that are
simplicial as cones in their linear spans. Therefore it is enough to prove
the proposition in the case that $f=\xi_C$ for such a cone $C$. By (A.2)
we may assume that $C$ spans $X$ and hence that $C$ is simplicial in $X$.
Then the proposition follows easily from Lemma A.1. 

The function $\psi_C$ is equal to $\psi_f$ for the conic function
$f=\xi_C$. We can also use the \it open \rm cone $\inC$ to 
obtain an equally useful variant $\varphi_{\inC}$ of $\psi_C$
by putting
$$
    \varphi_{\inC}:=\psi_g,$$
where $g$ is the conic function $\xi_{\inC}$. 

\proclaim{Lemma A.6} There is an equality 
$$
  \varphi_{\inC}(x,\l)=\sum_{F \in \Cal F} (-1)^{\dim(F)}
\xi_{\inC+\spn(F)}(x) \xi_{F^*}(\l). $$  \endproclaim

The cone $\inC+\spn(F)$ is the relative interior of the cone $C+\spn(F)$,
whose faces are in one-to-one correspondence with the faces $G$ of $C$
containing $F$, via the map $G \mapsto G+\spn(F)$. Applying (A.4) to
$C+\spn(F)$, we see that the right-hand side of the equality we are trying
to prove is equal to 
$$
   \sum_F (-1)^{\dim(F)} \sum_{\{G \in \Cal F | G \supset F \}}
(-1)^{\dim(C)-\dim(G)} \xi_{G+\spn(F)} (x) \xi_{F^*}(\l). \tag{A.12} $$
Applying (A.4) to $C$, we see that
$$
   \varphi_{\inC}=\sum _{G \in \Cal F} (-1)^{\dim(C)-\dim(G)} \psi_G, $$
and by writing out the definition of $\psi_G$, we see that
$\varphi_{\inC}(x,\l)$ is equal to the expression (A.12). This proves the
lemma. 

Note that the expression for $\varphi_{\inC}(x,\l)$ given by Lemma A.6 
is almost the same as
the expression for $\psi_C(x,\l)$ given by its definition; the only
difference is that the factor
$$
  \xi_{C+\spn(F)}(x) $$
appearing in the definition of $\psi_C(x,\l)$ is replaced by
$$
   \xi_{\inC+\spn(F)}(x) $$
in the expression for $\varphi_{\inC}(x,\l)$. Since $\inC+\spn(F)$ is the
relative interior of $C+\spn(F)$, we conclude  that for any conic function
$f$ on $X$, the biconic function $\psi_{f^*}$ associated to $f^*$ is
obtained from $\psi_f$ by applying the operation $*$ to $\psi_f$ in the
first variable. Moreover we conclude that
$$
    \varphi_{\inC}(x,\l)=\psi_C(x,\l) \text{ if $x$ is
$C$-regular.} \tag{A.13} $$

It is no surprise that $\varphi_{\inC}$ behaves just about the same way as
$\psi_C$. For example $\varphi_{\inC}$ satisfies the obvious analog of
(A.2) (replace $C,\tilde C$ in (A.2) by their relative interiors). Now
suppose that $C$ is a simplicial cone in $X$ and use the same notation as
in Lemma A.1. We need one more bit of notation: put 
$$ 
\inI_x = \{ i \in I \, | \, x_i > 0 \}. $$
Then it is easy to see that $\varphi_{\inC}$ satisfies the following analog of
Lemma A.1.
\proclaim{Lemma A.7} The number $\varphi_{\inC}(x,\l)$ is $0$ unless the
subsets $\inI_x$ and $I_\l$ of $I$ are complementary, in which case 
$$
   \varphi_{\inC}(x,\l)=(-1)^{|I_\l|}, $$
where $|I_\l|$ denotes the cardinality of $I_\l$. \endproclaim\heading B. Combinatorial lemma of Langlands \endheading
In this appendix we generalize a combinatorial lemma of Langlands
\cite{A1,Lemma 6.3}. Let $X$ be a finite dimensional real vector space, and
let $(\cdot,\cdot)$ be a positive definite symmetric bilinear form on $X$,
with associated metric
$$
   d(x,y)= (x-y,x-y)^{1/2}. $$

Let $C$ be a closed convex polyhedral cone in $X$. Let $x \in X$. Since $C$
is closed there exists a point $x_0 \in C$ that is closest to $x$, and
since $C$ is convex, the point $x_0 \in C$ is unique. Again by convexity
this closest point can be characterized as the unique point $x_0 \in C$
such that $$
  d(x,x_0) \le d(x,x_0+r(y-x_0)) $$
for all $y \in C$ and all real numbers $r \in [0,1]$.
For fixed $y \in C$ the truth of the inequality above for all $r \in [0,1]$
is equivalent to the inequality 
$$
   (x-x_0,y-x_0)\le 0. $$
Thus, since $C$ is a cone, $x_0 \in C$ is characterized by the property
that
$$
    (x-x_0,z) \le 0 $$ 
for all $z \in Z$, where $Z$ is the set of all elements in $X$ of the form
$y-rx_0$ for some $y \in C$ and some positive real number $r$. 

Let $F$ be the unique face of $C$ such that $x_0$ lies in the relative
interior $\intF$ of $F$. We claim that $Z$ is equal to $C+\spn(F)$.
Clearly $Z$ is contained in $C+\spn(F)$. Moreover, in order to prove the
reverse inclusion it is enough to show that $-F$ is contained in $Z$. Let
$x_1 \in F$. Since $x_0 \in \intF$, there exist $x_2 \in F$ and a positive
real number $r$ such that
$$
    x_1-x_0 = -r(x_2-x_0), $$
which shows that $-x_1 \in Z$, as desired.

We use the inner product $(\cdot,\cdot)$ to identify $X$ with its dual. In
particular we now view the dual cone $C^*$ as subset of $X$ itself: 
$$
  C^*=\{x \in X \, | \, (x,y) \ge 0 \text{ for all $ y \in C$} \}. $$
As in Appendix A we denote by $F^\perp$ the face $C^* \cap \spn(F)^\perp$ of
$C^*$ determined by the face $F$ of $C$. Of course $F^\perp$ is equal to
$$
   (C+\spn(F))^*. $$
We conclude that $x_0 \in \intF$ is the point of $C$ closest to $x$ if and
only if 
$$
   x-x_0 \in -F^\perp, $$
in which case $x_0$ is the orthogonal projection of $x$ on $\spn(F)$ and
$x-x_0$ is the orthogonal projection of $x$ on $\spn(F)^\perp$. In
particular the set of all points $x \in X$ such that the point $x_0$ in $C$
closest to $x$ lies in $\intF$ is equal to
$$
   \intF \oplus (-F^\perp). $$
We have written $\oplus$ rather than $+$ in order to emphasize that the
cones $\intF$ and $-F^\perp$ lie in the complementary subspaces $\spn(F)$ 
and $\spn(F)^\perp$ respectively. 
Let $\xi_{\intF \oplus (-F^\perp)}$ denote the characteristic
function of the subset $\intF \oplus (-F^\perp)$ of $X$. We conclude that
$$
   1 = \sum _{F \in \Cal F} \xi _{\intF \oplus (-F^\perp)}(x) \tag{B.1} $$
for all $x \in X$, where $\Cal F$ denotes the set of faces of $C$. 

The equality (B.1) is a (generalization of) a simple
 special case of Langlands's
combinatorial lemma. We will now use this special case to prove the general
case. 

We continue to identify $X$ with its dual, so that the function $\psi_C$ on
$X \times X^*$ (see Appendix A) becomes a function on $X \times X$. For any
subspace $Y$ of $X$ we denote by $p_Y$ the orthogonal projection map from
$X$ onto $Y$. Now we can state our generalization of the combinatorial lemma
of Langlands (take $C$ to be simplicial and use Lemma A.1 to recover the
usual form of the lemma).
\proclaim{Lemma B.1} For $x,y \in X$ there is an equality 
$$
   \sum_{F \in \Cal F} \psi_{C+\spn(F)} (-y,-p_{\spn(F)^\perp}(x)) \cdot
\xi_{\intF}(p_{\spn(F)}(x)) = (-1)^{\dim(C)}\xi_{\inC}(y). $$
\endproclaim
The faces of $C+\spn(F)$ are precisely the cones $G+\spn(F)$, where $G$
ranges through the set of faces of $C$ containing $F$; note that for such
$G$
$$
   \align
&\spn(G+\spn(F))=\spn(G) \\
&\dim(G+\spn(F))=\dim(G) \\
&(C+\spn(F))+\spn(G+\spn(F))=C+\spn(G). \endalign $$
Therefore by the definition of $\psi_{C+\spn(F)}$ the left-hand side of
the equality in the lemma is equal to
$$
 \sum_{F \in \Cal F} \sum_{G \in \Cal F: G \supset F} (-1)^{\dim(G)}
\xi_{C+\spn(G)} (-y) \cdot \xi_{(G+\spn(F))^*}(-p_{\spn(F)^\perp}(x))
\cdot \xi_{\intF}(p_{\spn(F)}(x)), $$
which by interchanging the order of summation we rewrite as
$$
  \sum_{G \in \Cal F} (-1)^{\dim(G)} \xi_{C+\spn(G)}(-y) \cdot \sum_{F \in
\Cal F(G)} \xi_{(G+\spn(F))^*} (-p_{\spn(F)^\perp}(x)) \cdot
\xi_{\intF}(p_{\spn(F)} (x)), $$
where $\Cal F(G)$ denotes the set of faces of $G$. The second sum is 1 by
(B.1) (applied to $G$). Therefore the double sum reduces to
$$
   \sum_{G \in \Cal F} (-1)^{\dim(G)} \xi_{C+\spn(G)}(-y). $$
Applying equality (A.7) to $C^*$, we see that this last expression is equal
to
$$
  (-1)^{\dim(C)}\xi_{\inC}(y). $$
This completes the proof of the lemma.

There are two special cases of Lemma B.1 that are worth noting. First, if
$y \in \inC$, it is easy to see that Lemma B.1 reduces to the equality
(B.1). Second, if $y \in -C$, then Lemma B.1 reduces to the following
result.
\proclaim{Corollary B.2} There is an equality
$$
   \sum_{F \in \Cal F} (-1)^{\dim(F^\perp)}\xi_{(F^\perp)^\circ \oplus
\intF} =\cases (-1)^{\dim(X)-\dim(C)} &\text{if $C$ is a linear subspace,}
\\ 0 &\text{otherwise.} \endcases $$
\endproclaim
Indeed, if $x \in C$, then $\xi_{(F^\perp)^*}(x)=1$ for all $F \in \Cal F$,
and therefore
$$ \align
 \psi_C(x,\l)&=\sum_{F \in \Cal F} (-1)^{\dim(F)} \xi_{F^*}(\l) \\
&=(-1)^{\dim(X)-\dim(C^*)} \xi_{\inC^*}(-\l) \endalign $$
(we used the equality (A.7)). Applying this to $C+\spn(F)$ instead of $C$,
we see that if $y \in -C$, then
$$\psi_{C+\spn(F)}(-y,\l)=(-1)^{\dim(X)-\dim(F^\perp)}
\xi_{(F^\perp)^\circ}(-\l). $$
Moreover, if $y \in -C$, then
$$
    \xi_{\inC}(y)=\cases 1 &\text{if $C$ is a linear subspace,} \\
0 &\text{otherwise.} \endcases $$
Therefore Lemma B.1 does reduce to Corollary B.2 when $y \in -C$. 
\proclaim{Corollary B.3} The sum
$$
    \sum _{F \in \Cal F} (-1)^{\dim(F)} \xi_{((F^\perp)^*)^\circ \oplus
(F^*)^\circ} (x) $$ 
is $0$ unless $C$ is a linear subspace and $x=0$. Here the dual cone
$(F^\perp)^*$ is taken inside the subspace 
$$
   \spn(F)^\perp=\spn(F^\perp) $$
and the dual cone $F^*$ is taken inside the subspace $\spn(F)$.
\endproclaim
This corollary can be derived from the previous one by applying the
operation $*$ and then the operation $\wedge$. Note that the formula 
$$
    \sum_{\{P_3:P_3 \supset P_1\}}
(-1)^{\dim(A_1/A_3)}\tau^3_1(H)\hat\tau_3(H)=0 $$
on p\. 940 of \cite{A1} is the special case of this corollary in which the
cone $C$ is the closed chamber in $\fA_{M_1}/\fA_G$ determined by $P_1$
(still using Arthur's notation). 
\end